\documentclass{modified}
\usepackage{amsmath,amsfonts,amsthm,xypic,graphics}

\newcommand{\Ker}{\operatorname{Ker}}

{}
\newtheorem{theorem}{Theorem}[chapter]
\newtheorem{corollary}[theorem]{Corollary}
\newtheorem{prop}[theorem]{Proposition}
\newtheorem{lemma}[theorem]{Lemma}
\newtheorem{defn}[theorem]{Definition}
\newtheorem{remark}[theorem]{Remark}
\newcommand{\nc}{\newcommand}
\nc{\FH}{\mathcal H} \nc{\QQ}{\mathbb Q} \nc{\CC}{\mathbb C}
\nc{\JJ}{\mathcal J} \nc{\KK}{\mathbb K} \nc{\RR}{\mathbb R}
\nc{\LL}{\mathcal L} \nc{\Ll}{\ell} \nc{\NN}{\mathbb N}
\nc{\ZZ}{\mathbb Z} \nc {\HH}{\mathbb H} \nc {\OO}{\mathcal{O}}
\nc{\lra}{\longrightarrow} \nc{\bdot}{\bullet} \nc{\w}{\omega}
\nc{\Jac}{\operatorname{\mathrm{Jac}}}
\nc{\Nij}{\operatorname{\mathrm{Nij}}} \nc{\TT}{\mathcal{T}}
\nc{\MM}{\mathcal{M}}\nc{\Har}{\mathcal{H}} \nc{\zed}{\mathcal{Z}}
\nc{\image}{\mathrm{Im}\ }

\newcommand{\Ann}{\operatorname{\mathrm{Ann}}}
\newcommand{\rk}{\operatorname{\mathrm{rk}}}
\newcommand{\Tr}{\operatorname{\mathrm{Tr}}}
\newcommand{\End}{\operatorname{\mathrm{End}}}
\newcommand{\so}{\mathfrak{so}}
\newcommand{\IP}[1]{\langle #1\rangle}

\newcommand{\del}{\partial}
\newcommand{\delbar}{\overline\partial}
\newcommand{\BB}{\mathcal B}

\theoremstyle{definition}
\newtheorem{example}[theorem]{Example}

\nc {\PP}{\mathbb P} \nc{\FF}{\mathcal{F}}

\nc{\EE}{\mathcal{E}} \nc{\WW}{\mathcal{W}}
\nc{\nabt}{\widetilde{\nabla}}
\nc{\eps}{\varepsilon}\nc{\IPS}[1]{#1} \nc{\DD}{\mathcal{D}}

\title{Generalized complex geometry}
\author{Marco Gualtieri}
\college{St John's College}
\degree{Doctor of Philosophy} \degreedate{November 2003}
\begin{document}
\baselineskip=14pt plus1pt
\setcounter{secnumdepth}{2} \setcounter{tocdepth}{3}

\maketitle

\clearpage
\thispagestyle{empty}
\cleardoublepage
\begin{acknowledgements}
During the four years I spent in Oxford, I always looked forward
to my weekly meetings with Nigel Hitchin, my supervisor.  In
coping with my interminable lists of questions in the early days
as well as my interminable lists of results in the later days,
Nigel was always remarkably patient and amazingly insightful, for
which I am very grateful.

Many thanks also to my friends and colleagues at the Maths
Institute, at St. John's, and at New College, for creating such an
interesting atmosphere in which to live and work.  Good luck, to
those who remain!

Special thanks to my parents Antoinette and Renzo for their love
and support.

Finally, I acknowledge the generous funding of the Rhodes Trust
and of the National Science and Engineering Research Council of
Canada which enabled me to complete this research.
\end{acknowledgements}

\clearpage
\thispagestyle{empty}
\cleardoublepage
\begin{center}
\thispagestyle{empty}
{\LARGE\bf Generalized complex geometry}\\
\vspace{20pt}
{\large Marco Gualtieri}\\
{\large Oxford University D.Phil. Thesis}\\
\vspace{20pt}

{\large \bf Abstract}\\\end{center} \vspace{10pt}

Generalized complex geometry is a new kind of geometrical
structure which contains complex and symplectic geometry as its
extremal special cases. In this thesis, we explore novel phenomena
exhibited by this geometry, such as the natural action of a
$B$-field.  We provide many examples of generalized complex
structures, including some on manifolds which admit no known
complex or symplectic structure. We prove a generalized Darboux
theorem which yields a local normal form for the geometry. We show
that there is a well-behaved elliptic deformation theory and
establish the existence of a Kuranishi-type moduli space.

We then introduce a Riemannian metric and define the concept of a
generalized K\"ahler manifold.  We prove that generalized K\"ahler
geometry is equivalent to a certain bi-Hermitian geometry first
discovered by physicists in the context of supersymmetric
sigma-models.  We then use this theorem together with our
deformation result to solve an outstanding problem in
4-dimensional bi-Hermitian geometry: we prove that there exists a
Riemannian metric on $\CC P^2$ which admits exactly two distinct
orthogonal complex structures with equal orientation.

In addition, we introduce the concept of a generalized complex
submanifold, and show that these sub-objects correspond precisely
with the predictions of physicists concerning D-branes in the
special cases of complex and symplectic manifolds.

\clearpage
\thispagestyle{empty}
\cleardoublepage

\begin{romanpages}
\tableofcontents
\end{romanpages}

\clearpage
\thispagestyle{empty}
\cleardoublepage
\chapter{Introduction}

A central theme of this thesis is that classical geometrical
structures which appear, at first glance, to be completely
different in nature, may actually be special cases of a more
general unifying structure. Of course, there is wide scope for
such generalization, as we may consider structures defined by
sections of any number of natural bundles present on manifolds.
What must direct us in deciding which tensor structures to study
is the presence of natural integrability conditions.

Good examples of such conditions include the closure of a
symplectic form, the Einstein or special holonomy constraint on a
Riemannian metric, the vanishing of the Nijenhuis tensor of a
complex structure, and the Jacobi identity for a Poisson bivector,
among many others.

This thesis describes a way of unifying complex and symplectic
geometry by taking seriously the idea that both of these
structures should be thought of, not as linear operations on the
tangent bundle of a manifold, but actually on the sum of the
tangent and cotangent bundles, $T\oplus T^*$.  Since the smooth
sections of $T\oplus T^*$ have a natural bracket operation called
the \emph{Courant bracket}, there are canonical integrability
conditions for certain linear structures on $T\oplus T^*$. Indeed,
any complex or symplectic structure determines a maximal isotropic
sub-bundle of $(T\oplus T^*)\otimes\CC$;  the requirement that
this sub-bundle be Courant involutive actually specializes to the
usual integrability conditions for these two structures.  This was
one of the observations which led Hitchin~\cite{Hitchin} to define
a generalized complex structure as an almost complex structure
$\JJ$ on $T\oplus T^*$ whose $+i$-eigenbundle $L$ is Courant
involutive. This new geometrical structure is, in a sense, the
complex analogue of a Dirac structure, a concept introduced by
Courant and Weinstein~\cite{Courant},\cite{CourWein} to unify
Poisson and symplectic geometry.

We begin, in Chapter~\ref{alg}, with a study of the natural
split-signature orthogonal structure which exists on the real
vector bundle $T\oplus T^*$.  A spin bundle for this orthogonal
bundle is shown always to exist and to be isomorphic to
$\wedge^\bullet T^*$, the bundle of differential forms. The
correspondence between maximal isotropic subspaces and pure
spinors then leads to the fact that a generalized complex
structure is determined by a canonical line sub-bundle of the
complex differential forms.  In the case of a complex manifold,
this line bundle is the usual canonical line bundle.  In the
symplectic case, however, this line bundle is generated by
$e^{i\omega}$, where $\omega$ is the symplectic form.

We proceed, in Chapter~\ref{cb}, to describe and study the Courant
bracket, which, while it is not a Lie bracket, does restrict, on
involutive maximal isotropic sub-bundles, to be a Lie bracket, and
thus endows the bundle $L$ with the structure of a Lie algebroid.
$L$ acquires not only a Lie bracket on its sections, but also an
exterior derivative operator $d_L:C^\infty(\wedge^k
L^*)\rightarrow C^\infty(\wedge^{k+1} L^*)$. The Courant
integrability of $L$ may also be phrased in terms of a condition
on the differential forms defining it, and we determine this
condition.  But perhaps the most important feature of the Courant
bracket is that, unlike the Lie bracket of vector fields, it
admits more symmetries than just diffeomorphisms.  The extra
symmetries are called $B$-field transformations and are generated
by closed 2-forms $B\in\Omega^2(M)$.  This means that given any
structure naturally defined in terms of the Courant bracket, like
generalized complex structures, a $B$-field transform produces
another one. This action of the 2-forms agrees precisely with the
2-form gauge freedom studied by physicists in the context of sigma
models.  We conclude the chapter with an investigation of the fact
that the Courant bracket itself may be deformed by a real closed
3-form $H$, and we describe what this means in the language of
gerbes.

The first two chapters have been organized to contain all the
necessary algebraic and differential-geometric machinery for the
rest of the thesis.  In Chapter~\ref{GC}, we come to the subject
at hand: generalized complex structures themselves.  We show that
topologically, the obstruction to the existence of a generalized
complex structure is the same as that for an almost complex
structure or a nondegenerate 2-form.  We then describe the
algebraic conditions on differential forms which makes them
generators for generalized complex structures.  This then allows
us to produce exotic examples of generalized complex structures;
indeed we exhibit examples on manifolds which admit no known
complex or symplectic structure.  We are even able to give an
example of a family of generalized complex structures which
interpolates between a complex and a symplectic structure, thus
connecting the moduli spaces of these two structures.

Still in Chapter~\ref{GC}, we prove a local structure theorem for
generalized complex manifolds, analogous to the Darboux theorem in
symplectic geometry and the Newlander-Nirenberg theorem in complex
geometry. We show that at each point of a $2n$-dimensional
generalized complex manifold, the structure is characterized
algebraically by an integer $k$, called the \emph{type}, which may
vary along the manifold and take values anywhere from $k=0$ to
$k=n$. It is lower semi-continuous, and a point where it is
locally constant is called a \emph{regular} point of the
generalized complex manifold.  Our local structure theorem states
that near any regular point of type $k$, the generalized complex
manifold is equivalent to a product of a complex space of
dimension $k$ with a symplectic space.  It is crucial to note,
however, that this equivalence is obtained not only by using
diffeomorphisms but also $B$-field transformations.  This is
consistent with the fact that, as symmetries of the Courant
bracket, $B$-field transformations should be considered on a par
with diffeomorphisms.  We end that chapter by defining
\emph{twisted} generalized complex structures, thus enabling us to
deform all our work by a real closed 3-form $H$.

Chapter~\ref{defm} contains our main analytical result: the
development of a Kuranishi deformation space for compact
generalized complex manifolds.  The deformation theory is governed
by the differential complex $(C^\infty(\wedge^k L^*),d_L)$
mentioned above, which we show is elliptic, and therefore has
finite-dimensional cohomology groups $H^k_L(M)$ on a compact
manifold $M$. In particular, integrable deformations correspond to
sections $\eps\in C^\infty(\wedge^2 L^*)$ satisfying the
Maurer-Cartan equation
\begin{equation*}
d_L\eps+\tfrac{1}{2}[\eps,\eps]=0.
\end{equation*}
There is an analytic obstruction map $\Phi:H^2_L(M)\rightarrow
H^3_L(M)$, and if this vanishes then there is a locally complete
family of deformations parametrized by an open set in $H^2_L(M)$.
In the case that we are deforming a complex structure, this
cohomology group turns out to be
\begin{equation*}
H^0(M,\wedge^2T)\oplus H^1(M,T)\oplus H^2(M,\mathcal{O}).
\end{equation*}
This is familiar as the ``extended deformation space'' of
Barannikov and Kontsevich~\cite{KontsevichBarannikov}, for which a
geometrical interpretation has been sought for some time.  Here it
appears naturally as the deformation space of a complex structure
in the generalized setting.

In chapter \ref{GeneralizedKahler}, we introduce a Riemannian
metric $G$ on $T\oplus T^*$ in such a way that it is compatible
with two commuting generalized complex structures $\JJ_1,\JJ_2$.
This enriched geometry generalizes classical K\"ahler geometry,
where the compatible complex and symplectic structures play the
role of $\JJ_1$ and $\JJ_2$.  In studying this geometry we
discover the result that it is equivalent to a geometry first
discovered by Gates, Hull and Ro\v{c}ek~\cite{Rocek} in their
study of supersymmetric sigma-models.  Our proof of this
equivalence occurs in section~\ref{relat}.  The latter geometry is
a bi-Hermitian geometry with an extra condition relating the pair
of Hermitian 2-forms.  The general subject of bi-Hermitian
geometry has been studied in depth by Apostolov, Gauduchon, and
Grantcharov~\cite{Gauduchon} in the four-dimensional case, and in
their paper they state that the main unsolved problem in the field
is to determine whether or not there exist bi-Hermitian structures
on complex surfaces which admit no hyperhermitian structure, for
example, $\CC P^2$.  Using our equivalence theorem together with
the deformation theorem of chapter~\ref{deform}, we are able to
prove that a bi-hermitian structure on $\CC P^2$ exists.  That is,
we show that there exists a Riemannian metric on $\CC P^2$
admitting exactly two distinct orthogonal complex structures.

Chapter $\ref{GeneralizedKahler}$ ends with the definition of
$\emph{twisted}$ generalized K\"ahler geometry, and we provide a
family of examples of twisted generalized K\"ahler manifolds: the
compact even-dimensional semi-simple Lie groups. We then define a
further restriction on K\"ahler geometry analogous to the
Calabi-Yau constraint: we define generalized Calabi-Yau metric
geometry.

Chapter~\ref{subm} introduces the concept of generalized complex
submanifold. The requirement that the property of being such a
submanifold should be invariant under $B$-field transformations
renders the definition highly non-trivial; it turns out that the
geometry of $T\oplus T^*$ demands that the submanifold carry a
2-form $F$; indeed if $(N,H)$ is a manifold together with a closed
3-form $H$ (the `twist'), we say that the pair $(M,F)$ of a
submanifold $M\subset N$ and a 2-form $F\in \Omega^2(M)$ is a
generalized submanifold of $(N,H)$ if $H|_M=dF$. We explore the
gerbe interpretation of this statement, and see how, in special
cases, $F$ should be interpreted as the curvature of a unitary
line bundle on $M$. This situation would sound very familiar to
physicists studying D-branes, which in certain cases are
represented by submanifolds equipped with line bundles.  While in
the case of complex manifolds, a generalized complex submanifold
is simply a complex submanifold equipped with a unitary
holomorphic line bundle (or more generally a closed $(1,1)$-form),
in the symplectic case we obtain not only Lagrangian submanifolds
equipped with flat line bundles, but also a special class of
co-isotropic submanifolds, equipped with a line bundle with
constrained curvature.  These are precisely the co-isotropic
A-branes recently discovered by Kapustin and
Orlov~\cite{Kapustin}. Finally, to conclude this chapter, we
indicate how, in the case of generalized Calabi-Yau metric
geometry, one defines the analog of calibrated special Lagrangian
submanifolds.

For many years physicists and mathematicians have studied the
mysterious links between complex and symplectic geometry predicted
by mirror symmetry.  Therefore, the explicit unification of these
two structures would seem an important step in understanding how
they are connected.  Indeed, this thesis is full of formulae and
concepts which appear in some form or other in the realm of mirror
symmetry.  In the final, very speculative, chapter of the thesis
we propose a vague picture of how mirror symmetry might be phrased
in the language of generalized complex geometry.

\clearpage
\thispagestyle{empty}
\cleardoublepage
\chapter{Linear algebra of $V\oplus V^*$}\label{alg}

Let $V$ be a real vector space of dimension $m$, and let $V^*$ be
its dual space.  Then $V\oplus V^*$ is endowed with the following
natural symmetric and skew-symmetric bilinear forms:
\begin{align*}
\IP{X+\xi,Y+\eta}_+&=\frac{1}{2}(\xi(Y)+\eta(X))\\
\IP{X+\xi,Y+\eta}_-&=\frac{1}{2}(\xi(Y)-\eta(X)),
\end{align*}
where $X,Y\in V$ and $\xi,\eta\in V^*$.  Both bilinear forms are
nondegenerate, and it is the symmetric one which is ubiquitous in
this thesis; for this reason we usually denote it by $\IP{,}$ and
refer to it as `the inner product'.  This symmetric inner product
has signature $(m,m)$ and therefore defines the non-compact
orthogonal group $O(V\oplus V^*)\cong O(m,m)$.  In addition to
these bilinear forms, $V\oplus V^*$ has a canonical orientation,
as follows.  The highest exterior power can be decomposed as
\begin{equation*}
\wedge^{2m}(V\oplus V^*)=\wedge^{m}V\otimes\wedge^{m}V^*,
\end{equation*}
and there is a natural pairing between $\wedge^kV^*$ and
$\wedge^kV$ given by
\begin{equation*}
(v^*,u)=\det(v_i^*(u_j)),
\end{equation*}
where $v^*=v^*_1\wedge\cdots\wedge v^*_k\in\wedge^kV^*$ and
$u=u_1\wedge\cdots\wedge u_k\in\wedge^kV$. Therefore we have a
natural identification $\wedge^{2m}(V\oplus V^*)=\RR$, and the
number $1\in\RR$ defines a canonical orientation on $V\oplus V^*$.
The Lie group preserving the symmetric bilinear form together with
the canonical orientation is of course the special orthogonal
group $SO(V\oplus V^*)\cong SO(m,m)$.

In this section we study the behaviour of maximal isotropic
subspaces of $V\oplus V^*$ and their description using pure
spinors. Further details can be found in the main reference for
this classical material, Chevalley's monograph~\cite{Chevalley},
which contains a chapter dealing exclusively with bilinear forms
of signature $(m,m)$. Also, we complexify the situation and
consider the \emph{real index} of a complex maximal isotropic
subspace.

\section{Symmetries of $V\oplus V^*$}\label{symm}

For all that follows, it is crucial to understand certain special
symmetries of $V\oplus V^*$.  The Lie algebra of the special
orthogonal group $SO(V\oplus V^*)$ is defined as usual:
\begin{align*}
\mathfrak{so}(V\oplus V^*)=\left\{T\ \big|\ \IP{Tx,y}+\IP{x,Ty}=0\
\ \forall\  x,y\in V\oplus V^*\right\}
\end{align*}
Using the splitting $V\oplus V^*$ we can decompose as follows:
\begin{equation}\label{comp}
T=\left(\begin{matrix}A&\beta\\B&-A^*\end{matrix}\right),
\end{equation}
where $A\in \End(V)$, $B:V\lra V^*$, $\beta:V^*\lra V$, and where
$B$ and $\beta$ are skew, i.e. $B^*=-B$ and $\beta^*=-\beta$.
Therefore we may view $B$ as a 2-form in $\wedge^2V^*$ via
$B(X)=i_X B$ and similarly we may regard $\beta$ as an element of
$\wedge^2 V$, i.e. a bivector.  This corresponds to the
observation that $\so(V\oplus V^*)=\wedge^2(V\oplus
V^*)=\End(V)\oplus \wedge^2 V^*\oplus\wedge^2V$.

By exponentiation, we obtain certain important orthogonal
symmetries of $T\oplus T^*$ in the identity component of
$SO(V\oplus V^*)$.
\begin{example}[$B$-transform]
First let $B$ be as above, and consider
\begin{equation}
\exp(B)=\left(\begin{matrix}1&\\B&1\end{matrix}\right),
\end{equation}
an orthogonal transformation sending $X+\xi\mapsto X+\xi+i_XB$. It
is useful to think of $\exp(B)$ a \emph{shear} transformation,
which fixes projections to $T$ and acts by shearing in the $T^*$
direction.  We will sometimes refer to this as a B-transform.
\end{example}
\begin{example}[$\beta$-transform]\label{betatrans}
Similarly, let $\beta$ be as above, and consider
\begin{equation}
\exp(\beta)=\left(\begin{matrix}1&\beta\\&1\end{matrix}\right),
\end{equation}
an orthogonal transformation sending $X+\xi\mapsto
X+\xi+i_\xi\beta$. Again, it is useful to think of $\exp(\beta)$ a
shear transformation, which fixes projections to $T^*$ and acts by
shearing in the $T$ direction.  We will sometimes refer to this as
a $\beta$-transform.
\end{example}
\begin{example}[$GL(V)$ action]
If we choose $A\in\so(V\oplus V^*)$ as above, then since
\begin{equation}
\exp(A)=\left(\begin{matrix}\exp A&\\&(\exp
A^*)^{-1}\end{matrix}\right),
\end{equation}
we see that we have a distinguished diagonal embedding of
$GL^+(V)$ into the identity component of $SO(V\oplus V^*)$. Of
course, we can extend this to a map
\begin{equation}
P\mapsto\left(\begin{matrix}P&\\&P^{*^{-1}}\end{matrix}\right)
\end{equation}
of the full $GL(V)$ into $SO(V\oplus V^*)$.  Note that $SO(V\oplus
V^*)$ has two connected components, and the two connected
components of $GL(V)$ do map into different components of
$SO(V\oplus V^*)$.
\end{example}

\section{Maximal isotropic subspaces}\label{maxiso}

A subspace $L<V\oplus V^*$ is \emph{isotropic} when $\IP{X,Y}=0$
for all $X,Y\in L$. Since we are in signature $(m,m)$, the maximal
dimension of such a subspace is $m$, and if this is the case, $L$
is called maximal isotropic. Maximal isotropic subspaces of
$V\oplus V^*$ are also called \emph{linear Dirac structures}
(see~\cite{Courant}).  Note that $V$ and $V^*$ are examples of
maximal isotropics. The space of maximal isotropics is
disconnected into two components, and elements of these are said
to have \emph{even} or \emph{odd} parity (sometimes called
helicity), depending on whether they share their connected
component with $V$ or not, respectively. This situation becomes
more transparent after studying the following two examples.
\begin{example}
Let $E\leq V$ be any subspace.  Then consider the subspace
\begin{equation*}
E\oplus\Ann( E) < V\oplus V^*,
\end{equation*}
where $\Ann (E)$ is the annihilator of $E$ in $V^*$. Then this is
a maximally isotropic subspace.
\end{example}
\begin{example}
Let $E\leq V$ be any subspace, and let $\varepsilon\in \wedge^2
E^*$.  Regarding $\varepsilon$ as a skew map $E\lra E^*$ via
$X\mapsto i_X\eps$, consider the following subspace, analogous to
the graph of $\eps$:
\begin{equation*}
L(E,\eps)=\left\{X+\xi\in E\oplus V^*\ :\
\xi\big|_{E}=\eps(X)\right\}.
\end{equation*}
Then if $X+\xi,Y+\eta\in L(E,\eps)$, we check that
\begin{align*}
\IP{X+\xi,Y+\eta}&=\tfrac{1}{2}(\xi(Y)+\eta(X))\\
&=\tfrac{1}{2}(\eps(Y,X)+\eps(X,Y))=0,
\end{align*}
showing that $L(E,\eps)$ is a maximal isotropic subspace.
\end{example}
Note that the second example specializes to the first by taking
$\eps=0$. Furthermore note that $L(V,0)=V$ and $L(\{0\},0)=V^*$.
It is not difficult to see that every maximal isotropic is of this
form:
\begin{prop}\label{Eeps}
Every maximal isotropic in $V\oplus V^*$ is of the form
$L(E,\eps)$.
\end{prop}
\begin{proof}
Let $L$ be a maximal isotropic and define $E=\pi_VL$, where
$\pi_V$ is the canonical projection $V\oplus V^*\rightarrow V$.
Then since $L$ is maximal isotropic, $L\cap V^*=\Ann(E)$.  Finally
note that $E^*=V^*/\Ann(E)$, and define $\eps:E\rightarrow E^*$
via $\eps: e\mapsto \pi_{V^*}(\pi_V^{-1}(e)\cap L)\in
V^*/\Ann(E)$. Then $L=L(E,\eps)$.
\end{proof}

The integer $k=\dim\Ann(E)=m-\dim\pi_V(L)$ is an invariant
associated to any maximal isotropic in $V\oplus V^*$.
\begin{defn}\label{type}
The \emph{type} of a maximal isotropic $L(E,\eps)$, is the
codimension $k$ of its projection onto $V$.
\end{defn}
Since a $B$-transform preserves projections to $V$, it does not
affect $E$:
\begin{align*}
\exp B(L(E,\eps))&=L(E,\eps+i^*B),
\end{align*}
where $i:E\hookrightarrow V$ is the natural inclusion. Hence
$B$-transforms do not change the type of the maximal isotropic. In
fact, we see that by choosing $B$ and $E$ accordingly, we can
obtain any maximal isotropic as a $B$-transform of $L(E,0)$.

On the other hand, $\beta$-transforms do modify projections to
$V$, and therefore may change the dimension of $E$.  To see how
this occurs more clearly, we express the maximal isotropic as a
generalized graph from $V^*\rightarrow V$, i.e. define
$F=\pi_{V^*}L$ and $\gamma\in\wedge^2F^*$ given by
$\gamma(f)=\pi_{V}(\pi_{V^*}^{-1}(f)\cap L)$. As before, define
\begin{equation*}
L(F,\gamma)=\left\{X+\xi\in V\oplus F\ :\
X\big|_{F}=\gamma(\xi)\right\}.
\end{equation*}
Then, as happened in the B-field case,
\begin{equation*}
\exp \beta(L(F,\gamma))=L(F,\gamma+j^*\beta),
\end{equation*}
where this time $j:F\hookrightarrow V^*$ is the inclusion. Now,
the projection $E=\pi_VL(F,\gamma)$ always contains $L\cap
V=Ann(F)$, and if we take the quotient of $E$ by this subspace we
obtain the image of $\gamma$ in $F^*=V/\Ann(F)$:
\begin{equation*}
\frac{E}{L\cap V}=\frac{E}{\Ann(F)}=\mathrm{Im}(\gamma).
\end{equation*}
Therefore, we obtain the dimension of $E$ as a function of
$\gamma$:
\begin{equation*}
\dim E=\dim L\cap V+\rk\gamma.
\end{equation*}
Because $\gamma$ is a skew form, its rank is even. A
$\beta$-transform sends $\gamma\mapsto\gamma+j^*\beta$, which also
has even rank, and therefore we see that a $\beta$-transform,
which is in the identity component of $SO(V\oplus V^*)$, can be
used to change the dimension of $E$, and hence the type of
$L(E,\eps)$, by an even number.
\begin{prop} Maximal isotropics
$L(E,\eps)$ of even parity are precisely those of even type; those
of odd type have odd parity.  The generic even maximal isotropics
are those of type 0, whereas the generic odd ones are of type 1.
The least generic type is $k=m$, of which there is only one
maximal isotropic: $V^*$.
\end{prop}

Before we move on to the description of maximal isotropics using
pure spinors, we indicate that alternative splittings for $V\oplus
V^*$ should be considered:
\begin{remark}
While the maximal isotropics $V$ and $V^*$ are distinguished in
the vector space $V\oplus V^*$, all our results about the linear
algebra of $V\oplus V^*$ are portable to the situation where $L$
and $L'$ are any maximal isotropics in $V\oplus V^*$ such that
$L\cap L'=0$. Then the inner product defines an isomorphism
$L'\cong L^*$, and we obtain $V\oplus V^*=L\oplus L^*$.
\end{remark}

\section{Spinors for $V\oplus V^*$: exterior
forms}\label{spinandforms}

Let $CL(V\oplus V^*)$ be the Clifford algebra defined by the
relation
\begin{equation}\label{clifford}
v^2=\IP{v,v},\ \ \forall v\in V\oplus V^*,
\end{equation}
The Clifford algebra has a natural representation on
$S=\wedge^\bullet V^*$ given by
\begin{equation}\label{cliff}
(X+\xi)\cdot\varphi = i_X\varphi + \xi\wedge\varphi,
\end{equation}
where $X+\xi\in V\oplus V^*$ and $\varphi\in\wedge^\bullet V^*$.
We verify that this defines an algebra representation:
\begin{align*}
(X+\xi)^2\cdot\varphi &=
i_X(i_X\varphi+\xi\wedge\varphi)+\xi\wedge(i_X\varphi+\xi\wedge\varphi)\\
&=(i_X\xi)\varphi\\
&=\IP{X+\xi,X+\xi}\varphi,
\end{align*}
as required.  This representation is the standard spin
representation, and so we see that for $V\oplus V^*$ there is a
canonical choice of spinors: the exterior algebra on $V^*$. Since
in signature $(m,m)$ the volume element $\omega$ of a Clifford
algebra satisfies $\omega^2=1$, the spin representation decomposes
into the $\pm1$ eigenspaces of $\omega$ (the positive and negative
helicity spinors):
\begin{equation*}
S=S^+\oplus S^-,
\end{equation*}
and this is easily seen to be equivalent to the decomposition
\begin{equation*}
\wedge^\bullet V^*=\wedge^\text{ev}V^*\oplus\wedge^\text{odd}V^*.
\end{equation*}
While the splitting $S=S^+\oplus S^-$ is not preserved by the
whole Clifford algebra, $S^\pm$ are irreducible representations of
the spin group, which sits in the Clifford algebra as
\begin{equation*}
Spin(V\oplus V^*)=\{v_1\cdots v_r\ \big|\ v_i\in V\oplus V^*,
\IP{v_i,v_i}=\pm 1\ \text{and}\ r\ \text{even}\},
\end{equation*}
and which is a double cover of $SO(V\oplus V^*)$ via the
homomorphism
\begin{align*}
\rho&:Spin(V\oplus V^*)\lra SO(V\oplus V^*)\\
\rho&(x)(v)=xv x^{-1}\ \ x\in Spin(V\oplus V^*),\ v\in V\oplus
V^*.
\end{align*}
Earlier we described certain symmetries of $V\oplus V^*$ by
exponentiating elements of $\so (V\oplus V^*)$. Since
$\mathfrak{so}(V\oplus V^*)=\wedge^2(V\oplus V^*)$ sits naturally
inside the Clifford algebra, we can see how its components act on
the spin representation.  Note that the derivative of $\rho$,
given by
\begin{equation*}
d\rho_x(v)=xv-vx = [x,v],\ \ x\in \mathfrak{so}(V\oplus V^*),\ \
v\in V\oplus V^*,
\end{equation*}
must be the usual representation of $\mathfrak{so}(V\oplus V^*)$
on $V\oplus V^*$.

Let $\{e_i\}$ be a basis for $V$ and $\{e^i\}$ be the dual basis.
We calculate the action of $B$- and $\beta$-transforms on $S$, as
well as the more complicated case of the $GL(V)$ action:

\begin{example}[$B$-transform]
If $B=\frac{1}{2}B_{ij}e^i\wedge e^j$, $B_{ij}=-B_{ji}$ is a
2-form acting on $V\oplus V^*$ via $X+\xi\mapsto i_XB$, then its
image in the Clifford algebra is $\frac{1}{2}B_{ij}e^je^i$, since
\begin{align*}
e^i\wedge e^j: e_i\mapsto e^j,
\end{align*}
and
\begin{align*}
d\rho_{e^je^i}(e_i)&=(e^je^i)e_i-e_i(e^je^i)\\
&=(e^ie_i+e_ie^i)e^j=e^j.
\end{align*}
Its spinorial action on an exterior form $\varphi\in\wedge^\bullet
V^*$ is
\begin{equation*}
B\cdot\varphi =
\frac{1}{2}B_{ij}e^j\wedge(e^i\wedge\varphi)=-B\wedge\varphi.
\end{equation*}
Therefore, exponentiating, we obtain
\begin{equation}
e^{-B}\varphi=(1-B+\tfrac{1}{2}B\wedge B + \cdots)\wedge\varphi.
\end{equation}
\end{example}
\begin{example}[$\beta$-transform]
If $\beta=\frac{1}{2}\beta^{ij}e_i\wedge e_j$,
$\beta^{ij}=-\beta^{ji}$ is a 2-vector acting on $V\oplus V^*$ via
$X+\xi\mapsto i_\xi\beta$, then its image in the Clifford algebra
is $\frac{1}{2}\beta^{ij}e^je^i$, and its spinorial action on a
form $\varphi$ is
\begin{equation*}
\beta\cdot\varphi =
\frac{1}{2}\beta^{ij}i_{e_j}(i_{e_i}\varphi)=i_{\beta}\varphi.
\end{equation*}
Therefore, exponentiating, we obtain
\begin{equation}
e^\beta\varphi=(1+i_\beta+\tfrac{1}{2}i_\beta^2 + \cdots)\varphi.
\end{equation}
\end{example}
Next we must understand the inverse image of the diagonally
embedded $GL(V)\subset SO(V\oplus V^*)$ under the covering $\rho$.
The group $SO(V\oplus V^*)\cong SO(n,n)$ has maximal compact
subgroup $S(O(n)\times O(n))$ and hence has two connected
components, each with fundamental group $\ZZ_2\times\ZZ_2$.  Hence
the double cover $Spin(n,n)$, which has maximal compact
$S(Pin(n,0)\times_{\ZZ_2}Pin(0,n))$, has two connected components,
each with fundamental group $\ZZ_2$.  Note that the image of the
injective map $\rho_{*}:\pi_1(Spin(n,n))\lra\pi_1(SO(n,n))$ is the
diagonal subgroup $\{1,(\sigma,\sigma)\}$, where $\sigma$ is the
nontrivial element of $\pi_1(SO(n))$. Now, consider the natural
inclusion of $GL(V)$ in $SO(V\oplus V^*)$ via
\begin{equation*}
P\mapsto\left(\begin{matrix}P&\\&P^{*^{-1}}\end{matrix}\right).
\end{equation*}
Topologically, this maps the maximal compact of $GL(V)$, which is
a disjoint union of two copies of $SO(n)$, diagonally into
$S(O(n)\times O(n))$.  This means that upon restriction to
$GL(V)\subset SO(V\oplus V^*)$, the Spin double cover is a
\emph{trivial} covering. Therefore there are two distinguished
$GL(V)$ subgroups of $Spin(V\oplus V^*)$ (depending on which
branch is chosen for the non-identity component of $GL(V)$) which
we denote $GL(V)_1$ and $GL(V)_2$, both mapping isomorphically to
$GL(V)$ via the covering map. Their intersection is the connected
group $GL^+(V)$.  We are interested to see how these general
linear groups act on $\wedge^\bullet V^*$ via the spin
representation.
\begin{example}[$GL^+(V)$ action]
The action of exponentiating elements of $\End(V)<\so(V\oplus
V^*)$ is slightly more complicated, and becomes crucial when
constructing the spinor bundle later on.

If $A=A_i^je^i\otimes e_j$ is an endomorphism of $V$ then as an
element of $\so(V\oplus V^*)$ it acts on $V\oplus V^*$ via
$X+\xi\mapsto A(X)-A^*(\xi)$. That is to say, $e^i\otimes e_j$
maps $e_i\mapsto e_j$ and $e^j\mapsto -e^i$. Its image in the
Clifford algebra is $\frac{1}{2}A_i^j(e_je^i-e^ie_j)$, and its
spinorial action on $\varphi$ is
\begin{align*}
A\cdot\varphi&=\tfrac{1}{2}A_i^j(e_je^i-e^ie_j)\cdot\varphi\\
&=\tfrac{1}{2}A_i^j(i_{e_j}(e^i\wedge\varphi)-e^i\wedge i_{
e_j}\varphi)\\
&=\tfrac{1}{2}A_i^j\delta^i_j\varphi - A_i^je^i\wedge
i_{e_j}\varphi\\
&=-A^*\varphi + \tfrac{1}{2}(\Tr A) \varphi,
\end{align*}
where $\varphi\mapsto -A^*\varphi=-A^j_i e^i\wedge i_{e_j}\varphi$
is the usual action of $End(V)$ on $\wedge^\bullet V^*$.  Hence,
by exponentiation, the spinorial action of $GL^+(V)$ on
$\wedge^\bullet V^*$ is by
\begin{equation*}
g\cdot\varphi=\sqrt{\det g}(g^*)^{-1}\varphi,
\end{equation*}
i.e. as a $GL^+(V)$ representation the spinor representation
decomposes as
\begin{equation*}
S=\wedge^\bullet V^*\otimes (\det V)^{1/2}.
\end{equation*}
\end{example}
As we have seen, the $GL^+(V)$ action may be extended to a full
$GL(V)$ action in two natural ways, which we now describe.
\begin{example}
Let $A_{\pm}=\pm\lambda e^1\otimes e_1+e^2\otimes e_2 +\cdots +
e^n\otimes e_n$, where $\lambda$ is a positive real number.  Then
clearly $A_{\pm}\in GL^{\pm}(V)$.  We now describe the elements in
$Spin(V\oplus V^*)$ covering these.

Using the Clifford algebra norm, we see that $N(e_1+e^1)=1$ and
that $N(\lambda^{-1/2}e_1\pm\lambda^{1/2}e^1)=\pm 1$ (Note that we
take positive square roots only.) Hence we form the $Spin(V\oplus
V^*)$ elements
\begin{align*}
a_{\pm}&=(e_1+e^1)(\lambda^{-1/2}e_1\pm\lambda^{1/2}e^1)\\
&=\lambda^{-1/2}e^1e_1\pm\lambda^{1/2}e_1e^1.
\end{align*}
We now check that these elements actually cover $A_{\pm}$.  Note
that $a_{\pm}^{-1}=\lambda^{1/2}e^1e_1\pm\lambda^{-1/2}e_1e^1$:
\begin{align*}
\rho(a_{\pm})(e_1)&=a_{\pm} e_1 a^{-1}_{\pm}\\
&=(\lambda^{-1/2}e^1e_1\pm\lambda^{1/2}e_1e^1)e_1(\lambda^{1/2}e^1e_1\pm\lambda^{-1/2}e_1e^1)\\
&=\pm\lambda e_1,
\end{align*}
as required.  Hence the two elements of $Spin(V\oplus V^*)$ which
cover the element $A_{-}$ are $\pm a_-$.  By convention, let us
say that $a_-$ is in $GL^-(V)_1\subset Spin(V\oplus V^*)$ and that
$-a_-$ is in $GL^-(V)_2\subset Spin(V\oplus V^*)$. Now we may see
how they act via the Spin representation.  We take the unique
Koszul decomposition of the general form
$\varphi=e^1\wedge\varphi_0 + \varphi_1$ where
$i_{e_1}\varphi_0=i_{e_1}\varphi_1=0$.
\begin{align*}
a_- \cdot \varphi &=
(\lambda^{-1/2}e^1e_1-\lambda^{1/2}e_1e^1)\cdot\varphi\\
&=\lambda^{-1/2}e^1\wedge\varphi_0 - \lambda^{1/2}(\varphi-e^1\wedge \varphi_0)\\
&=-\lambda^{1/2}(-\lambda^{-1}e^1\wedge\varphi_0 +\varphi_1).\\
&=-(-\det A_-)^{1/2}A^{*^{-1}}_-\varphi
\end{align*}
\end{example}
Hence we deduce that the action of $g\in GL^-(V)$ is as follows:
\begin{equation}\label{spin12}
g\cdot\varphi=
\begin{cases}-(-\det g)^{1/2} g^{*^{-1}}\varphi &
\text{ if $g\in GL^-(V)_1$},\\
(-\det g)^{1/2}g^{*^{-1}}\varphi & \text{ if $g\in GL^-(V)_2$}.
\end{cases}
\end{equation}

\section{The bilinear form on spinors}

There is a bilinear form on spinors which behaves well under the
spin representation. We describe it here, following the treatment
of Chevalley~\cite{Chevalley}.  For $V\oplus V^*$ this bilinear
form coincides with the Mukai pairing of forms~\cite{Mukai}.

Since we have the splitting $V\oplus V^*$ into maximal isotropics,
the exterior algebras $\wedge^\bullet V$ and $\wedge^\bullet V^*$
are subalgebras of $CL(V\oplus V^*)$.  In particular, $\det V$ is
a distinguished line inside $CL(V\oplus V^*)$, and it generates a
left ideal with the property that upon choosing a generator
$f\in\det V$, every element of the ideal has a unique
representation as $sf$, $s\in\wedge^\bullet V^*$.  This defines an
action of the Clifford algebra on $\wedge^\bullet V^*$ by
\begin{equation*}
(\rho(x)s)f=xsf\ \ \ \forall x\in CL(V\oplus V^*),
\end{equation*}
which is the same action as that defined by~(\ref{cliff}).

Having expressed the spin representation in this way, we proceed
to write down the bilinear form. Let $\alpha$ be the main
antiautomorphism of the Clifford algebra, i.e. that determined by
the tensor map $v_1\otimes\cdots\otimes v_k\mapsto v_k\otimes
\cdots \otimes v_1$. Now let $s,t\in \wedge^\bullet V^*$ be
spinors and consider the bilinear expression
\begin{equation*}
\begin{split}
\xymatrix{S\otimes S\ar[r]^{\IPS{(,)}}&\det V^*}\\
\IPS{(s,t)}=(\alpha(s)\wedge t)_{\text{top}},
\end{split}
\end{equation*}
where $()_{\text{top}}$ indicates taking the top degree component
of the form.

We can express $(\IPS{,})$ in the following way, using any
generator $f\in\det V$:
\begin{equation}\label{beta}
\begin{split}
(i_{f}\IPS{(s,t)})f&=(i_{f}(\alpha(s)\wedge t))f\\
&=(\rho(\alpha(f)) (\alpha(s)t) )f\\
&=\alpha(f)\alpha(s)tf\\
&=\alpha(sf)tf.
\end{split}
\end{equation}
From this description the behaviour of $\IPS{}$ under the action
of the Clifford algebra becomes clear. Let $v\in V\oplus V^*$:
\begin{align*}
(i_f\IPS{(v\cdot s,v\cdot t)})f&=\alpha(vsf) vtf\\
&=\alpha(sf)\alpha(v)vtf\\
&=\IP{v,v}i_f\IPS{(s,t)}f,
\end{align*}
showing that $\IPS{(v\cdot s,v\cdot t)}=\IP{v,v}\IPS{(s,t)}$, so
that in particular $\IPS{(g\cdot s,g\cdot t)}=\pm\IPS{(s,t)}$ for
any $g \in Spin(V\oplus V^*)$.
\begin{prop}
The bilinear form $\xymatrix{S\otimes S\ar[r]^{\IPS{(,)}}&\det
V^*}$ is invariant under the identity component of Spin:
\begin{equation*}
\IPS{(x\cdot s,x\cdot t)}=\IPS{(s,t)}\ \ \forall\  x \in
Spin_0(V\oplus V^*).
\end{equation*}
For example, $\IPS{(\exp B \cdot s,\exp B\cdot s)}=\IPS{(s,t)}$,
for any $B\in \wedge^2 V^*$.
\end{prop}

This bilinear form is non-degenerate, and can be symmetric or
skew-symmetric depending on the dimension of $V$:
\begin{equation*}
\IPS{(s,t)}=(-1)^{m(m-1)/2}\IPS{(t,s)}.
\end{equation*}
\begin{prop}
The nondegenerate bilinear form $\IPS{(,)}$ is symmetric for
$m\equiv 0\ \text{or}\ 1\pmod 4$ and skew-symmetric otherwise.
\end{prop}
The behaviour of $\IPS{(,)}$ with respect to the decomposition
$S=S^+\oplus S^-$ depends on $m=\dim V$ in the following way,
which is clear from the description of $\beta$ in terms of wedge
product.
\begin{prop}
If $m\equiv 0\pmod 2$ then $\IPS$ is zero on $S^+\times S^-$ (and
hence $S^-\times S^+$); if $m\equiv 1\pmod 2$ then $\IPS$ is zero
on $S^+\times S^+$ and $S^-\times S^-$.
\end{prop}

\begin{example}
Suppose $V$ is 4-dimensional; then $\IPS{}$ is symmetric, and the
even spinors are orthogonal to the odd spinors.  The inner product
of even spinors $\rho=\rho_0+\rho_2+\rho_4$ and
$\sigma=\sigma_0+\sigma_2+\sigma_4$ is given by
\begin{align*}
\IPS{(\rho,\sigma)}&=((\rho_0 - \rho_2 + \rho_4)\wedge(\sigma_0 +
\sigma_2+ \sigma_4))_4\\
&=\rho_0\sigma_4 - \rho_2\sigma_2 + \rho_4\sigma_0.
\end{align*}
The inner product of odd spinors $\rho=\rho_1+\rho_3$ and
$\sigma=\sigma_1+\sigma_3$ is given by
\begin{align*}
\IPS{(\rho,\sigma)}&=((\rho_1 - \rho_3)\wedge(\sigma_1 +
\sigma_3))_4\\
&=\rho_1\sigma_3 - \rho_3\sigma_1.
\end{align*}
\end{example}

\section{Pure spinors}
Let $\varphi$ be any nonzero spinor. Then we define its \emph{null
space} $L_\varphi<V\oplus V^*$, as follows:
\begin{equation}
L_\varphi=\{v\in V\oplus V^*\ :\ v\cdot\varphi=0\},
\end{equation}
and it is clear that $L_\varphi$ depends equivariantly on
$\varphi$ under the spin representation:
\begin{equation}
L_{g\cdot\varphi}=\rho(g)L_\varphi\ \ \forall\ g\in Spin(V\oplus
V^*).
\end{equation}
The key property of null spaces is that they are isotropic: if
$v,w\in L_\varphi$, then
\begin{equation}
2\IP{v,w}\varphi = (vw+wv)\cdot\varphi=0,
\end{equation}
implying that $\IP{v,w}=0\ \forall\ v,w\in L_\varphi$.
\begin{defn}
A spinor $\varphi$ is called \emph{pure} when $L_\varphi$ is
maximally isotropic, i.e. has dimension $m$.
\end{defn}
\begin{example}Let $1\in \wedge^\bullet V^*$ be the unit spinor.
Then the null space is
\begin{equation*}
\{X+\xi\in V\oplus V^*\ :\ (i_X+\xi\wedge)1=0\}= V,
\end{equation*}
and $V<V\oplus V^*$ is maximally isotropic, equal to $L(V,0)$ as
we saw in section~\ref{maxiso}. Hence $1$ is a pure spinor.  Of
course, we may apply any spin transformation to $1$ to obtain more
pure spinors; for instance, let $B\in \wedge^2 V^*$ and form
$\varphi=e^B\wedge 1 = e^B$, which has maximal null space
\begin{equation*}
N_\varphi = \{X-i_XB\ :\ X\in V\} = L(V,-B).
\end{equation*}
\end{example}
\begin{example}
Let $\theta\in \wedge^\bullet V^*$ be a nonzero 1-form. Then its
null space is
\begin{equation*}
\{X+\xi\in V\oplus V^*\ :\ (i_X+\xi\wedge)\theta=0\}=
\ker\theta\oplus \langle\theta\rangle =L(\ker\theta,0),
\end{equation*}
which is also maximal isotropic. Hence $\theta$, and therefore
$e^B\theta$ for any $B\in\wedge^2 V^*$, is a pure spinor.
\end{example}

We will refer to Chevalley~\cite{Chevalley} for the main
properties of pure spinors, and summarize the results here. Every
maximal isotropic subspace of $V\oplus V^*$ is represented by a
unique pure line in the spin bundle $S$. This line must lie in
$S^+$ for even maximal isotropics and in $S^-$ for odd ones.  The
cone of pure spinors in $\PP(S^+)$ and $\PP(S^-)$ is defined by a
set of quadratic equations; in the even case these are the
equations relating the different degree components of $e^B$ (the
generic element), and in the odd case they are the equations
relating the different degree components of $e^B\theta$, for
$\theta$ a 1-form.

The intersection properties of maximal isotropics can also be
obtained from the pure spinors, using the bilinear form
$\IPS{(,)}$:
\begin{prop}[\cite{Chevalley}, III.2.4.]\label{indexzero}
Maximal isotropics $L,L'$ satisfy $\dim L\cap L'=0$ if and only if
their pure spinor representatives $\varphi,\varphi'$ satisfy
\begin{equation*}
\IPS{(\varphi,\varphi')}\neq 0.
\end{equation*}
\end{prop}
The bilinear form also provides us with the operation of
``squaring'' the spinor, as follows.  The bilinear form determines
an isomorphism
\begin{equation*}
S\lra S^*\otimes\det V^*,
\end{equation*}
and, tensoring with $\text{id}:S\rightarrow S$, we obtain an
isomorphism
\begin{equation*}
S\otimes S\xrightarrow{\ \Phi\ }S\otimes S^*\otimes\det
V^*=\wedge^\bullet(V\oplus V^*)\otimes\det V^*.
\end{equation*}
Chevalley shows that this map takes a pure spinor line to the
determinant line of the maximal isotropic it defines:
\begin{prop}[\cite{Chevalley}, III.3.2.]\label{squaredet}
Let $U_L<\wedge^\bullet V^*$ be a pure spinor line representing
the maximal isotropic $L<V\oplus V^*$.  Then $\Phi$ determines an
isomorphism
\begin{equation*}
U_L\otimes U_L\xrightarrow{\ \Phi\ } \det L\otimes \det V^*,
\end{equation*}
where $\det L<\wedge^m(V\oplus V^*)$ is the determinant line.
\end{prop}

In the remainder of this section, we will provide an expression
for the pure spinor line associated to any maximal isotropic
$L(E,\eps)$.
\begin{lemma}
Let $E\leq V$ be any subspace of codimension $k$. The maximal
isotropic $L(E,0)=E\oplus \Ann (E)$ is associated to the pure
spinor line $\det(\Ann(E))<\wedge^k V^*$.
\end{lemma}
\begin{proof}
If $\varphi=\theta_1\wedge\cdots\wedge\theta_k$ is a nonzero
element of $\det(\Ann(E))$, then clearly $(X+\xi)\cdot\varphi=0$
if and only if $X\in E$ and $\xi\in\Ann(E)$, as required.
\end{proof}
As we saw in section \ref{maxiso}, any maximal isotropic
$L(E,\eps)$ may be expressed as the $B$-transform of $L(E,0)$ for
$B$ chosen such that $i^*B=\eps$.  Although $\eps$ is not in
$\wedge^2 V^*$, we may abuse notation and write
\begin{equation}
L(E,\eps)=\exp(\eps) ( L(E,0)),
\end{equation}
where in this equation $\eps$ is understood to represent any
2-form $B\in \wedge^2 V^*$ such that $i^*B=\eps$, where
$i:E\hookrightarrow V\otimes\CC$ is the inclusion. Passing now to
the spinorial description of maximal isotropics, and making use of
the previous lemma, we obtain a description of any pure spinor:
\begin{prop}\label{evenodd}
Let $L(E,\eps)$ be any maximal isotropic. Then the pure spinor
line $U_L$ defining it is given by
\begin{equation}
U_L=\exp(\eps) (\det\Ann (E)).
\end{equation}
To be more precise, let $(\theta_1,\ldots,\theta_k)$ be a basis
for $\Ann(E)$, and let $B\in\wedge^2 V^*$ be any 2-form such that
$i^*B=\eps$.  Then the following spinor represents the maximal
isotropic $L(E,\eps)$:
\begin{equation}\label{represent}
\varphi_L = c\exp(B)\theta_1\wedge\cdots\wedge\theta_k,\ \ c\neq
0,
\end{equation}
and any pure spinor can be expressed this way.  Note that even
maximal isotropics are represented by even forms and odd maximal
isotropics by odd forms.
\end{prop}

\section{Complexification and the real index}

The natural inner product $\IP{,}$ extends by complexification to
$(V\oplus V^*)\otimes\CC$, and all of our results concerning
maximal isotropics and spinors for $V\oplus V^*$ can be extended
by complexification to $(V\oplus V^*)\otimes\CC$.  We summarize
our results in this new context.
\begin{prop}\label{algtype}
Let $V$ be a real vector space of dimension $m$.  A maximal
isotropic subspace $L<(V\oplus V^*)\otimes\CC$ of type
$k\in\{0,\ldots,m\}$ is specified equivalently by the following
data:
\begin{itemize}
\item A complex subspace $L<(V\oplus V^*)\otimes\CC$, maximal
isotropic with respect to $\IP{,}$, and such that
$E=\pi_{V\otimes\CC}L$ has complex dimension $m-k$;
\item A complex subspace $E<V\otimes\CC$ such that $\dim_\CC E = m-k$, together with a complex
2-form $\eps\in\wedge^2E^*$;
\item A complex spinor line
$U_L<\wedge^\bullet(V^*\otimes\CC)$ generated by
\begin{equation}\label{formform}
\varphi_L = c\exp(B+i\omega)\theta_1\wedge\cdots\wedge\theta_k,
\end{equation}
where $(\theta_1,\ldots,\theta_k)$ are linearly independent
complex 1-forms in $V^*\otimes\CC$, $B$ and $\omega$ are the real
and imaginary parts of a complex 2-form in
$\wedge^2(V^*\otimes\CC)$, and $c\in\CC$ is a nonzero scalar.
\end{itemize}
\end{prop}

A new ingredient which appears when considering the complexified
situation is the complex conjugate, which acts on all associated
structures $L$, $E$, and $U_L$.  We use it to define the concept
of \emph{real index}, introduced in~\cite{Trautman}.

\begin{defn}
Let $L<(V\oplus V^*)\otimes\CC$ be a maximal isotropic subspace.
Then $L\cap\overline{L}$ is real, i.e. the complexification of a
real space: $L\cap\overline{L}=K\otimes\CC$, for $K<V\oplus V^*$.
The real index $r$ of the maximal isotropic $L$ is defined by
\begin{equation*}
r=\dim_\CC L\cap \overline{L}=\dim_\RR K.
\end{equation*}
For example, a real maximal isotropic $L<V\oplus V^*$ has real
index $m$.
\end{defn}

Since the parity of $L$ is determined by its intersection with the
real subspace $V\otimes\CC<(V\oplus V^*)\otimes\CC$, it is clear
that $L$ and $\overline{L}$ must have the same parity, implying
that $\dim L\cap\overline{L}\equiv m\pmod 2$, i.e.
\begin{equation}\label{realindex}
r\equiv m \pmod 2,
\end{equation}
showing that the real index must be even or odd depending on the
dimension of $V$.

\section{Functorial property of Dirac structures}

Maximal isotropics (linear Dirac structures) have the interesting
functorial property, noticed by Weinstein and used profitably
in~\cite{BurRadko}, that they can be pulled back and pushed
forward along any linear map $f:V\rightarrow W$ of vector spaces.
In particular, if $L<V\oplus V^*$ is a maximal isotropic then
\begin{equation}\label{pushforward}
f_*L=\{f(X)+\eta\in W\oplus W^*\ :\ X+f^*\eta\in L\}
\end{equation}
is a maximal isotropic subspace of $W\oplus W^*$.  Similarly, if
$M<W\oplus W^*$ is maximal isotropic then
\begin{equation}\label{pullback}
f^*M=\{X+f^*\eta\in V\oplus V^*\ :\ f(X)+\eta\in M\}
\end{equation}
is a maximal isotropic subspace of $V\oplus V^*$.

\section{The spin bundle for $T\oplus T^*$}

In this section we will transport our algebraic work on $V\oplus
V^*$ to a manifold.  Let $M$ be a smooth manifold of real
dimension $m$, with tangent bundle $T$.  Then consider the direct
sum of the tangent and cotangent bundles $T\oplus T^*$.  This
bundle is endowed with the same canonical bilinear forms and
orientation we described on $V\oplus V^*$.  Therefore, while we
are aware of the fact that $T\oplus T^*$ is associated to a
$GL(m)$ principal bundle, we may also view it as having natural
structure group $SO(m,m)$.

It is well-known that an oriented bundle with Euclidean structure
group $SO(n)$ admits spin structure if and only if the second
Stiefel-Whitney class vanishes, i.e. $w_2(E)=0$.  The situation
for bundles with metrics of indefinite signature is more
complicated, and was worked out by Karoubi in~\cite{Karoubi}. We
summarize his results:

If an orientable bundle $E$ has structure group $SO(p,q)$, we can
always reduce the structure group to its maximal compact subgroup
$S(O(p)\times O(q))$.  This reduction is equivalent to the choice
of a maximal positive definite subbundle $E^+<E$, which allows us
to write $E$ as the direct sum $E=E^+\oplus E^-$, where
$E^-=(E^+)^\perp$ is negative definite.
\begin{prop}[\cite{Karoubi}, 1.1.26]
The $SO(p,q)$ bundle $E$ admits $Spin(p,q)$ structure if and only
if $w_2(E^+)=w_2(E^-)$.
\end{prop}
In the special case of $T\oplus T^*$, which has signature $(m,m)$,
the positive and negative definite bundles $E^\pm$ project
isomorphically via $\pi_T:T\oplus T^*\rightarrow T$ onto the
tangent bundle.  Hence the condition $w_2(E^+)=w_2(E^-)$ is always
satisfied for $T\oplus T^*$, yielding the following result.
\begin{prop}
The $SO(m,m)$ bundle $T\oplus T^*$ always admits $Spin(m,m)$
structure.
\end{prop}
As usual, the exact sequence
\begin{equation*}
\xymatrix{1\ar[r]&\ZZ_2\ar[r]&Spin(m,m)\ar[r]&SO(m,m)\ar[r]&1}
\end{equation*}
informs us that the set of spin structures is an affine space
under $H^1(M,\ZZ_2)$, the group of real line bundles, which act on
the associated spinor bundles by tensor product.

While we have given an argument for the existence of spin
structure using characteristic classes, it is clear from the
results of section~\ref{spinandforms} that the spin bundles can be
constructed.  In that section we showed that there are two
subgroups of $Spin(V\oplus V^*)$, namely $GL(V)_1$ and $GL(V)_2$,
which map isomorphically to $GL(V)\subset SO(V\oplus V^*)$ under
the spin homomorphism. Therefore we can use the principal $GL(m)$
principal bundle of frames in $T$ to form an associated
$Spin(m,m)$ principal bundle, and this can be done in two
different ways, producing two spin structures.  Let the associated
spinor bundles be $S_1$ and $S_2$.  Since these bundles are
associated to the original $GL(m)$ bundle of frames, by
decomposing the spin representation according to $GL(m)_1$ and
$GL(m)_2$ we can express $S_1$ and $S_2$ in terms of well-known
associated bundles to the frame bundle.

In section~\ref{spinandforms} we did just this, and obtained the
result that the spin representation decomposes under $GL(V)_1$ as
\begin{equation*}
S=\wedge^\bullet V^*\otimes \det V |\det V|^{-1/2},
\end{equation*}
and under $GL(V)_2$ as
\begin{equation*}
S=\wedge^\bullet V^*\otimes |\det V|^{1/2}.
\end{equation*}
Note that $|\det V|^{1/2}$ indicates the representation of $GL(V)$
sending $g\mapsto |\det g|^{1/2}$. Alternatively we could think of
$|\det V|^{1/2}$ as the 1-dimensional vector space of
$1/2$-densities on $V^*$, i.e. maps $v:\det V^*\rightarrow\RR$
such that $v(\lambda\omega)=|\lambda|^{1/2}v(\omega)$ for all
$\lambda\in\RR$ and $\omega\in\det V^*$.

Therefore when we form the spin bundles $S_1$ and $S_2$ by
association, we obtain
\begin{align*}
S_1&=\wedge^\bullet T^*\otimes\det T |\det T|^{-1/2}\\
S_2&=\wedge^\bullet T^*\otimes|\det T|^{1/2}.
\end{align*}
Since $|\det T|^{1/2}$ is isomorphic to the trivial bundle, we see
that there is always a choice of spin structure such that the spin
bundle is (non-canonically) isomorphic to the exterior algebra
$\wedge^\bullet T^*$.  In any case, we will be primarily
interested not with sections of the spin bundle itself but of its
projectivisation. Therefore we can use the fact that the
projectivisation of any spin bundle for $T\oplus T^*$ is
canonically isomorphic to the projectivised differential forms:
\begin{equation}
\PP(S)=\PP(\wedge^\bullet T^*).
\end{equation}

The bilinear form $\IPS$ is inherited by the spin bundle, and
further by the projectively isomorphic $\wedge^\bullet T^*$, in
the form
\begin{equation*}
\IPS(\varphi,\psi)=(\alpha(\varphi)\wedge\psi)_{\text{top}},
\end{equation*}
which is invariant under the action of $Spin_0(T\oplus T^*)$ and
is covariant under diffeomorphisms. '

In previous sections we have studied the correspondence between
maximal isotropics in $V\oplus V^*$ and certain lines in the spin
bundle.  In the same way, maximal isotropic subbundles of $T\oplus
T^*$ correspond to sections of the projectivised bundle of
differential forms, or equivalently line subbundles of
$\wedge^\bullet T^*$.  If we are fortunate and the line bundle is
trivial, then it is possible to represent the maximal isotropic
subbundle by a global differential form.  For example, the maximal
isotropic $T<T\oplus T^*$ is represented by the line generated by
the differential form $1\in\wedge^\bullet T^*$, whereas the
maximal isotropic $T^*$, represented by the line $\det
T^*<\wedge^\bullet T^*$, cannot be given by a global form if $M$
is non-orientable.

\clearpage
\thispagestyle{empty}
\cleardoublepage
\chapter{The Courant bracket}\label{cb}

The Courant bracket is a natural bracket operation on the smooth
sections of $T\oplus T^*$. It was first introduced in its present
form by T. Courant~\cite{Courant}, in the context of his work with
Weinstein~\cite{CourWein}.  It also was implicit in the
contemporaneous work of Dorfman~\cite{Dorfman}. Courant and
Weinstein used it to define a new geometrical structure called a
\emph{Dirac} structure, which successfully unifies Poisson
geometry and presymplectic geometry (the geometry defined by a
real closed 2-form) by expressing each structure as a maximal
isotropic subbundle of $T\oplus T^*$. The integrability condition,
namely that the subbundle be closed under the Courant bracket,
specializes to the usual integrability conditions in the Poisson
and presymplectic cases.  In the same way, we will use the Courant
bracket to define the integrability of generalized complex
structures.

We begin this section with an introduction to Lie algebroids, a
class of vector bundles with structure that closely resembles that
of the tangent bundle. Lie algebroids are particularly useful for
at least two reasons. First, they provide a sufficiently general
framework to accommodate a unified treatment of many kinds of
geometry, including Poisson, foliated, (pre)symplectic, and as we
shall show, complex and CR geometry, in addition to many new
(generalized) types of geometry. Second, they provide a way to
handle, in a smooth (non-singular) fashion, structures which at
first glance appear to acquire singularities at certain loci in
the manifold.

We then introduce the Courant bracket on $T\oplus T^*$, describing
its basic properties, and showing that it fails to fit into the
framework of Lie algebroid theory. Indeed, when its properties are
systematized, we obtain the axioms of a Courant algebroid, first
introduced in~\cite{LWX}.  We follow the treatment in Roytenberg's
thesis~\cite{Roytenberg} for this material.  The Courant algebroid
structure of $T\oplus T^*$ is useful to us for two main reasons.
First, it provides a source of new Lie algebroids by restriction
to subbundles (as in the case of Dirac structures, mentioned
above). Second, its natural group of symmetries includes not only
diffeomorphisms, but also closed 2-forms, which act in a way
familiar to physicists as the action of the $B$-field. As a
consequence, all geometrical structures defined in terms of the
Courant bracket can be transformed naturally by a $B$-field.

Finally, we investigate the notion of `twisting' the Courant
bracket, and the relation of this to the theory of gerbes with
connection.

Throughout this section we will make extensive use of the
following identities relating Lie derivative and interior product:
\begin{equation*}
\LL_X=i_Xd+di_X,\ \ \ \ \LL_{[X,Y]}=[\LL_X,\LL_Y],\ \ \ \
i_{[X,Y]}=[\LL_X,i_Y],
\end{equation*}
where $X,Y$ are vector fields.  We follow the convention that for
any differential form $\rho$, the interior product is a
contraction with the first argument: $i_X\rho= \rho(X,\cdots)$.

\section{Lie algebroids}

A Lie algebroid, first defined by Pradines in~\cite{Pradines} and
explored in depth by Mackenzie in~\cite{Mackenzie}, is a vector
bundle $L$ on a smooth manifold $M$, equipped with a Lie bracket
$[,]$ on $C^\infty(L)$ and a smooth bundle map $a:L\rightarrow T$,
called the \emph{anchor}. The anchor must induce a Lie algebra
homomorphism $a:C^\infty(L)\rightarrow C^\infty(T)$, i.e.
\begin{equation}\label{LA1}
a([X,Y])=[a(X),a(Y)]\ \ \forall X,Y\in C^\infty(L),
\end{equation}
and the following Leibniz rule must be satisfied:
\begin{equation}\label{LA2}
[X,fY]=f[X,Y]+(a(X)f)Y\ \ \forall X,Y\in C^\infty(L), f\in
C^\infty(M).
\end{equation}

\begin{example}[The tangent bundle]
The tangent bundle is itself a Lie algebroid, taking the identity
map as anchor.  It is useful to think of a Lie algebroid as a
generalization of the tangent bundle.
\end{example}
\begin{example}[Foliations]
Any integrable sub-bundle of the tangent bundle defines a Lie
algebroid, choosing $a$ to be the inclusion map.
\end{example}
\begin{example}[The Atiyah sequence]
Let $\pi:P\rightarrow M$ be a principal $G$-bundle on the manifold
$M$. Then $G$-invariant vector fields on $P$ are given by sections
of the vector bundle $TP/G\rightarrow M$. This bundle has a Lie
algebroid structure defined by the Lie bracket on $C^\infty(TP)$
and the surjective anchor $\pi_*$, which defines an exact sequence
of vector bundles on $M$:
\begin{equation*}
\xymatrix{0\ar[r]&\mathfrak{g}\ar[r]&TP/G\ar[r]^{\ \
\pi_*}&T\ar[r]& 0},
\end{equation*}
where $\mathfrak{g}$ is the adjoint bundle associated to $P$.
\end{example}

The notion of Lie algebroid can obviously be complexified,
defining a \emph{complex} Lie algebroid by requiring $L$ to be a
complex bundle and $a:L\rightarrow T\otimes\CC$ a complex map,
satisfying complexified conditions~(\ref{LA1}),~(\ref{LA2}). While
the theory of Dirac structures uses real Lie algebroids, in this
thesis complex Lie algebroids are particularly important.  We will
produce many examples of complex Lie algebroids in what follows,
but we have some immediately at hand.
\begin{example}[Complex structures]
If $M$ is a complex manifold then $T_{1,0}<T\otimes\CC$ is a
complex bundle closed under the Lie bracket. Using the inclusion
map as anchor, $T_{1,0}$ is a complex Lie algebroid.
\end{example}
\begin{example}[CR structures]
A CR (Cauchy-Riemann) structure on a real $2n-1$-dimensional
manifold is a complex $n-1$-dimensional sub-bundle $L<T\otimes\CC$
which satisfies $\dim L\cap\overline{L}=0$ and which is closed
under the Lie bracket.  $L$ is then a complex Lie algebroid, with
the inclusion map as anchor.
\end{example}

We now define several structures which exist naturally on Lie
algebroids by generalizing the well-known case of the tangent
bundle. We begin with the fact that the Lie bracket on vector
fields has a natural $\ZZ$-graded extension to multivector fields
$C^\infty(\wedge^\bullet T)$, called the Schouten bracket. This
can be generalized to the context of Lie algebroids.

\begin{defn}
Let $L$ be a Lie algebroid.  The Schouten bracket acting on
sections $X_1\wedge\cdots\wedge X_p\in C^\infty(\wedge^p L)$,
$Y_1\wedge\cdots\wedge Y_q\in C^\infty(\wedge^q L)$ is as follows:
\begin{equation*}
[X_1\wedge\cdots\wedge X_p,Y_1\wedge\cdots\wedge
Y_q]=\sum_{i,j}(-1)^{i+j}[X_i,Y_j]\wedge X_1\wedge\cdots
\hat{X_i}\cdots\wedge X_p\wedge
Y_1\wedge\cdots\hat{Y_j}\cdots\wedge Y_q,
\end{equation*}
and $[X,f]=-[f,X]=a(X)f$ for $X\in C^\infty(L)$ and $f\in
C^\infty(M)$. This bracket makes $C^\infty(\wedge^\bullet L)$ into
a graded Lie algebra where the degree $k$ component is
$C^\infty(\wedge^{k+1}L)$.  That is,
\begin{align*}
[A,B]&=-(-1)^{(a-1)(b-1)}[B,A]\ \ \text{and}\\
[A,[B,C]]&=[[A,B],C]+(-1)^{(a-1)(b-1)}[B,[A,C]],
\end{align*}
for all $A\in C^\infty(\wedge^aL),B\in C^\infty(\wedge^bL),C\in
C^\infty(\wedge^cL)$.

Furthermore, if $A\in C^\infty(\wedge^{a}L)$, then
$ad_A=[A,\cdot]$ is a derivation of degree $a-1$ of the exterior
multiplication on $C^\infty(\wedge^\bullet L)$:
\begin{align*}
ad_A(B\wedge C)=ad_A(B)\wedge C + (-1)^{(a-1)b}B\wedge ad_A(C).
\end{align*}
The two interacting graded algebra structures of
$C^\infty(\wedge^\bullet L)$, namely the exterior product and the
Schouten bracket, make it into a Poisson
superalgebra\footnote{Since the two graded structures differ in
parity, it is an `odd Poisson superalgebra', otherwise known as a
\emph{Gerstenhaber algebra}.}.
\end{defn}

In addition to the Schouten bracket of vector fields, smooth
manifolds are equipped with the exterior derivative operator $d$,
a derivation of degree 1 of the algebra of differential forms. The
exterior derivative operator can be defined in terms of the Lie
bracket, and for this reason we may generalize it to Lie
algebroids as follows:

\begin{defn}
The Lie algebroid derivative $d_L$ is a first order linear
operator from $C^\infty(\wedge^k L^*)$ to $C^\infty(\wedge^{k+1}
L^*)$ defined by
\begin{align*}
d_L\sigma(X_0,\ldots,X_k)=&\sum_i
(-1)^ia(X_i)\sigma(X_0,\ldots,\hat{X_i},\ldots,X_k)
\\&+\sum_{i<j}(-1)^{i+j}\sigma([X_i,X_j],X_0,\ldots,\hat{X_i},\ldots,\hat{X_j},\ldots,X_k),
\end{align*}
where $\sigma\in C^\infty(\wedge^k L^*)$ and $X_i\in C^\infty(L)$.
Its principal symbol $s(d_L):T^*\otimes\wedge^k L^*\rightarrow
\wedge^{k+1}L^*$ is given by $a^*:T^*\rightarrow L^*$ composed
with wedge product, i.e.
\begin{equation*}
s_\xi(d_L)=a^*(\xi)\wedge\cdot\ \ ,
\end{equation*}
where $\xi\in T^*$.  The operator $d_L$ satisfies $d_L^2=0$ due to
the Jacobi identity for $[,]$, and therefore
$(C^\infty(\wedge^\bullet L^*),d_L)$ is a natural differential
complex associated with any Lie algebroid.
\end{defn}

By analogy, we define interior product and Lie derivative for Lie
algebroids:
\begin{defn}\label{Liederivative}
Let $X\in C^\infty(L)$.  Then the interior product $i_X$ is the
degree $-1$ derivation on $C^\infty(\wedge^\bullet L^*)$ defined
by $i_X\sigma=\sigma(X,\cdots)$, and the Lie derivative $\LL^L_X$
is defined by the Cartan formula
\begin{equation*}
\LL^L_X=d_Li_X+i_Xd_L.
\end{equation*}
\end{defn}

So far, we have described the intrinsic algebraic structures
present in Lie algebroids.  It is important to understand what
these structures imply for the underlying geometry of the
manifold.  In particular, as described by Courant~\cite{Courant},
every Lie algebroid induces a generalized foliation in the sense
of Sussmann~\cite{Sussmann}, which we now explain.

A foliated manifold $M$ is one which has been expressed as a
disjoint union of subsets called \emph{leaves}.  A leaf is a
connected submanifold (injective immersion) $l\subset M$ such that
any point $p\in l$ has a neighbourhood $U\subset M$ where the
connected component of $p$ in $l\cap U$ is an embedded submanifold
of $M$.  A usual foliation is one where all leaves have the same
dimension, whereas a generalized foliation allows the dimension of
the leaves to vary. The second main result in Sussmann's
paper~\cite{Sussmann} describes necessary and sufficient
conditions on a distribution $\Delta$ for it to be integrable into
such a generalized foliation (by which is meant that at any point
$m\in M$, the tangent space to the leaf through $m$ is precisely
$\Delta(m)$).  Essentially, Sussmann's theorem is a more powerful
version of the classical Frobenius integrability theorem.  We
needn't use Sussmann's full theorem, but only a corollary of his
main result:

\begin{theorem}[\cite{Sussmann}, Theorem 8.1.]
Let $\Delta_D$ be a distribution spanned by a collection $D\subset
C^\infty(T)$ of smooth vector fields.  This is called a
\emph{smooth distribution}.  Such a distribution is said to be of
\emph{finite type} if, for every $m\in M$, there exist smooth
vector fields $X^1,\ldots,X^n\in D$ such that
\begin{itemize}
\item $X^1(m),\ldots,X^n(m)$ span $\Delta_D(m)$ and
\item For every $X\in D$, there exists a neighbourhood $U$ of $m$ and smooth functions $c^{i}_k\in C^\infty(U)$ such
that for all $i$,
\begin{equation*}
[X,X^i]=\sum_kc^{i}_kX^k.
\end{equation*}
\end{itemize}
If $\Delta_D$ is of finite type, then it is integrable to a
generalized foliation as described above.
\end{theorem}

Any real Lie algebroid $L$ with anchor $a$ produces a distribution
$\Delta_D=a(L)$, and since it is the image of a smooth bundle map,
it is spanned by the smooth vector fields $D=a(C^\infty(L))$.
Hence it is a smooth distribution.  Furthermore, for any point
$m\in M$ we may choose a local basis of sections $X^1,\ldots,
X^n\in C^\infty(U,L)$ in some neighbourhood of $m$.  Then
$a(X^1),\ldots,a(X^n)$ certainly span $\Delta_D$ in $U$ and by the
Lie algebroid property~(\ref{LA1}), we see that
\begin{equation*}
[a(X^i),a(X^j)]=a([X^i,X^j])=a(\sum_kc^{ij}_kX^k)=\sum_kc^{ij}_ka(X^k),
\end{equation*}
for some $c^{ij}_k\in C^\infty(U)$, implying that $\Delta_D$ is of
finite type.  Hence by the theorem, we conclude that a real Lie
algebroid induces a generalized foliation on the manifold $M$.
Furthermore, since the rank of the smooth bundle map $a$ is a
lower semi-continuous function (every point has a neighbourhood in
which the rank \emph{does not decrease}), we conclude that the
dimension of the leaf is also a lower semi-continuous function.
\begin{prop}[\cite{Courant}, Theorem 2.1.3.]\label{smoothdist}
If $L$ is a real Lie algebroid on $M$ with anchor $a$, then
$\Delta=a(L)$ is a smooth integrable distribution in the sense of
Sussmann, implying that $M$ can be expressed locally as a disjoint
union of embedded submanifolds (called leaves) such that at any
point $m\in M$, the tangent space to the leaf through $m$ is
precisely $\Delta(m)$. Furthermore, the dimension of the leaf,
$\dim\Delta$, is a lower semi-continuous function on the manifold.
\end{prop}

If $L$ is a complex Lie algebroid, then let $E<T\otimes\CC$ denote
the image under the anchor, i.e. $E=a(L)$.  The complex
distribution $E$ induces two real distributions
$\Delta\subset\Theta\subset T$, defined by $E+\overline
E=\Theta\otimes\CC$ and $E\cap\overline E=\Delta\otimes\CC$. While
$\Theta$ need not be integrable, or even involutive, certainly
$\Delta$ is involutive, but we haven't enough information to
decide its Sussmann integrability.  To establish an analogous
result to the previous proposition we will need an extra
assumption.

In particular, consider the case where
$E+\overline{E}=T\otimes\CC$.  Because of this, the real bundle
map $i(a-\overline{a}):L\rightarrow T$ is surjective, hence the
kernel is a smooth real sub-bundle $K<L$:
\begin{equation*}
K=\left\{X\in L\ :\ a(X)=\overline{a(X)}\right\}.
\end{equation*}
Projecting $K$ to $T$ via the anchor $a$, we obtain precisely the
distribution $\Delta$. As a result, we see that $\Delta$ is a
smooth distribution.  It is easy to check now that it is of finite
type:  for any point $m\in M$, let $X^1,\ldots,X^n$ be a local
basis of sections for $C^\infty(U,K)$ in some neighbourhood $U$ of
$m$. Then $a(X^1),\ldots,a(X^n)$ span $\Delta$ in $U$.
Furthermore, since
$a([X^i,X^j])=[a(X^i),a(X^j)]=\overline{[a(X^i),a(X^j)]}=\overline{a([X^i,X^j])}$,
we see that $[X^i,X^j]\in C^\infty(U,K)$, and therefore
$[X^i,X^j]=\sum_k c^{ij}_k X^k$, giving
\begin{equation*}
[a(X^i),a(X^j)]=\sum_k c^{ij}_k a(X^k),
\end{equation*}
which implies that $\Delta$ is of finite type. Hence we obtain the
following result.
\begin{prop}\label{cxliedist}
Let $L$ be a complex Lie algebroid on $M$ with anchor $a$, and
such that $E+\overline{E}=T\otimes\CC$, where $E=a(L)$.  Let
$\Delta$ be the real distribution defined by $\Delta\otimes\CC=
E\cap\overline{E}$.  Then $\Delta$ is a smooth integrable
distribution in the sense of Sussmann, defining a generalized
foliation of $M$. Furthermore, the dimension of the leaf,
$\dim\Delta$, is a lower semi-continuous function on the manifold.
\end{prop}

Besides the generalized foliation, a complex Lie algebroid of this
type induces a \emph{transverse complex structure} on this
foliation in a sense which we now describe.

A complex distribution $E<T\otimes\CC$ of constant complex
codimension $k$ on a real $n$-manifold $M$ is \emph{integrable}
if, in some neighbourhood $U$ of each point $m\in M$, there exist
complex functions $f_{1},\ldots, f_{k}\in C^\infty(U,\CC)$ such
that $\{df_{1},\ldots,df_{k}\}$ are linearly independent at each
point in $U$ and annihilate all complex vector fields lying in
$E$.  By the Newlander-Nirenberg theorem~\cite{N-N}, we know that
$E$ is integrable if the following conditions are satisfied:
\begin{itemize}
\item
$E$ is involutive (closed under Lie bracket), and
\item $\dim E\cap \overline{E}$ is constant, and
\item $E+\overline{E}$ is involutive as well.
\end{itemize}
In this situation, the functions $f_{1},\ldots f_{k}$ are complex
coordinates transverse to the foliation determined by
$E\cap\overline{E}$.  In other words, every point $m\in M$ has a
neigbourhood isomorphic, as a smooth manifold with complex
distribution, to an open set in $\RR^{n-2k}\times\CC^{k}$, which
has natural complex distribution spanned by $\{ \partial/{\partial
x_1},\ldots,\partial/{\partial x_{n-2k}},\partial/{\partial
z_1},\ldots,\partial/{\partial z_{k}} \}$.

In our situation of a complex Lie algebroid satisfying
$E+\overline{E}=T\otimes\CC$, if we restrict our attention to
points $m\in M$ which are \emph{regular}, in the sense that the
leaf dimension is constant in a neighbourhood $U$ of $m$, then the
above conditions are satisfied, and we obtain a transversal
complex structure (at regular points).
\begin{prop}\label{transversecx}
Let $L$ be a complex Lie algebroid on the real $n$-manifold $M$
with anchor $a$, and such that $E+\overline{E}=T\otimes\CC$, where
$E=a(L)$. Let $m\in M$ be a regular point for the Lie algebroid,
i.e. a point where $k=\dim E\cap\overline{E}$ is locally constant.
Then in some neighbourhood $U$ of $m$, there exist complex
functions $z_{1},\ldots, z_{k}\in C^\infty(U,\CC)$ such that
$\{dz_1,\ldots,dz_{k}\}$ are linearly independent at each point in
$U$ and annihilate all complex vector fields lying in $E$, i.e. we
have a transverse complex structure to the foliation, at regular
points.
\end{prop}
\begin{remark}
Note that a transverse complex structure on a foliation implies
that any smooth section of the foliation in a domain of regular
points inherits an integrable complex structure. The coordinates
$\{z_i\}$ can be restricted to the section and serve as complex
coordinates.
\end{remark}

Before we proceed to study the Courant bracket, which will provide
us with many more examples of complex Lie algebroids, we provide
the definition of a structure called a \emph{Lie bialgebroid},
which was introduced by Mackenzie and Xu~\cite{MackXu} as the
infinitesimal object corresponding to a Poisson groupoid. The Lie
algebroids we will study will appear naturally in bialgebroid
pairs, and this will be crucially important when we study the
deformation theory of generalized complex structures.

\begin{defn}[Lie bialgebroid]
Let $L$ be a Lie algebroid and suppose its dual bundle $L^*$ also
has the structure of a Lie algebroid.  Then $(L,L^*)$ is a Lie
bialgebroid if the Lie algebroid derivative
$d_L:C^\infty(L^*)\rightarrow C^\infty(\wedge^2L^*)$ is a
derivation of the Schouten bracket $[,]_{L^*}$ on $C^\infty(L^*)$,
in the sense that
\begin{equation*}
d_L[X,Y]=[d_LX,Y]+[X,d_LY].
\end{equation*}
The Lie bialgebroid condition is self-dual, in the sense that
$(L,L^*)$ is a Lie bialgebroid if and only if $(L^*,L)$ is.
Furthermore, the Lie bialgebroid condition is equivalent to
requiring that $d_L$ is a derivation of degree $1$ of the graded
algebra $(C^\infty(\wedge^\bullet(L^*)),[,]_{L^*})$.  Both these
facts are proven in~\cite{K-S}.
\end{defn}
\begin{example}
The most obvious Lie bialgebroid is simply $(T,T^*)$, where we
take the usual Lie algebroid structure on the tangent bundle $T$
and the trivial structure on $T^*$ (zero bracket and anchor). Then
the exterior derivative is certainly a derivation of the trivial
bracket.
\end{example}

\section{The Courant bracket and Courant algebroids}

The Courant bracket is a skew-symmetric bracket defined on smooth
sections of $T\oplus T^*$, given by
\begin{equation*}
[X+\xi,Y+\eta]=[X,Y]+\LL_X\eta-\LL_Y\xi-\frac{1}{2}d(i_X\eta-i_Y\xi),
\end{equation*}
where $X+\xi,Y+\eta\in C^\infty(T\oplus T^*)$.

Note that on vector fields the Courant bracket reduces to the Lie
bracket $[X,Y]$; in other words, if $\pi:T\oplus T^*\rightarrow T$
is the natural projection,
\begin{equation}\label{C1}
\pi([A,B])=[\pi (A),\pi (B)],
\end{equation}
for any $A,B\in C^\infty(T\oplus T^*)$.  On the other hand, on
1-forms the Courant bracket vanishes. We will address this basic
asymmetry of the bracket in section~\ref{ISO}.

The Courant bracket is not a Lie bracket, since it fails to
satisfy the Jacobi identity. Therefore, although the projection
map $\pi$ acts as a Lie algebroid anchor, $(T\oplus T^*,[,])$ is
not a Lie algebroid. However it is interesting to examine how it
fails to be a Lie algebroid. The Jacobiator is a trilinear
operator which measures the failure to satisfy the Jacobi
identity:
\begin{equation*}
\Jac(A,B,C)=[[A,B],C]+[[B,C],A]+[[C,A],B],
\end{equation*}
where $A,B,C\in C^\infty(T\oplus T^*)$. The Jacobiator can be
usefully expressed as the derivative of a quantity which we will
call the \emph{Nijenhuis} operator, for reasons which will become
clear later.  For this reason, one can say that the Courant
bracket satisfies the Jacobi identity \emph{up to an exact term}.
We now prove this and two other basic properties of the Courant
bracket; these results are implicit in~\cite{LWX}; we provide
proofs here which will be useful for later development.
\begin{prop}
\begin{equation}\label{C2}
\Jac(A,B,C)=d\left(\Nij(A,B,C)\right),
\end{equation}
where $\Nij$ is the Nijenhuis operator:
\begin{equation*}
\Nij(A,B,C)=\tfrac{1}{3}\left(\IP{[A,B],C}+\IP{[B,C],A}+\IP{[C,A],B}\right).
\end{equation*}
Here $\IP{,}$ is the inner product on $T\oplus T^*$ introduced in
the previous section.
\end{prop}
\begin{proof}
To prove this result, we introduce the Dorfman bracket operation
$\circ$ on $T\oplus T^*$ which is not skew, but whose skew
symmetrization is the Courant bracket:
\begin{equation*}
(X+\xi)\circ (Y+\eta) = [X,Y] + \LL_X\eta - i_Yd\xi.
\end{equation*}
The difference between the two brackets is as follows:
\begin{equation*}
[A,B]=A\circ B-d\IP{A,B},
\end{equation*}
and of course $[A,B]=\tfrac{1}{2}(A\circ B- B\circ A)$. The
advantage of the Dorfman bracket is that it satisfies a kind of
Leibniz rule\footnote{In fact, the Dorfman bracket makes
$C^\infty(T\oplus T^*)$ into a \emph{Loday} algebra.}:
\begin{equation*}
A\circ(B\circ C)=(A\circ B)\circ C + B\circ (A\circ C),
\end{equation*}
which is easily proved, setting $A=X+\xi,\ B=Y+\eta,\ C=Z+\zeta$:
\begin{align*}
(A\circ B)&\circ C + B\circ (A\circ C)\\
&= [[X,Y],Z]+[Y,[X,Z]] + \LL_{[X,Y]}\zeta -
i_Zd(\LL_X\eta-i_Yd\xi)+\LL_Y(\LL_X\zeta-i_Zd\xi)-i_{[X,Z]}d\eta\\
&=[X,[Y,Z]]+\LL_X\LL_Y\zeta-\LL_Xi_Zd\eta - \LL_Yi_Zd\xi + i_Zdi_Yd\xi\\
&=[X,[Y,Z]]+\LL_X(\LL_Y\zeta-i_Zd\eta) - i_{[Y,Z]}d\xi\\
&=A\circ(B\circ C),
\end{align*}
as required.  Now note that
\begin{align*}
[[A,B],C]&=[A,B]\circ C -d\IP{[A,B],C}\\
&=(A\circ B - d\IP{A,B})\circ C - d\IP{[A,B],C}\\
&=(A\circ B)\circ C - d\IP{[A,B],C},
\end{align*}
where we have used the fact that $\mu\circ C=0$ whenever $\mu$ is
a closed 1-form.  Finally, we express the Jacobiator as follows
(c.p. indicates cyclic permutations):
\begin{align*}
\Jac(A,B,C) &= [[A,B],C]+\text{c.p.}\\
&=\tfrac{1}{4}((A\circ B)\circ C - C\circ (A\circ B)-(B\circ A)\circ C +C\circ(B\circ A)+\text{c.p.})\\
&=\tfrac{1}{4}(A\circ(B\circ C)-B\circ(A\circ C)-C\circ(A\circ B) - B\circ(A\circ C)+A\circ(B\circ C)+C\circ(B\circ A)+\text{c.p})\\
&=\tfrac{1}{4}(A\circ(B\circ C)-B\circ(A\circ C)+\text{c.p.})\\
&=\tfrac{1}{4}((A\circ B)\circ C +\text{c.p.})\\
&=\tfrac{1}{4}([[A,B],C]+d\IP{[A,B],C} + \text{c.p.})\\
&=\tfrac{1}{4}(\Jac(A,B,C) + 3d (\Nij(A,B,C))),
\end{align*}
which implies that $\Jac(A,B,C)=d(\Nij(A,B,C))$, as required.
\end{proof}

The next proposition describes the failure of the Courant bracket
to satisfy the second Lie algebroid axiom~(\ref{LA2}).
\begin{prop} Let $f\in C^\infty(M)$. Then the Courant bracket
satisfies
\begin{equation}\label{C3}
[A,fB]=f[A,B]+(\pi(A)f)B - \IP{A,B}df.
\end{equation}
\end{prop}
\begin{proof}
Let $A=X+\xi$ and $B=Y+\eta$, so that
\begin{align*}
[X+\xi,f(Y+\eta)]&=[X,fY]+\LL_X
f\eta-\LL_{fY}\xi-\tfrac{1}{2}d(i_X(f\eta)-i_{fY}\xi)\\
&=f[X+\xi,Y+\eta]+(Xf)Y + (Xf)\eta -
(i_Y\xi)df-\tfrac{1}{2}(i_X\eta-i_Y\xi)df\\
&=f[X+\xi,Y+\eta]+(Xf)(Y+\eta) -\IP{X+\xi,Y+\eta}df,
\end{align*}
as required.\end{proof} We see yet again that the Courant bracket
differs from being a Lie algebroid by exact terms. Both
properties~(\ref{C2}) and~(\ref{C3}) of the Courant bracket
demonstrate that it is intimately linked to the natural inner
product $\IP{,}$, a fact which we will employ frequently.  There
is a further property which highlights the relationship between
$[,]$ and $\IP{,}$, which we will find useful:
\begin{prop}Differentiation of the natural inner product can be
expressed in terms of the Courant bracket thus:
\begin{equation}\label{C4}
\pi(A)\IP{B,C}=\IP{[A,B]+d\IP{A,B},C}+\IP{B,[A,C]+d\IP{A,C}}.
\end{equation}
\end{prop}
\begin{proof}
In terms of the Dorfman bracket, we wish to prove that
\begin{equation*}
\pi(A)\IP{B,C}=\IP{A\circ B,C}+\IP{B,A\circ C};
\end{equation*}
set $A=X+\xi,\ B=Y+\eta,\ C=Z+\zeta$.  Then we have
\begin{align*}
\IP{A\circ B,C}+\IP{B,A\circ C}&=\tfrac{1}{2}(i_{[X,Y]}\zeta+i_{Z}(\LL_X\eta-i_Yd\xi) + i_{[X,Z]}\eta+i_{Y}(\LL_X\zeta-i_Zd\xi) )\\
&=\tfrac{1}{2}(L_Xi_Y\zeta+L_Xi_Z\eta)\\
&=\tfrac{1}{2}i_Xd(i_Y\zeta+i_Z\eta)\\
&=\pi(A)\IP{B,C},
\end{align*}
as required.
\end{proof}

The fundamental properties~(\ref{C1}), (\ref{C2}), (\ref{C3}), and
(\ref{C4}) make $(T\oplus T^*,\ \IP{,},\ [,],\ \pi)$ into the
motivating example of a \emph{Courant algebroid}, the definition
of which is the first main result of the paper~\cite{LWX}.
\begin{defn}[\cite{LWX}, Definition 2.1]
A Courant algebroid is a vector bundle $E$ equipped with a
nondegenerate symmetric bilinear form $\IP{,}$ as well as a
skew-symmetric bracket $[,]$ on $C^\infty(E)$, and with a smooth
bundle map $\pi:E\rightarrow T$ called the \emph{anchor}. This
induces a natural differential operator
$\mathcal{D}:C^\infty(M)\rightarrow C^\infty(E)$ via the
definition $\IP{\DD f,A}=\tfrac{1}{2}\pi(A)f$ for all $f\in
C^\infty(M)$ and $A\in C^\infty(E)$. These structures must be
compatible in the following sense:
\begin{itemize}
\item[C1)] $\pi([A,B])=[\pi(A),\pi(B)]\ \ \ \forall A,B\in C^\infty(E)$,
\item[C2)] $\Jac(A,B,C)=\DD\left(\Nij(A,B,C)\right)\ \ \ \forall A,B,C\in
C^\infty(E)$,
\item[C3)] $[A,fB]=f[A,B]+(\pi(A)f)B - \IP{A,B}\DD f,\ \ \ \forall A,B\in C^\infty(E),\ f\in
C^\infty(M)$,
\item[C4)] $\pi\circ\DD=0$, i.e. $\IP{\DD f,\DD g}=0\ \ \ \forall
f,g\in C^\infty(M)$.
\item[C5)] $\pi(A)\IP{B,C}=\IP{[A, B]+\DD\IP{A,B},C}+\IP{B,[A,
C]+\DD\IP{A,C}}\ \ \ \forall A,B,C\in C^\infty(E)$,
\end{itemize}
where the Jacobiator $\Jac(\cdot,\cdot,\cdot)$ and the Nijenhuis
operator $\Nij(\cdot,\cdot,\cdot)$ are as defined before.
\end{defn}
\begin{remark}
The definition of $\DD$ implies that it satisfies a Leibniz rule,
i.e. $\DD(fg)=f\DD(g)+\DD(f)g$.  In~\cite{Uchino} it is noted that
this Leibniz property, together with axioms $(C1)$ and $(C5)$,
imply not only the definition of $\DD$ but also axioms $(C3)$ and
$(C4)$.
\end{remark}
\begin{remark}
As with Lie algebroids, we define the concept of \emph{complex}
Courant algebroid in the obvious way.
\end{remark}

\section{Symmetries of the Courant bracket; the B-field}

The Lie bracket of smooth vector fields is a canonically defined
structure on a manifold; that is to say, it is invariant under
diffeomorphisms.  In fact, there are no other symmetries of the
tangent bundle preserving the Lie bracket.
\begin{prop}
Let $(f,F)$ be an automorphism of the tangent bundle $\pi:TM\lra
M$ of a smooth manifold $M$, i.e. a pair of diffeomorphisms
$f:M\lra M$, $F:TM\lra TM$ such that the diagram
\begin{equation*}
\xymatrix{TM\ar[r]^{F}\ar[d]_{\pi}&TM\ar[d]^\pi\\M\ar[r]^f&M}
\end{equation*}
commutes and $F$ is a linear map on each fibre.  Suppose also that
$F$ preserves the Lie bracket, i.e. $F([X,Y])=[F(X),F(Y)]$ for all
vector fields $X,Y$.  Then $F$ must equal $f_*$, the derivative of
$f$.
\end{prop}
\begin{proof}
Note that $(f,f_*)$ is an automorphism of the tangent bundle
preserving the Lie bracket. Therefore, setting $G=f_*^{-1}\circ
F$, the pair $(Id,G)$ is also an automorphism preserving the Lie
bracket.  In particular, for any vector fields $X,Y$ and $h\in
C^\infty(M)$ we have $G([hX,Y])=[G(hX),G(Y)]$, or expanding,
\begin{equation*}
G([hX,Y])=G(h[X,Y] - Y(h)X)=hG([X,Y])-Y(h)G(X),
\end{equation*}
while, on the other hand,
\begin{align*}
[G(hX),G(Y)])&=h[G(X),G(Y)] - G(Y)(h)G(X)\\
&=hG([X,Y]) - G(Y)(h)G(X),
\end{align*}
so that $Y(h)G(X)=G(Y)(h)G(X)$ for all $X,Y,h$.  This can only
hold when $G(Y)=Y$ for all vector fields $Y$, i.e. $G=Id$,
yielding finally that $F=f_*$.
\end{proof}

In the case of $T\oplus T^*$, however, the situation is not so
simple.  While the Courant bracket and natural inner product are
invariant under diffeomorphisms, they have an additional symmetry,
which we call a \emph{B-field transformation}.  Let $B$ be a
smooth 2-form and view it as a map $T\lra T^*$ via interior
product $X\mapsto i_XB$.  This is the natural Lie algebra action
of $\wedge^2 T^*<\mathfrak{so}(T\oplus T^*)$ on $T\oplus T^*$.
Then the invertible bundle map given by exponentiating $B$, namely
\begin{equation*}
e^B=\left(\begin{matrix}1&0\\B&1\end{matrix}\right):\ X+\xi\mapsto
X+\xi+i_XB
\end{equation*}
is orthogonal, since $B^*=-B$ implies that $(e^B)^*e^B=e^{B-B}=1$.
It is useful to think of $e^B$ a \emph{shear} transformation,
which fixes projections to $T$ and acts by shearing in the $T^*$
direction.
\begin{prop}\label{2form}
The map $e^B$ is an automorphism of the Courant bracket if and
only if $B$ is closed, i.e. $dB=0$.
\end{prop}
\begin{proof}
Let $X+\xi,Y+\eta\in C^\infty(T\oplus T^*)$ and let $B$ be a
smooth 2-form.  Then
\begin{align*}
[e^B(X&+\xi),e^B(Y+\eta)]\\
&=[X+\xi+i_XB,Y+\eta+i_YB]\\
&=[X+\xi,Y+\eta] + [X,i_YB] + [i_XB,Y]\\
&=[X+\xi,Y+\eta] + L_Xi_YB-\frac{1}{2}di_Xi_YB
-L_Yi_XB+\frac{1}{2}di_Yi_XB\\
&=[X+\xi,Y+\eta]+L_Xi_YB-i_YL_XB + i_Yi_XdB\\
&=[X+\xi,Y+\eta]+i_{[X,Y]}B + i_Yi_XdB\\
&=e^B([X+\xi,Y+\eta])+i_Yi_XdB.
\end{align*}
Therefore we see that $e^B$ is an automorphism of the Courant
bracket if and only if $i_Yi_XdB=0$ for all $X,Y$, which happens
precisely when $dB=0$.
\end{proof}

A natural question we may ask at this point is whether B-field
transforms and diffeomorphisms are the only orthogonal
automorphisms of the Courant bracket.
\begin{prop}
Let $(f,F)$ be an orthogonal automorphism of the direct sum
$T\oplus T^*$ for a smooth manifold $M$. Suppose also that $F$
preserves the Courant bracket, i.e. $F([A,B])=[F(A),F(B)]$ for all
sections $A,B\in C^\infty(T\oplus T^*)$. Then $F$ must be the
composition of a diffeomorphism of $M$ and a B-field
transformation.  To be more precise, the group of orthogonal
Courant automorphisms of $T\oplus T^*$ is the semidirect product
of $\text{Diff}(M)$ and $\Omega^2_{\text{closed}}(M)$.
\end{prop}
\begin{proof}
Note that if $f$ is a diffeomorphism, the map
$f_c=\left(\begin{smallmatrix}f_*&0\\0&(f^*)^{-1}\end{smallmatrix}\right)$
is an orthogonal automorphism of $T\oplus T^*$ preserving the
Courant bracket. Therefore, setting $G=f_c^{-1}\circ F$, the pair
$(Id,G)$ is also an orthogonal automorphism preserving the Courant
bracket. In particular, for any sections $A,B\in C^\infty(T\oplus
T^*)$ and $h\in C^\infty(M)$ we have $G([hA,B])=[G(hA),G(B)]$, or
expanding,
\begin{align*}
G([hA,B])&=G(h[A,B] - (B_Th)A-\IP{A,B}dh)\\
&=hG([A,B])-(B_Th)G(A) - \IP{A,B}G(dh),
\end{align*}
while, on the other hand,
\begin{align*}
[G(hA),G(B)])&=h[G(A),G(B)] - (G(B)_Th)G(A)-\IP{G(A),G(B)}dh\\
&=hG([A,B]) - (G(B)_Th)G(A)-\IP{G(A),G(B)}dh.
\end{align*}
Setting these equal and using orthogonality, we obtain
\begin{equation*}
(B_Th)G(A) + \IP{A,B}G(dh)= (G(B)_Th)G(A)+\IP{A,B}dh.
\end{equation*}
Choose $A=X, B=Y$, where $X,Y\in C^\infty(T)$ so that
$\IP{A,B}=0$. Then we have that $Y(h)G(X)=(G(Y)_Th)G(X)$ for all
$X,Y,h$. This can only hold when $G(Y)_T=Y$ for all vector
fields $Y$, implying that $G=\left(\begin{smallmatrix}1&*\\
*&*\end{smallmatrix}\right)$. With this in mind, the previous
equation becomes
\begin{equation*}
\IP{A,B}G(dh)=\IP{A,B}dh,
\end{equation*}
which implies that $G=\left(\begin{smallmatrix}1&0\\
*&1\end{smallmatrix}\right)$.  Orthogonality then forces $G=\left(\begin{smallmatrix}1&0\\
B&1\end{smallmatrix}\right)=e^B$ where $B$ is a skew 2-form, and
to preserve the Courant bracket $B$ must be closed.  Hence we have
that $F=f_c\circ e^B$, as required.
\end{proof}
\begin{remark}
All our results concerning symmetries of the Courant algebroid
$(T\oplus T^*,\ \IP{,},\ [,],\ \pi)$ hold for the complexified
situation, but although any closed complex 2-form will act as a
symmetry, we restrict the terminology ``B-field'' to only those
2-forms which are \emph{real}.
\end{remark}

\section{Dirac structures}

It is clear from our investigation of the Courant bracket that it
fails to be a Lie algebroid due to exact terms \emph{involving the
inner product $\IP{,}$}.  For this reason, if we were to find a
sub-bundle $L<(T\oplus T^*)\otimes\CC$ which was involutive
(closed under the Courant bracket) as well as being isotropic,
then the anomalous terms would vanish, and $(L,[,],\pi)$ would
define a Lie algebroid. Furthermore, we could take the image of
such a sub-bundle under a B-field symmetry, obtaining another Lie
algebroid. Even beyond this, there may be orthogonal
transformations of $T\oplus T^*$ which, while they may not be
symmetries of the entire Courant structure, may take $L$ to a Lie
algebroid nonetheless.  In these ways, we will manufacture many
Lie algebroids as sub-bundles of $(T\oplus T^*)\otimes\CC$.

In fact, the Courant bracket itself places a tight constraint on
which proper sub-bundles may be involutive a priori:
\begin{prop}
If $L<T\oplus T^*$ is involutive then $L$ must either be an
isotropic subbundle, or a bundle of type $\Delta\oplus T^*$ for
$\Delta$ a nontrivial involutive sub-bundle of $T$.  Similarly for
the complex case $L<(T\oplus T^*)\otimes\CC$.
\end{prop}
\begin{proof}
Suppose that $L<T\oplus T^*$ is involutive, but not an isotropic
subbundle, i.e. there exists $X+\xi\in C^\infty(L)$ which is not
null at some point $m\in M$, i.e. $\xi(X)_m\neq 0$. Then for any
$f\in C^\infty(M)$,
\begin{align*}
[X+\xi,f(X+\xi)]&= (Xf)(X+\xi)-\xi(X)df,
\end{align*}
implying that $df_m\in L$ for all $f$, i.e. $T_m^*\leq L_m$. Since
$T^*_m$ is isotropic, this inclusion must be proper, i.e.
$L_m=\Delta_m\oplus T^*_m$, where $\Delta_m=\ker
\pi_{T^*}|_L:L_m\rightarrow T^*_m$.  Hence the rank of $L$ must
exceed the maximal dimension of an isotropic sub-bundle, which for
a real $n$-manifold is simply $n$.  This implies that $T^*_m<L_m$
at every point $m$, and hence that
$\Delta=\ker\pi_{T^*}|_L:L\rightarrow T^*$ is a smooth sub-bundle
of $T$, which must itself be involutive.  Hence $L=\Delta\oplus
T^*$, as required.
\end{proof}

In the maximal isotropic case, there is a remarkable equivalence
between certain natural conditions on $L$ and the involutive
condition.
\begin{prop}
Let $L$ be a maximal isotropic sub-bundle of $T\oplus T^*$ (or its
complexification). Then the following are equivalent:
\begin{itemize}
\item $L$ is involutive,
\item $\Nij\big|_{L}=0$,
\item $\Jac\big|_L=0$.
\end{itemize}
\end{prop}
\begin{proof}
If $L$ is involutive then it is clear that $\Nij|_L=0$ and since
$\Jac(A,B,C)=d(\Nij(A,B,C))$, it is clear that this implies
$\Jac|_L=0$ as well.  What remains to show is that $\Jac|_L=0$
implies that $L$ is involutive.

Suppose then that $\Jac|_L=0$ but that $L$ is not involutive, so
that there exist $A,B,C\in C^\infty(L)$ such that
$\IP{[A,B],C}\neq 0$. Then for all $f\in C^\infty(M)$,
\begin{align*}
0=\Jac(A,B,fC)&=d(\Nij(A,B,fC))\\
&=\tfrac{1}{3}\IP{[A,B],C}df,
\end{align*}
which is a contradiction.  Hence $L$ must be involutive.
\end{proof}
\begin{defn}[Dirac structure]
A real, maximal isotropic sub-bundle $L<T\oplus T^*$ is called an
almost \emph{Dirac} structure.  If $L$ is involutive, then the
almost Dirac structure is said to be integrable, or simply a
\emph{Dirac} structure.  Similarly, a maximal isotropic and
involutive complex sub-bundle $L<(T\oplus T^*)\otimes\CC$ is
called a \emph{complex Dirac} structure.  All the Lie algebroids
we will be considering will be Dirac structures.  Note that this
definition still works if $T\oplus T^*$ is replaced with any real
or complex Courant algebroid.  Thus we may speak of Dirac
structures in an arbitrary Courant algebroid.
\end{defn}

\begin{remark}
Note that for an isotropic sub-bundle $L$, the restricted
Nijenhuis operator $\Nij|_L$ is actually tensorial and is a
section of $\wedge^3 L^*$. Hence the above theorem indicates that
the integrability of a Dirac structure $L$ is determined by the
vanishing of a tensor field $\Nij|_L$.  In particular, any almost
Dirac structure on a 2-dimensional surface $\Sigma$ is integrable
for this reason.
\end{remark}

We will now provide several main examples of Dirac structures.
\begin{example}[Symplectic geometry]\label{symp}
The tangent bundle $T$ is itself maximal isotropic and involutive,
hence defines a Dirac structure.  To this basic Dirac structure we
can apply any closed (possibly complex) 2-form $\omega\in
\Omega^2_{cl}(M)$ to obtain another involutive maximal isotropic.
Indeed, the maximal isotropic subspace
\begin{equation*}
e^\omega(T)=\{X+i_X\omega\ :\ X\in T\}
\end{equation*}
is involutive if and only if $d\omega=0$ (see
Proposition~\ref{2form}). Immediately we see that pre-symplectic
geometry, i.e. the geometry defined by a closed 2-form, can be
described by a Dirac structure.
\end{example}

\begin{example}[Poisson geometry]\label{pois}
Similarly, the cotangent bundle $T^*$ is maximal isotropic and
involutive (the Courant bracket vanishes on $T^*$), defining a
Dirac structure.  We may apply a bivector field $\beta\in
C^\infty(\wedge^2 T)$ to this basic Dirac structure, obtaining
\begin{equation*}\label{poisson}
L_\beta=e^\beta(T^*)=\{i_\xi\beta+\xi \ :\ \xi\in T^*\}.
\end{equation*}
Since $\Nij|_{L_\beta}$ is tensorial, it suffices to check the
involutivity of sections of the form $i_\xi\beta+\xi$ where
$\xi=df$ for $f\in C^\infty(M)$.  The bivector $\beta$ determines
a bracket on functions
\begin{equation*}
\{f,g\}=\beta(df,dg),
\end{equation*}
and it is a straightforward calculation that
\begin{equation*}
\Nij(i_{df}\beta+df,i_{dg}\beta+dg,i_{dh}\beta+dh)=\{\{f,g\},h\}+\{\{g,h\},f\}+\{\{h,f\},g\},
\end{equation*}
showing that $\Nij|_{L_\beta}=0$ if and only if the bracket
$\{,\}$ satisfies the Jacobi identity, that is, $\beta$ is a
Poisson structure, which is equivalent to $[\beta,\beta]=0$.  In
this case we see that while $\beta$ is not a symmetry of the
Courant bracket, it does take $T^*$ to an integrable Dirac
structure.
\end{example}
\begin{example}[Foliated geometry]
Let $\Delta<T$ be a smooth distribution of constant rank. Then
form the maximal isotropic subbundle
\begin{equation*}
\Delta\oplus \Ann(\Delta)<T\oplus T^*.
\end{equation*}
This almost Dirac structure is Courant involutive if and only if
$\Delta$ is an integrable distribution.  By the theorem of
Frobenius, this produces a foliation on the manifold.  From this
point of view, a foliation can also be described by a (real) Dirac
structure.
\end{example}
\begin{remark}
Note that any involutive sub-bundle of type $\Delta\oplus T^*$
contains a maximal isotropic involutive sub-bundle
$\Delta\oplus\Ann(\Delta)$.
\end{remark}
\begin{example}[Complex geometry]\label{cxdirac}
An almost complex structure $J\in\End(T)$ determines a complex
distribution, given by the $-i$-eigenbundle $T_{0,1}<T\otimes\CC$
of $J$.  Forming the maximal isotropic space
\begin{equation*}
L_J=T_{0,1}\oplus\Ann(T_{0,1})=T_{0,1}\oplus T^*_{1,0},
\end{equation*}
we see that if $L$ is Courant involutive, then since the vector
component of the Courant bracket is simply the Lie bracket, this
implies $T_{0,1}$ is Lie involutive, i.e. $J$ is integrable.
Conversely, if $J$ is integrable, then, letting $X+\xi,Y+\eta\in
C^\infty(T_{0,1}\oplus T^*_{1,0})$, we have
\begin{equation*}
[X+\xi,Y+\eta]=[X,Y]+i_X\delbar\eta-i_Y\delbar\xi,
\end{equation*}
which is clearly a section of $T_{0,1}\oplus T^*_{1,0}$. Hence $L$
is Courant involutive if and only if $J$ is integrable. In this
way, integrable complex structures can also be described by
(complex) Dirac structures.
\end{example}

Although the preceding examples are only a few of the possible
Dirac structures, they demonstrate that four completely separate
classical geometrical structures, with very different-looking
integrability conditions, are unified when considered as Dirac
structures. As we shall see, there are many additional advantages
to describing these geometries in this way.  We have already
encountered one: the fact that from this point of view, there is a
natural action of $\Omega^2_\text{cl}(M)$ on the geometries.
Another is the fact that under suitable conditions, almost Dirac
structures may be pushed forward and pulled back along smooth maps
$f:M\rightarrow N$ between manifolds.

\section{Lie bialgebroids and the Courant bracket}\label{ISO}

In the preceding section, we observed that the Courant bracket
interpolates between the condition $d\omega=0$ for a 2-form and
the condition $[\beta,\beta]=0$ for a bivector.  This curious
behaviour is related to the asymmetry of the Courant bracket,
which we now address.

It was observed in~\cite{LWX} that Courant algebroids can be
constructed out of Lie bialgebroids via a generalization of the
Drinfel'd double construction.  In particular, given a Lie
bialgebroid $(L,L^*)$, one can define an inner product and bracket
on the sections of $L\oplus L^*$ making this bundle into a Courant
algebroid, and both $L,L^*$ into Dirac structures in this Courant
algebroid. In the same paper, the converse is shown: given any two
Dirac structures $L,L'$ in a Courant algebroid which are
transversal, then the inner product may be used to identify
$L'=L^*$, and $(L,L^*)$ is a Lie bialgebroid.

With this in mind, it is clear that the asymmetrical form of the
Courant bracket on $T\oplus T^*$ is related to the asymmetry of
the Lie bialgebroid $(T,T^*)$.  By choosing a different pair of
transverse Dirac structures $L,L'<(T\oplus T^*)\otimes\CC$, the
Courant bracket may appear more even-handed with respect to $L$
and $L^*$.

\begin{theorem}[\cite{LWX}, Theorem 2.5]
Let $(L,L^*)$ be a Lie bialgebroid.  Then we have the following
inner product $\IP{,}$ on the bundle $L\oplus L^*$:
\begin{equation*}
\IP{A+\alpha,B+\beta}=\tfrac{1}{2}(\alpha(B)+\beta(A)).
\end{equation*}
We define also the following skew-symmetric bracket operation on
$C^\infty(L\oplus L^*)$:
\begin{align*}
[A+\alpha,B+\beta]&= [A,B] + \LL_\alpha B - \LL_\beta A
-\tfrac{1}{2}d_{L^*}(i_A\beta-i_B\alpha)\\
&+[\alpha,\beta] + \LL_A\beta- \LL_B\alpha +
\tfrac{1}{2}d_{L}(i_A\beta-i_B\alpha),
\end{align*}
where here the Lie derivative and interior product operators are
as in Definition~\ref{Liederivative}.  If $a,a_*$ are the anchors
for $L,L^*$, we define the bundle map
\begin{equation*}
\pi=a+a_*:L\oplus L^*\rightarrow T.
\end{equation*}
With these definitions, the structure $(L\oplus
L^*,[,],\IP{,},\pi)$ is a Courant algebroid.  The operator
$\DD:C^\infty(M)\rightarrow C^\infty(L\oplus L^*)$ clearly becomes
\begin{equation*}
\DD=d_L+d_{L^*}.
\end{equation*}
\end{theorem}

\begin{theorem}[\cite{LWX}, Theorem 2.6]\label{bial}
Let $(E,[,],\IP{,},\pi)$ be a Courant algebroid, and let $L,L'<E$
be Dirac sub-bundles transverse to each other, i.e. $E=L\oplus
L'$.  Then $L'=L^*$ using the inner product, and $(L,L')$ is a Lie
bialgebroid.  Applying the construction in the previous theorem to
$(L,L')$, we recover the original Courant algebroid structure on
$E$.
\end{theorem}

Given a splitting of a Courant algebroid $E=L\oplus L^*$ into a
Lie bialgebroid, we could attempt to create new Dirac structures
as we did for $T\oplus T^*$ in examples~\ref{symp} and~\ref{pois},
as graphs of elements in $\wedge^2L$ or $\wedge^2L^*$.  Without
loss of generality, let $\eps\in\wedge^2L^*$ and consider the
maximal isotropic sub-bundle
\begin{equation*}
L_\eps=\{A+i_A\eps\ :\ A\in L\}.
\end{equation*}
The condition for the involutivity of the almost Dirac structure
$L_\eps$ is the final main result of~\cite{LWX}:
\begin{theorem}[\cite{LWX}, Theorem 6.1]\label{masterequation}
The almost Dirac structure $L_\eps$, for $\eps\in\wedge^2L^*$, is
integrable if and only if $\eps$ satisfies the generalized
Maurer-Cartan equation
\begin{equation}\label{MC}
d_L\eps+\tfrac{1}{2}[\eps,\eps]=0.
\end{equation}
Here $d_L:C^\infty(\wedge^k L^*)\rightarrow
C^\infty(\wedge^{k+1}L^*)$ and $[,]$ is the Lie algebroid bracket
on $L^*$.
\end{theorem}

From this result we can finally understand how the Courant bracket
allows us to interpolate between $d\omega=0$ and
$[\beta,\beta]=0$: due to the asymmetric bialgebroid structure on
$(T,T^*)$, we see that $[,]$ vanishes on $T^*$ while $d_{T^*}=0$
on $T$.  Hence the Maurer-Cartan equation on a 2-form $B\in
\wedge^2 T^*$ is simply $dB=0$, whereas on a bivector
$\beta\in\wedge^2T$ it is $\tfrac{1}{2}[\beta,\beta]=0$.

The results of this section will be particularly important when we
study the deformation theory of generalized complex structures:
finding solutions to the Maurer-Cartan equation will correspond to
finding integrable deformations of the generalized complex
structure.  In this way, we solve an open problem stated
in~\cite{LWX}, namely, to find an interpretation of
equation~(\ref{MC}) in terms of the deformation theory of a
geometric structure.  Similar problems arise in the work of
Barannikov and Kontsevich~\cite{KontsevichBarannikov}.

\section{Pure spinors and integrability}\label{filtr}

Since each maximal isotropic sub-bundle $L<(T\oplus
T^*)\otimes\CC$ corresponds to a pure line sub-bundle of the spin
bundle $U<\wedge^\bullet T^*\otimes\CC$, it stands to reason that
the involutivity of $L$ corresponds to some integrability
condition on $U$.  We now determine this condition.

Recall that $L<(T\oplus T^*)\otimes\CC$ is the annihilator of $U$
under Clifford multiplication:
\begin{equation*}
L=\{X+\xi\in (T\oplus T^*)\otimes\CC\ :\ (X+\xi)\cdot U=0 \}.
\end{equation*}

The Clifford algebra $CL(V\oplus V^*),\ \dim V=m$ is a
$\ZZ/2\ZZ$-graded, $\ZZ$-filtered algebra with the following
filtered direct summands:
\begin{align*}
\RR&=CL^0<CL^2<\cdots<CL^{2m}=CL^+(V\oplus V^*)\\
V\oplus V^*&=CL^1<CL^3<\cdots<CL^{2m-1}=CL^-(V\oplus V^*),
\end{align*}
where $CL^{2k}$ is spanned by products of even numbers of not more
than $2k$ elements of $V\oplus V^*$, and $CL^{2k-1}$ is spanned by
products of odd numbers of not more than $2k-1$ elements. The
Clifford multiplication respects this graded filtration structure.

By Clifford multiplication by $U$, we obtain filtrations of the
even and odd exterior forms (here $2n$ is the real dimension of
the manifold):
\begin{align}\label{filtration}
U&=U_0<U_2<\cdots<U_{2n}=\wedge^{ev/odd} T^*\otimes\CC,\\
L^*\cdot U&=U_1<U_3<\cdots<U_{2n-1}=\wedge^{odd/ev} T^*\otimes\CC,
\end{align}
where $ev/odd$ is chosen according to the parity of $U$ itself,
and $U_k$ is defined as $CL^k\cdot U$. Note that, using the inner
product, we have the canonical isomorphism $L^*=((T\oplus
T^*)\otimes\CC)/L$, and so $U_1$ is isomorphic to $L^*\otimes
U_0$.
\begin{theorem}\label{integspinor}
The almost Dirac structure $L$ is Courant involutive if and only
if the exterior derivative $d$ satisfies
\begin{equation*}
d(C^\infty(U_0))\subset C^\infty(U_1),
\end{equation*}
i.e. $L$ is involutive if and only if, for any local
trivialization $\rho$ of $U$, there exists a section $X+\xi\in
C^\infty((T\oplus T^*)\otimes\CC)$ such that
\begin{equation}\label{integspin}
d\rho=i_X\rho+\xi\wedge\rho.
\end{equation}
Note that this condition is invariant under rescaling of $\rho$ by
a smooth function.
\end{theorem}
\begin{proof}
Let $L$ be the almost Dirac structure, and let $\rho$ be a
trivialization of $U$ over some open set.  Then we prove in the
next paragraph that
\begin{equation}\label{dspinor}
[A,B]\cdot\rho=A\cdot B\cdot d\rho,
\end{equation}
for any sections $A,B\in C^\infty(L)$.  Hence $L$ is involutive if
and only if $A\cdot B\cdot d\rho=0\ \forall A,B\in C^\infty(L)$,
which is true if and only if $d\rho$ is in $C^\infty(U_1)$, as
elements of $U_k$ are precisely those which are annihilated by
$k+1$ elements in $L$.

We now prove the identity. Since $L$ is isotropic, the Courant and
Dorfman brackets agree when restricted to $L$; we use the latter
for simplicity. Let $A=X+\xi$ and $B=Y+\eta$, so that
$i_X\rho=-\xi\wedge\rho$ and $i_Y\rho=-\eta\wedge\rho$. Then
\begin{align*}
i_{[X,Y]}\rho&=[\LL_X,i_Y]\rho\\
&=\LL_X(-\eta\wedge\rho)-i_Y(d(-\xi\wedge\rho)+i_Xd\rho)\\
&=-\LL_X\eta\wedge\rho-\eta\wedge(i_Xd\rho+d(-\xi\wedge\rho)) -
i_Y(-d\xi\wedge\rho+\xi\wedge d\rho+i_Xd\rho)\\
&=(-\LL_X\eta+i_Yd\xi)\wedge\rho -
(i_Y+\eta\wedge)(i_X+\xi\wedge)d\rho\\
&=(-\LL_X\eta+i_Yd\xi)\wedge\rho +A\cdot B\cdot d\rho,
\end{align*}
showing that
\begin{equation*}
[A,B]\cdot\rho = (A\circ B)\cdot\rho = A\cdot B\cdot d\rho,
\end{equation*}
as required.
\end{proof}

\begin{remark}
For many examples of Dirac structures, it is possible to choose
local trivializations $\rho$ of the pure spinor line which are
closed: $d\rho=0$, which obviously satisfies the above
integrability condition.  However this is not always the case; the
more general integrability condition stated above is the
appropriate one.
\end{remark}

It should be mentioned that if $L$ is integrable, then the
exterior derivative, which takes $U$ to $U_1=L^*\otimes U$, is an
example of a \emph{Lie algebroid connection} on the line bundle
$U$ with respect to the Lie algebroid $L$.  The notion of Lie
algebroid connection is developed in~\cite{Fernandes}.
\begin{defn}[Lie algebroid connection]\label{algebroidconnection}
Let $L$ be a Lie algebroid, $E$ a vector bundle and
$D:C^\infty(E)\lra C^\infty(L^*\otimes E)$ a linear operator such
that
\begin{equation*}
D(fs)=(d_Lf)\otimes s + fDs
\end{equation*}
for any $f\in C^\infty(M)$ and $s\in C^\infty(E)$. Then $D$ is
called a generalized connection, or Lie algebroid connection.

If $D$ is such a connection, it can be extended in the usual way
to a sequence $D_L:C^\infty(\wedge^k L^*\otimes E)\lra
C^\infty(\wedge^{k+1}L^*\otimes E)$ via the rule:
\begin{equation*}
D_L(\mu\otimes s)=d_L\mu\otimes s + (-1)^{|\mu|}\mu\wedge Ds,
\end{equation*}
and $D_L^2\in C^\infty(\wedge^2 L^*\otimes\End(E))$ is then the
curvature of the connection.
\end{defn}

In this way we obtain a natural connection structure on the spinor
line $U$ of any Dirac structure.  This will be particularly useful
when we consider generalized complex structures in the next
section.

Before we proceed to our final remarks concerning the Courant
bracket, we should give a simple application of the above
integrability condition to products of Dirac structures.

\begin{prop}
Let $(M_1,L_1)$, $(M_2,L_2)$ be two manifolds equipped with Dirac
structures. Then $\pi_1^*L_1\oplus \pi_2^*L_2$, where $\pi_i$ are
the canonical projections, is a Dirac structure on $M_1\times
M_2$.
\end{prop}
\begin{proof}
Choose trivializations $\rho_1,\rho_2$ for the spinor lines
$U_1,U_2$ in open sets around points $m_1\in M_1,m_2\in M_2$
respectively. Then form $\rho=\pi_1^*\rho_1 \wedge \pi_2^*\rho_2$,
defined in a neighbourhood of $(m_1,m_2)$.  The annihilator of
this spinor is the maximal isotropic space $\pi_1^*L_1\oplus
\pi_2^*L_2$, where $\pi^*$ is the pull-back of Dirac structures
defined in Equation~\ref{pullback}.  Furthermore, if we have
$d\rho_i=\alpha_i\cdot\rho_i$ for $\alpha_1,\alpha_2\in
C^\infty(T\oplus T^*)$, then we see that
\begin{equation*}
d\rho=(\alpha_1\pm\alpha_2)\cdot\rho,
\end{equation*}
showing that $\rho$ satisfies the integrability
condition~(\ref{integspin}), implying the integrability of the
Dirac structure.
\end{proof}

\section{The twisted Courant bracket}\label{twistedcour}

As was noticed by \v{S}evera, and developed by him and Weinstein
in~\cite{SeveraWeinstein}, the Courant bracket on $T\oplus T^*$
can be `twisted' by a real\footnote{None of our results depend on
$H$ being real, however since $[,]$ is a real quantity, it is
reasonable to restrict ourselves to real twistings.}, closed
3-form $H$, in the following way: given a 3-form $H$, define
another bracket $[,]_H$ on $T\oplus T^*$, by
\begin{equation*}
[X+\xi,Y+\eta]_H=[X+\xi,Y+\eta]+i_Yi_XH.
\end{equation*}
Then, defining $\Nij_H$ and $\Jac_H$ using the usual formulae but
replacing $[,]$ with $[,]_H$, one calculates that if
$A=X+\xi$,$B=Y+\eta$, and $C=Z+\zeta$, then
\begin{equation*}
\Nij_H(A,B,C)=\Nij(A,B,C)+H(X,Y,Z),
\end{equation*}
and that
\begin{equation*}
\Jac_H(A,B,C)=d(\Nij_H(A,B,C)) + i_Zi_Yi_XdH.
\end{equation*}
We conclude that $[,]_H$ defines a Courant algebroid structure on
$T\oplus T^*$ (using the same inner product and anchor) if and
only if the extraneous term vanishes, i.e. $dH=0$.

Therefore, given any closed 3-form $H$, one can study maximal
isotropic sub-bundles of $T\oplus T^*$ which are involutive with
respect to $[,]_H$: these are called \emph{twisted Dirac
structures}.  Reexamining Proposition~\ref{2form}, we obtain the
following relationship between $2$-forms and the twisted Courant
brackets.
\begin{prop}
If $b$ is a 2-form then we have
\begin{equation*}
[e^b(W),e^b(Z)]_H = e^b[W,Z]_{H+db}\ \ \ \forall\ W,Z\in
C^\infty(T\oplus T^*),
\end{equation*}
showing that $e^b$ is a symmetry of $[,]_H$ if and only if $db=0$.
These are the B-field transforms.
\end{prop}
\begin{proof} See Proposition~\ref{2form}.
\end{proof}

Note that the tangent bundle $T$ is not involutive with respect to
the bracket $[,]_H$ unless $H=0$.  If $H$ is exact, i.e. there
exists $b\in\Omega^2(M)$ such that $db=H$, then clearly the bundle
\[
e^{-b}(T)=\{X-i_Xb\ :\ X\in T\}
\]
is closed with respect to $[,]_H$.  In general, a sub-bundle $L$
is involutive for $[,]_H$ if and only if $e^{-b}L$ is involutive
for $[,]_{H+db}$.  Using this observation, one reduces the study
of $H$-Dirac structures, for $H$ exact, to the study of ordinary
Dirac structures.  For $H$ such that $[H]\in H^3(M,\RR)$ is
nonzero, $H$-Dirac strucures represent geometries genuinely
separate from ordinary Dirac structures. The main object of study
in~\cite{SeveraWeinstein} is $H$-twisted Poisson structure.
\begin{example}[Twisted Poisson geometry]
Let $H\in\Omega^3_{cl}(M)$, and let $\beta\in C^\infty(\wedge^2
T)$ be such that
\begin{equation*}
[\beta,\beta]=\beta^*(H),
\end{equation*}
where on the right hand side we are pulling back by
$\beta:T^*\rightarrow T$. Then $\beta$ is called a $H$-twisted
Poisson structure, i.e. the graph $L_{\beta}=\{i_\xi\beta+\xi\ :\
\xi\in T^* \}$ is involutive with respect to the $H$-twisted
Courant bracket. See~\cite{SeveraWeinstein} for details.
\end{example}

Modifying the integrability condition of a maximal isotropic
sub-bundle by introducing a twist must also modify the
integrability condition on the level of the spinor line
$U<\wedge^\bullet T^*\otimes\CC$.  We now determine what this is.
\begin{prop}\label{twistprop}
The almost Dirac structure $L$ is involutive with respect to the
$H$-twisted Courant bracket if and only if the operator $d_H=d +
H\wedge\cdot$ satisfies
\begin{equation*}
d_H(C^\infty(U_0))\subset C^\infty(U_1),
\end{equation*}
i.e. $L$ is $H$-involutive if and only if, for any local
trivialization $\rho$ of $U$, there exists a section $X+\xi\in
C^\infty((T\oplus T^*)\otimes\CC)$ such that
\begin{equation}
d_H\rho=i_X\rho+\xi\wedge\rho.
\end{equation}
\end{prop}
\begin{proof}
The key result of the proof of Theorem~\ref{integspinor} was that
\begin{equation*}
[A,B]\cdot\rho=A\cdot B\cdot d\rho.
\end{equation*}
From this we see immediately that
\begin{equation*}
[A,B]_H\cdot\rho = A\cdot B\cdot d\rho + i_Yi_XH\wedge\rho,
\end{equation*}
where $A=X+\xi$ and $B=Y+\eta$.  Using the fact that
$i_X\rho+\xi\wedge\rho=i_Y\rho+\eta\wedge\rho=0$, we obtain
\begin{equation*}
[A,B]_H\cdot\rho = A\cdot B\cdot (d\rho + H\wedge\rho),
\end{equation*}
which is what we require: using the same reasoning as the
untwisted case, $\rho$ determines a $H$-Dirac structure if and
only if $d_H\rho=i_X\rho+\xi\wedge\rho$ for some section $X+\xi\in
C^\infty(T\oplus T^*)$.
\end{proof}

\section{Relation to gerbes}\label{gerbes}

At this point we would like to give some indication of how gerbes
are related to Courant algebroids; in particular how the
$H$-twisted Courant bracket can be viewed as a ``twist by a
gerbe'' when $[H]$ is an integral cohomology class.  This section
is independent of the rest of the thesis and is mostly a
collection of remarks intended to demonstrate that the algebraic
structures introduced in this chapter have geometric
underpinnings. The gerbe interpretation is particularly relevant
for physicists working with sigma models, for whom $H$ is known as
the Neveu-Schwarz 3-form flux: in current theories $[H]$ is
required to be integral for the Lagrangian to be well-defined.
Among other things, the gerbe interpretation provides a reason for
the fact that $B$-field transformations with $[B]\in H^2(M,\ZZ)$
should be interpreted as gauge transformations, e.g. for questions
of moduli they should be quotiented out.

The Courant bracket on $T\oplus T^*$ is just ``level 1'' of a
hierarchy of brackets on the bundles $T\oplus \wedge^p T^*$,
$p=0,1,\ldots$, defined by the same formula
\begin{equation*}
[X+\sigma,Y+\tau]=[X,Y]+\LL_X\tau-\LL_Y\sigma
-\tfrac{1}{2}d(i_X\tau-i_Y\sigma),
\end{equation*}
where now $\sigma,\tau\in C^\infty(\wedge^p T^*)$.  In his
thesis~\cite{Roytenberg}, Roytenberg showed that for $p=1$, the
Courant bracket defines an $L_\infty$ algebra of level $2$.   The
exact nature of the full algebraic structure for all $p$ remains a
work in progress. Using identical arguments to those used for
$T\oplus T^*$, one sees that symmetries of the bracket are given
by closed $p+1$ forms, and that the bracket may be twisted by a
closed $p+2$ form.

\paragraph{Level 0:}
It will be fruitful to begin our investigation with the case
$p=0$, i.e. the bracket on $T\oplus 1$ given by
\begin{equation*}
[X+f,Y+g]=[X,Y]+Xg-Yf,
\end{equation*}
where $X,Y\in C^\infty(T)$ and $f,g\in C^\infty(M)$.  Unlike the
bracket on $T\oplus T^*$, this one actually satisfies the Jacobi
identity, and makes $T\oplus 1$ into a Lie algebroid, with anchor
the obvious projection $\pi:T\oplus 1\lra T$.

This bracket operation should be familiar: it is the Lie bracket
structure on the Atiyah sequence associated to the trivial $S^1$
principal bundle $P=S^1\times M$.  In particular, if
$X+f\partial_t,Y+g\partial_t$ are $S^1$-invariant vector fields on
$P$ (here $t$ is the coordinate on $S^1$ and $f,g$ have no
dependence on $t$), then their Lie bracket is
\begin{equation*}
[X+f\partial_t,Y+g\partial_t]=[X,Y] + (Xg-Yf)\partial_t.
\end{equation*}
In this way we achieve a geometrical interpretation of the
untwisted Courant bracket for $p=0$.  The bracket, however, may be
twisted by a closed 2-form $F$:
\begin{equation*}
[X+f,Y+g]_F=[X+f,Y+g]+i_Yi_XF,
\end{equation*}
and we would like to interpret this twisting in some geometrical
way.

The interpretation using $S^1\times M$ can be generalized since
there is an Atiyah sequence associated to any $S^1$ principal
bundle $P$:
\begin{equation*}
\xymatrix{0\ar[r]&1\ar[r]^{j}&TP/S^1\ar[r]^{\ \ \pi_*}&T\ar[r]&
0}.
\end{equation*}
By choosing a splitting of this extension, we choose an
isomorphism $TP/S^1\cong T\oplus 1$, and we may then transport the
Lie algebroid structure on $TP/S^1$ to $T\oplus 1$. As is well
known, choosing a splitting $\nabla$ for this sequence is
equivalent to choosing a connection on $P$:
\begin{equation*}
\xymatrix{0\ar[r]& 1\ar@<0.5ex>[r]^j&
TP/S^1\ar@<0.5ex>[l]^{s}\ar@<0.5ex>[r]^{\ \ \pi_*}&
T\ar@<0.5ex>[l]^{\ \ \nabla}\ar[r]& 0},
\end{equation*}
and the curvature of the connection measures the failure of
$\nabla$ to be a Lie algebra morphism:
\begin{equation*}
F^\nabla(X,Y)=s(\nabla(X),\nabla (Y)),
\end{equation*}
where $X,Y\in C^\infty(T)$.  Now let us calculate the Lie bracket
on $T\oplus 1$ induced by this choice of splitting:
\begin{align*}
[X+f,Y+g]&=\pi_*[\nabla
(X),\nabla(Y)]+s([\nabla(X),j(g)]-[\nabla(Y),j(f)]) +
s([\nabla(X),\nabla(Y)])\\
&=[X,Y]+Xg-Yf+F^\nabla(X,Y),
\end{align*}
where we have used the Lie algebroid axioms (satisfied by
$TP/S^1$) to establish the final equality.  Hence we see that from
a line bundle with connection we obtain a twisted Lie algebroid
bracket on $T\oplus 1$; in this case, the twisting 2-form
satisfies $[F^\nabla/{2\pi}]\in H^2(M,\ZZ)$. Note that choosing a
different connection $\nabla'=\nabla + A$, for $A\in
C^\infty(T^*)$, modifies the curvature, and therefore the twisting
2-form, by an exact form:
\begin{equation*}
F^{\nabla'}=F^{\nabla} + dA
\end{equation*}
So, we see that while we do obtain a geometrical interpretation of
certain twisted Lie algebroid brackets on $T\oplus 1$, it is only
those for which the twisting 2-form is \emph{integral} which can
be interpreted in terms of line bundles with connection.
\begin{prop}
The twisted Courant bracket $[,]_F$ on $T\oplus 1$ can be obtained
from a principal $S^1$ bundle with connection when $[F/2\pi]\in
H^2(M,\ZZ)$.
\end{prop}
From this point of view, the symmetries of the bracket also become
clear. Changing the splitting to $\nabla'=\nabla+A$ corresponds to
mapping $X+f\mapsto X+f-i_XA$, and this preserves the bracket as
long as the curvature is unchanged, i.e. $A$ is a closed 1-form.
This could be thought of as taking the tensor product with the
trivial bundle equipped with a flat connection $d+A$, $dA=0$.

If $[F/2\pi]$ is integral, then some of these symmetries actually
 come from gauge transformations of the underlying $S^1$ bundle, i.e. those
for which $A=f^{-1}df$, for $f$ an $S^1$-valued function.  In
other words, those symmetries for which $[A]\in H^1(M,\ZZ)$ are
gauge transformations modulo constant gauge transformations.
\begin{prop}
Let $[,]_F$ be a Courant bracket on $T\oplus 1$ derived from a
principal $S^1$ bundle with connection.  Its symmetries correspond
to tensoring with the trivial $S^1$ bundle with a flat connection
$d+A$, of which those with $[A]\in H^1(M,\ZZ)$ derive from gauge
transformations.
\end{prop}

In the case that $[F/2\pi]$ is not integral, it is of course still
possible to interpret the twisted bracket as deriving from an
exact Lie algebroid
\begin{equation*}
\xymatrix{0\ar[r]&1\ar[r]^{j}&E\ar[r]^{a}&T\ar[r]& 0}
\end{equation*}
equipped with a splitting $\nabla$. Then the symmetries of the
bracket are simply automorphisms of the exact sequence preserving
the Lie bracket structure.

However, it is possible to go further, and interpret this
non-integral case not as a line bundle, but as the following
generalization of a line bundle: a trivialization, or section, of
an $S^1$ gerbe. In particular we need a trivialization (with
connection) of an $S^1$ gerbe with flat connection. This is
directly analogous to the fact that we can interpret a closed
1-form $A$ as given by a trivialization (nonvanishing section) $s$
of a trivial line bundle with connection $(L,\nabla)$ in the
manner
\begin{equation*}
\nabla(s)=A\otimes s.
\end{equation*}
We will follow an argument similar to that used by Hitchin
in~\cite{HitSL}, and we use the \v{C}ech description of gerbes he
espouses, which was developed in the thesis of
Chatterjee~\cite{Chatterjee}.

In the same way that an $S^1$ principal bundle can be specified by
a cocycle $g_{\alpha\beta}:U_\alpha\cap U_\beta\rightarrow S^1$ in
$\check{C}^1(M,C^\infty(S^1))$, an $S^1$ gerbe can be specified by
a cocycle $g_{\alpha\beta\gamma}:U_\alpha\cap U_\beta\cap
U_\gamma\rightarrow S^1$ in $\check{C}^2(M,C^\infty(S^1))$. A
connection on a gerbe is specified by 1-forms $A_{\alpha\beta}$
and 2-forms $B_\alpha$ satisfying the following deRham-\v{C}ech
conditions:
\begin{align*}
iA_{\alpha\beta}+iA_{\beta\gamma}+iA_{\gamma\alpha}&=g^{-1}_{\alpha\beta\gamma}dg_{\alpha\beta\gamma},\\
B_\beta-B_\alpha&=dA_{\alpha\beta}
\end{align*}
and we see immediately that there is a globally-defined 3-form $H$
(with $[H/2\pi]$ integral), known as the \emph{curvature} of the
gerbe connection and defined by
\begin{equation*}
H\big|_{U_\alpha}=dB_\alpha.
\end{equation*}

A trivialization of a gerbe is given by
$h_{\alpha\beta}\in\check{C}^1(M,C^\infty(S^1))$ satisfying
\begin{equation*}
h_{\alpha\beta}h_{\beta\gamma}h_{\gamma\alpha}=g_{\alpha\beta\gamma},
\end{equation*}
and the existence of such a section implies that the gerbe is
trivial, i.e. $[g_{\alpha\beta\gamma}]=0$.

If the connection on the gerbe is flat, then since $dB_\alpha=0$
we can write (for a suitable cover refinement)
$B_\alpha=da_{\alpha}$, and then since
$dA_{\alpha\beta}=B_\beta-B_\alpha=d(a_\beta-a_\alpha)$, we obtain
functions $f_{\alpha\beta}$ such that
$A_{\alpha\beta}-a_\beta+a_\alpha = df_{\alpha\beta}$.  This then
implies that
\begin{equation*}
idf_{\alpha\beta}+idf_{\beta\gamma}+idf_{\gamma\alpha}=g^{-1}_{\alpha\beta\gamma}dg_{\alpha\beta\gamma}.
\end{equation*}

A trivialization $h_{\alpha\beta}$ of such a flat gerbe produces
immediately a cocycle
\begin{equation*}
\mathfrak{a}_{\alpha\beta}=-i(\frac{dh_{\alpha\beta}+idf_{\alpha\beta}h_{\alpha\beta}}{h_{\alpha\beta}}),
\end{equation*}
(analogous to $A$ in $\nabla(s)=iA\otimes s$) which can be used to
manufacture an extension
\begin{equation}\label{ext}
\xymatrix{0\ar[r]&1\ar[r]&E\ar[r]&T\ar[r]& 0},
\end{equation}
by gluing $(T\oplus 1)|_{U_\alpha}$ to $(T\oplus 1)|_{U_\beta}$
using the automorphism
\begin{equation*}
\begin{pmatrix}
1&0\\\mathfrak{a}_{\alpha\beta}&1
\end{pmatrix},
\end{equation*}
where $\mathfrak{a}_{\alpha\beta}$ acts on $T$ by contraction
$X\mapsto i_X\mathfrak{a}_{\alpha\beta}$.  This extension inherits
a Lie bracket since the one forms $\mathfrak{a}_{\alpha\beta}$ are
closed.

The trivialization $h_{\alpha\beta}$ is \emph{with connection}
relative to the gerbe connection if we have 1-forms $A_\alpha$
such that
\begin{equation*}
iA_\beta-iA_\alpha = h^{-1}_{\alpha\beta} dh_{\alpha\beta} +
idf_{\alpha\beta}.
\end{equation*}
This is precisely the data determining a splitting $\nabla$ of the
extension~(\ref{ext}), and $dA_\alpha$ determines the global
closed 2-form $F$. Hence we see that the twisted Courant bracket
on $T\oplus 1$ may be viewed as a trivialization with connection
of a flat gerbe. Chatterjee studies this situation (but calls a
trivialization an ``object'') and concludes~(\cite{Chatterjee},
Prop. 3.2.5) that any two such trivializations with connection are
equivalent if and only if they differ by a trivial line bundle
with flat connection, corresponding to the fact that the
symmetries of $[,]_F$ are given by closed 1-forms.

\paragraph{Level 1:}
In the same way that $[,]_F$ on $T\oplus 1$ can be understood in
terms of the Atiyah sequence of a line bundle when $[F/2\pi]$ is
integral, the twisted Courant bracket $[,]_H$ on $T\oplus T^*$ can
be understood in terms of gerbes when $[H/2\pi]$ is integral.  In
the non-integral case, one would need to pass to $2$-gerbes, which
we will not address.

Any twisted Courant bracket on $T\oplus T^*$ can be obtained by
choosing a splitting of an exact Courant algebroid, of the form
\begin{equation*}
\xymatrix{0\ar[r]&T^*\ar[r]^{j}&E\ar[r]^{\pi}&T\ar[r]& 0}.
\end{equation*}
In~\cite{Hitchin}, Hitchin demonstrates how, in the integral case,
such a split algebroid can be naturally obtained from a gerbe with
connection.  The first condition on the connection data, namely
\begin{equation*}
A_{\alpha\beta}+A_{\beta\gamma}+A_{\gamma\alpha}=g^{-1}_{\alpha\beta\gamma}dg_{\alpha\beta\gamma},
\end{equation*}
implies that $dA_{\alpha\beta}$ is a cocycle, which can be used to
produce an extension by gluing $(T\oplus T^*)|_{U_\alpha}$ to
$(T\oplus T^*)|_{U_\beta}$ using the automorphism
\begin{equation*}
\begin{pmatrix}
1&0\\dA_{\alpha\beta}&1
\end{pmatrix},
\end{equation*}
where $dA_{\alpha\beta}$ acts on $T$ by contraction $X\mapsto
i_XdA_{\alpha\beta}$.  This extension inherits a bracket structure
since the 2-forms $dA_{\alpha\beta}$ are closed.  The second
condition on the connection data, namely
\begin{equation*}
B_\beta-B_\alpha=dA_{\alpha\beta},
\end{equation*}
defines a splitting $(\nabla,s)$ of the extension as in the Lie
algebroid case, with curvature
\begin{equation*}
s([\nabla(X),\nabla(Y)])=i_Yi_X H
\end{equation*}
given by $H$, the gerbe curvature.  In this way, we see that a
gerbe connection gives rise to a generalized Atiyah sequence
together with a splitting.
\begin{prop}
The twisted Courant bracket $[,]_H$ on $T\oplus T^*$ can be
obtained from a $S^1$ gerbe with connection when $[H/2\pi]\in
H^3(M,\ZZ)$.
\end{prop}
The advantage of this point of view, as in the ``level 0'' case,
is that we gain an understanding of the symmetries of the Courant
bracket.  In the ``level 0'' case, when $[F/2\pi]$ is integral,
trivializations of a flat line bundle act as symmetries of
$[,]_F$, and the difference of two such trivializations is a gauge
transformation. In the ``level 1'' case, again when $[H/2\pi]$ is
integral, trivializations (with connection) of a flat gerbe act as
symmetries of $[,]_H$ (B-field transforms), and the difference of
two such trivializations, a line bundle with connection, acts as a
gauge transformation (integral B-fields).

In this way, we see that trivializations (with connection) of a
flat gerbe not only give rise to twisted Lie algebroid structures
on $T\oplus 1$ but also to symmetries of $T\oplus T^*$.

\clearpage
\thispagestyle{empty}
\cleardoublepage
\chapter{Generalized complex structures}\label{GC}

In this chapter we introduce our main object of study: generalized
complex structures.  The idea originates with Nigel
Hitchin~\cite{Hitchin}, and is an extension of the use of Dirac
structures and Lie algebroids to incorporate complex geometry,
and, as we shall see, many other new forms of geometry. In this
way, symplectic and complex geometry can be viewed as extremal
cases of a more general structure.

In view of recent work in mirror symmetry, it is clear that there
are deep connections between the holomorphic and symplectic
categories.  This makes the explicit unification of both
structures all the more intriguing; indeed generalized complex
structures provide a natural framework in which to discuss mirror
symmetry.  In particular, the concepts of the \emph{B-field}, the
\emph{extended deformation space} of Kontsevich~\cite{Kont}, as
well as the newly-discovered \emph{coisotropic D-branes} of
Kapustin~\cite{Kapustin}, not to mention the target space
bi-Hermitian geometry discovered by Gates, Hull and
Ro\v{c}ek~\cite{Rocek}, all have natural interpretations in terms
of generalized complex geometry.

In the first section we will introduce the algebraic nature of
generalized complex structures.  In section \ref{sec2} we will
transport the algebraic structure to a manifold and investigate
the topological implications of having a generalized complex
structure.  In section \ref{sec3}, we impose the integrability
condition, and describe how it unifies complex and symplectic
geometry. In section \ref{sec4}, we describe how a generalized
complex structure affects the differential forms on a manifold. In
section \ref{sec5}, we provide some examples of manifolds which
are generalized complex, and yet have no known complex or
symplectic structures. In section \ref{sec6}, we describe a family
of generalized complex structures interpolating between a complex
and a symplectic structure, demonstrating that in some cases, the
moduli space of complex and symplectic structures may be connected
through generalized complex structures.  In section \ref{sec7} we
prove a Darboux-type theorem describing the local form of a
generalized complex structure in a regular neighbourhood.  In the
final section of the chapter, we discuss the jumping phenomenon,
where generalized complex structures may change algebraic type
along loci in the manifold.

\section{Linear generalized complex structures}\label{sec1}
We begin by defining the notion of generalized complex structure
on a real vector space.  We will use the well known structures of
complex and symplectic geometry to guide us.

Let $V$ be a real, finite dimensional vector space. A
\emph{complex} structure on $V$ is an endomorphism $J:V\lra V$
satisfying $J^2=-1$.  By comparison, a \emph{symplectic} structure
on $V$ is a nondegenerate skew form $\omega\in \wedge^2V$.  We
may, however, view $\omega$ as a map $V\lra V^*$ via interior
product:
\begin{equation*}
\omega: v\mapsto i_v\omega,\ \ v\in V.
\end{equation*}
With this in mind, a symplectic structure on $V$  can be defined
as an isomorphism $\omega: V\lra V^*$ satisfying
$\omega^*=-\omega$.  Note that we are using an asterisk to denote
the linear dual of a space or mapping, so that $\omega^*$ maps
$(V^*)^*=V$ to $V^*$.

In attempting to include both these structures in a higher
algebraic structure, we will consider endomorphisms of the direct
sum $V\oplus V^*$.  Recall from section~\ref{alg} that $V\oplus
V^*$ may be identified with its dual space using the natural inner
product $\IP{,}$.

\begin{defn}
A generalized complex structure on $V$ is an endomorphism $\JJ$ of
the direct sum $V\oplus V^*$ which satisfies two conditions.
First, it is complex, i.e. $\JJ^2=-1$; and second, it is
symplectic, i.e. $\JJ^*=-\JJ$.
\end{defn}

\begin{prop}
Equivalently, we could define a generalized complex structure on
$V$ as a complex structure on $V\oplus V^*$ which is orthogonal in
the natural inner product.
\end{prop}
\begin{proof}
If $\JJ^2=-1$ and $\JJ^*=-\JJ$, then $\JJ^*\JJ=1$, i.e. $\JJ$ is
orthogonal.  Conversely, if $\JJ^2=-1$ and $\JJ^*\JJ=1$, then
$\JJ^*=-\JJ$.
\end{proof}
The usual complex and symplectic structures are embedded in the
notion of a generalized complex structure in the following way.
Consider the endomorphism
\begin{equation*}
\JJ_J=\left(\begin{matrix}-J&0\\0&J^*\end{matrix}\right),
\end{equation*}
where $J$ is a usual complex structure on $V$, and the matrix is
written with respect to the direct sum $V\oplus V^*$. Then we see
that $\JJ_J^2=-1$ and $\JJ_J^*=-\JJ_J$, i.e. $\JJ_J$ is a
generalized complex structure. Similarly, consider the
endomorphism
\begin{equation*}
\JJ_\omega=\left(\begin{matrix}0&-\omega^{-1}\\\omega&0\end{matrix}\right),
\end{equation*}
where $\omega$ is a usual symplectic structure.  Again, we observe
that $\JJ_\omega$ is a generalized complex structure.  Therefore
we see that diagonal and anti-diagonal generalized complex
structures correspond to complex and symplectic structures,
respectively.  Now we wish to understand any intermediate
structures which mix $V$ and $V^*$; indeed we wish to understand
the space of all generalized complex structures for $V$.

The important observation is that specifying $\JJ$ is equivalent
to specifying a maximal isotropic subspace of $(V\oplus
V^*)\otimes\CC$:
\begin{prop}
A generalized complex structure on $V$ is equivalent to the
specification of a maximal isotropic complex subspace $L<(V\oplus
V^*)\otimes\CC$ \emph{of real index zero}, i.e. such that
$L\cap\overline{L}=\{0\}$.
\end{prop}
\begin{proof}
If $\JJ$ is a generalized complex structure, then let $L$ be its
$+i$-eigenspace in $(V\oplus V^*)\otimes\CC$. Then if $x,y\in L$,
$\IP{x,y}=\IP{\JJ x,\JJ y}$ by orthogonality and $\IP{\JJ x,\JJ
y}=\IP{ix,iy}=-\IP{x,y}$, implying that $\IP{x,y}=0$. Therefore
$L$ is isotropic and half-dimensional, i.e. maximally isotropic.
Also, $\overline{L}$ is the $-i$-eigenspace of $\JJ$ and thus
$L\cap\overline{L}=\{0\}$. Conversely, given such an $L$, simply
define $\JJ$ to be multiplication by $i$ on $L$ and by $-i$ on
$\overline{L}$. This real transformation then defines a
generalized complex structure on $V\oplus V^*$.
\end{proof}
This means that studying generalized complex structures is
equivalent to studying complex maximal isotropics with real index
zero, which is the most generic possible real index.  The real
index zero condition may also be expressed in terms of the data
$(E,\eps)$ (see section~\ref{maxiso} for notation):
\begin{prop}\label{indexzeros}
The maximal isotropic $L(E,\eps)$ has real index zero if and only
if $E+\overline E=V\otimes\CC$ and $\eps$ is such that the real
skew 2-form
$\omega_\Delta=\mathrm{Im}(\eps\big|_{E\cap\overline{E}})$ is
nondegenerate on $E\cap\overline E=\Delta\otimes\CC$.
\end{prop}
\begin{proof}
Let $L$ have real index zero.  Then since $(V\oplus V^*)\otimes\CC
= L\oplus\overline L$, we see that $E+\overline E = V\otimes\CC$.
Also, if $0\neq X\in\Delta$ such that $(\eps-\overline\eps)(X)=0$
then there exists $\xi\in V\otimes\CC$ such that $X+\xi\in
L\cap\overline L$, which is a contradiction. Hence $\omega_\Delta$
is nondegenerate.

Conversely, assume $E+\overline E=V\otimes\CC$ and that
$\omega_\Delta$ is nondegenerate.  Suppose $0\neq X+\xi\in
L\cap\overline L$; then $\xi|_E=\eps(X)$ and $\xi|_{\overline E} =
\overline\eps(X)$, so that $(\eps-\overline\eps)(X)=0$, which
implies $X=0$.  But then $\xi|_E=\xi|_{\overline E}=0$, hence
$\xi=0$ as well, a contradiction.  Hence $L\cap\overline L =
\{0\}$. This completes the proof.
\end{proof}

The result~(\ref{realindex}) about the parity of the real index
indicates that generalized complex structures may only exist on
even-dimensional spaces:
\begin{prop}
The vector space $V$ admits a generalized complex structure if and
only if it is even dimensional.
\end{prop}
\begin{proof}
An even-dimensional real vector space always admits complex and
symplectic structures, and these are generalized complex
structures as we have just shown.  For the converse, equation
(\ref{realindex}) implies that the real index must be congruent to
$\dim(V)$ modulo 2, showing that generalized complex structures
exist only on vector spaces of even dimension.  We also include
the following alternative argument:

Let $\JJ$ be a generalized complex structure on $V$. Since the
natural inner product on $V\oplus V^*$ is indefinite, we can find
a null vector $x\in V\oplus V^*$, i.e. $\IP{x,x}=0$. Since $\JJ$
is an orthogonal complex structure, $\JJ x$ is also null and is
orthogonal to $x$. Therefore $x,\JJ x$ span an isotropic subspace
$N<V\oplus V^*$. We can iteratively enlarge the isotropic subspace
$N$ by adding a pair of vectors, consisting of a null vector $x'$
orthogonal to $N$ together with $\JJ x'$, to the spanning set
until $N$ becomes maximally isotropic.  Since the inner product
has split signature, $N$ will finally have dimension $\dim V$.
Thus $V$ must be even dimensional.
\end{proof}

In view of this, let $2n$ be the dimension of $V$, and let $\JJ$
be a generalized complex structure on $V\oplus V^*$.  We may now
properly describe the $G$-structure determined by a generalized
complex structure $\JJ$:
\begin{prop}
A generalized complex structure on $V\oplus V^*$, for $\dim V=2n$,
is equivalent to a reduction of structure from $O(2n,2n)$ to
$U(n,n)=O(2n,2n)\cap GL(2n,\CC)$.
\end{prop}
Any other generalized complex structure can be obtained by
conjugating $\JJ$ by an element of the orthogonal group
$O(2n,2n)$.  The stabilizer of $\JJ$ for this action is
$U(n,n)=O(2n,2n)\cap GL(2n,\CC)$. Therefore the space of
generalized complex structures on $V$ is given by the coset space
\begin{equation*}
S_\JJ\cong\frac{O(2n,2n)}{U(n,n)}.
\end{equation*}
The Lie group $O(2n,2n)$ is homotopic to $O(2n)\times O(2n)$ and
so has four connected components, while $U(n,n)\sim U(n)\times
U(n)$ has only one.  Therefore the space $S_\JJ$ has four
connected components.  We can distinguish between these by noting
that a generalized complex structure induces an orientation on
maximal positive-definite complex subspaces of $(V\oplus
V^*,\JJ)$, and another on maximal negative-definite subspaces,
yielding four possibilities, one for each component of $S_\JJ$.
These orientations pair to give an orientation on the total space
$V\oplus V^*$, which may agree or disagree with the canonical
orientation on $V\oplus V^*$. In this way we have identified a
simple $\ZZ/2\ZZ$ invariant of a generalized complex structure,
which may be calculated by taking the top power of $\JJ$, thinking
of it as a 2-form on $V\oplus V^*$.
\begin{prop}
A generalized complex structure $\JJ$ is said to have \emph{even}
or \emph{odd} parity depending on whether its induced orientation
on $V\oplus V^*$ is $\pm 1$, i.e. $\JJ$ is even if
$\frac{1}{(2n)!}\JJ^{2n}=1$ and odd if
$\frac{1}{(2n)!}\JJ^{2n}=-1$.   This is equivalent to whether the
parity of its $+i$-eigenbundle $L$ (as a maximal isotropic
subspace) is even or odd.
\end{prop}
\begin{proof}
Suppose $L$ is even, i.e. it is in the same component of the space
of maximal isotropics as $V\otimes\CC$.  Since $O(4n,\CC)$ acts
transitively on the space of maximal isotropics, there exists
$R\in SO(4n,\CC)$ such that $R(L)=V\otimes\CC$. Then $R(\overline
L)$ intersects trivially with $V\otimes\CC$, and so is the graph
of a bivector $\beta:V^*\otimes\CC\rightarrow V\otimes\CC$ so that
$e^{-\beta} R(\overline L)=V^*\otimes\CC$ and $e^{-\beta}
R(L)=V\otimes\CC$. Hence, conjugating $\JJ$ by $e^{-\beta}R\in
SO(4n,\CC)$ produces the transformation
\begin{equation*}
(e^{-\beta}R)\JJ (e^{-\beta}R)^{-1}=\begin{pmatrix} i\mathbf{1}& \\
&-i\mathbf{1}
\end{pmatrix},
\end{equation*}
written in the $V\oplus V^*$ splitting.  This transformation then
can be thought of as a 2-form $\kappa=i\sum_{k=1}^{2n}e^k\wedge
e_k$ for $e_k,e^k$ dual bases for $V$, $V^*$ respectively.
Computing, we obtain $\tfrac{1}{(2n)!}\kappa^{2n}=1$, showing that
$\kappa$, and therefore $\JJ$, induces the canonical orientation
on $V\oplus V^*$.

On the other hand, if $L$ is odd, then it is obtained from
$V\otimes\CC$ by applying an orientation-reversing element of
$O(4n,\CC)$; therefore the orientation induced by $\JJ$ is
opposite to the canonical orientation.
\end{proof}
Since a generalized complex structure is given by the maximal
isotropic $L$, then by Proposition~\ref{algtype}, it can equally
be specified by the spinor line determining it,
$U_L<\wedge^\bullet T^*\otimes \CC$, generated by
\begin{equation}\label{form}
\varphi_L = \exp(B+i\omega)\theta_1\wedge\cdots\wedge\theta_k,
\end{equation}
where $(\theta_1,\ldots,\theta_k)$ are linearly independent
complex 1-forms in $V^*\otimes\CC$ spanning
$\Ann(\pi_{V\otimes\CC}(L))$, and $B,\omega$ are the real and
imaginary parts of a complex 2-form in $\wedge^2(V^*\otimes\CC)$.
As we showed, the parity of $\varphi_L$ as a differential form is
the same as the parity of $L$ as a maximal isotropic.  Thus,
expressing the generalized complex structure as a form renders its
parity manifest.  We call the complex line $U_L$ the
\emph{canonical line} of the generalized complex structure.

The additional constraint that $L$ is of real index zero imposes
an additional constraint on the line $U_L$, namely:
\begin{theorem}\label{spincond}
Every maximal isotropic in $V\oplus V^*$ corresponds to a pure
spinor line generated by
\begin{equation*}
\varphi_L=\exp(B+i\omega)\Omega,
\end{equation*}
where $B,\omega$ are real 2-forms and
$\Omega=\theta_1\wedge\cdots\wedge\theta_k$ for some linearly
independent complex 1-forms $(\theta_1,\ldots,\theta_k)$.  The
integer $k$ is called the \emph{type} of the maximal isotropic, as
in section~\ref{maxiso}.

The maximal isotropic is of real index zero if and only if
\begin{equation}\label{nondegeneracy}
\omega^{n-k}\wedge\Omega\wedge\overline\Omega\neq 0,
\end{equation}
or in other words
\begin{itemize}
\item
$(\theta_1,\ldots,\theta_k,\overline\theta_1,\ldots,\overline\theta_k)$
are linearly independent, and
\item $\omega$ is nondegenerate when restricted to the real
$(2n-2k)$-dimensional subspace $\Delta\leq V$ defined by
$\Delta=\Ker(\Omega\wedge\overline{\Omega})$.
\end{itemize}
\end{theorem}
\begin{proof}
By Corollary~\ref{indexzero}, we see that $\dim L\cap\overline
L=0$ if and only if $\IPS{(\varphi_L,\overline\varphi_L)}\neq 0$,
yielding
\begin{align*}
0\neq\IPS{(e^{B+i\omega}\Omega,e^{B-i\omega}\overline\Omega)}&=\IPS{(e^{2i\omega}\Omega,\overline\Omega)}\\
&=\tfrac{(-1)^{2n-k}(2i)^{n-k}}{(n-k)!}\omega^{n-k}\wedge
\Omega\wedge\overline\Omega,
\end{align*}
as required.  Alternatively, note that this result is a direct
consequence of Proposition~\ref{indexzeros}, since $B+i\omega$ has
been chosen so that $i^*(B+i\omega)=\eps$, where
$i:E\hookrightarrow V\otimes\CC$ is the inclusion.
\end{proof}
\begin{remark}
While $O(4n,\CC)$ acts transitively on the space of maximal
isotropics of $V\otimes\CC$, the real group $O(2n,2n)$ preserves
the real index $r=\dim L\cap\overline L$, and this integer labels
the orbits of $O(2n,2n)$, as shown in~\cite{Trautman}.  In
particular, note that we may apply $B$-field transforms to
generalized complex structures to obtain new ones: the condition
$r=0$ is preserved since $B$ is chosen to be real.
\end{remark}

We will now describe in detail certain examples of generalized
complex structures.
\begin{example}[Symplectic type $(k=0)$] The generalized complex
structure determined by a symplectic structure
\begin{equation*}
\JJ_\omega=\left(\begin{matrix}0&-\omega^{-1}\\\omega&0\end{matrix}\right)
\end{equation*}
determines a maximal isotropic
\begin{equation*}
L=\{X-i\omega(X)\ :\ X\in V\otimes\CC\}
\end{equation*}
and a spinor line generated by
\begin{equation*}
\varphi_L=e^{i\omega},
\end{equation*}
showing that symplectic structures are always of even parity. This
generalized complex structure has type $k=0$, where we recall that
$k$ is the codimension of the projection of $L$ to $V\otimes\CC$.
With respect to the stratification on maximal isotropics according
to type, the case $k=0$ is the generic stratum.  Since a $B$-field
transformation does not affect projections to $V\otimes\CC$, it
preserves type.  So we may transform this example by a $B$-field
and obtain another generalized complex structure of type $0$:
\begin{equation*}
e^{-B}\JJ_\omega e^{B}=
\left(\begin{matrix}-\omega^{-1}B&-\omega^{-1}\\\omega+B\omega^{-1}B&B\omega^{-1}\end{matrix}\right),
\end{equation*}
\begin{equation*}
e^{-B}(L)=\{X-(B+i\omega)(X)\ :\ X\in V\otimes\CC\}
\end{equation*}
\begin{equation*}
\varphi_{e^{-B}L}=e^{B+i\omega}.
\end{equation*}
We will call this a B-symplectic structure; by~(\ref{form}) we see
that any generalized complex structure of type $k=0$ is a B-field
transform of a symplectic structure.
\end{example}

It is an important observation that for a generalized complex
structure, $(V\oplus V^*)\otimes\CC=L\oplus\overline L$, and
projecting to $V\otimes\CC$, it follows that $V\otimes\CC =
E+\overline{E}$, where $E=\pi_{V\otimes\CC}(L)$.  For this reason,
the type of the generalized complex structure cannot exceed $n$.
Our next example is of this extremal type.
\begin{example}[Complex type $(k=n)$]The generalized complex
structure determined by a complex structure
\begin{equation*}
\JJ_J=\left(\begin{matrix}-J&0\\0&J^*\end{matrix}\right)
\end{equation*}
determines a maximal isotropic
\begin{equation*}
L=V_{0,1}\oplus V^*_{1,0}
\end{equation*}
(where $V_{1,0}=\overline{V_{0,1}}$ is the $+i$-eigenspace of
$J$), as well a spinor line generated by
\begin{equation*}
\varphi_L=\Omega^{n,0},
\end{equation*}
where $\Omega^{n,0}$ is any generator of the $(n,0)$-forms for the
complex $n$-dimensional space $(V,J)$. Hence we see that a complex
structure is of even parity when $n$ is even and of odd parity
when $n$ is odd.  This example may be transformed by a B-field,
yielding the equivalent representations
\begin{equation*}
e^{-B}\JJ_Je^B=\left(\begin{matrix}-J&0\\BJ+J^*B&J^*\end{matrix}\right)
\end{equation*}
\begin{equation*}
e^{-B}(L)=\{X+\xi- i_XB\ :\ X+\xi\in V_{0,1}\oplus V^*_{1,0}\},
\end{equation*}
\begin{equation*}
\varphi_{e^{-B}L}=e^B\Omega^{n,0}.
\end{equation*}
Notice that only the $(0,2)$ component of the real 2-form $B$ has
any effect in this transformation.  According to~(\ref{form}), the
most general form of a type $n$ generalized complex structure is
given by
\begin{equation*}
\varphi=e^{B+i\omega}\Omega,
\end{equation*}
but since $\Omega\wedge\overline\Omega\neq 0$, the n-form $\Omega$
determines a complex structure for which it is of type $(n,0)$.
Using this complex structure, we see that only the $(0,2)$
component $c=(B+i\omega)^{0,2}$ is effective.  Hence, defining the
real 2-form $B'=c+\overline{c}$, we obtain $\varphi=e^{B'}\Omega$,
showing that any generalized complex structure of type $n$ is the
B-field transform of a complex structure.
\end{example}
\begin{example}[Products]
There is a natural notion of \emph{direct sum} of generalized
complex structures; that is, if $\JJ_1$ and $\JJ_2$ are
generalized complex structures on $V_1,V_2$ respectively, then
$\JJ_1\oplus \JJ_2$ is a generalized complex structure on
$V_1\oplus V_2$. This corresponds to taking direct sums
$\pi_1^*L_1\oplus \pi_2^*L_2$ of maximal isotropics (where $\pi_i$
are the natural projections onto $V_i$), or forming the wedge
product of the defining spinor lines, i.e.
$\pi_1^*\varphi_1\wedge\pi^*_2\varphi_2$. Clearly the
nondegeneracy condition~(\ref{nondegeneracy}) is preserved, and by
the properties of wedge product it is clear that parity is
additive with respect to the product of generalized complex
structures.  An example of such a product structure is the product
of a complex structure on $V$ and a symplectic structure on $W$:
the product structure can be described by the line in
$\wedge^\bullet (V\oplus W)^*$ generated by
\begin{equation*}
\varphi=e^{i\omega}\Omega,
\end{equation*}
where $\omega\in \wedge^2 V^*$ is the symplectic structure on $V$
and $\Omega\in \wedge^{n,0}W^*$ defines the complex structure on
$W$.
\end{example}
\begin{theorem}
Any linear generalized complex structure of type $k$ can be
(noncanonically) expressed as a B-field transform of the direct
sum of a complex structure of complex dimension $k$ and a
symplectic structure of real dimension $2n-2k$.
\end{theorem}
\begin{proof}
The general form of a generalized complex structure is
\begin{equation*}
\varphi_L=e^{B+i\omega}\Omega,
\end{equation*}
where the data satisfy the nondegeneracy
condition~(\ref{nondegeneracy}).  If we choose a subspace $N\leq
V$ transverse to $\Delta$, then $V=N\oplus \Delta$ and while
$\Delta$ carries a symplectic structure
$\omega_0=\omega|_{\Delta}$, $N$ inherits a complex structure
determined by $\Omega|_N$. The space $\wedge^2 V^*$ then
decomposes as
\begin{equation*}
\wedge^2V^*=\bigoplus_{p+q+r=2}\wedge^p\Delta^*\otimes\wedge^qN^*_{1,0}\otimes\wedge^rN^*_{0,1},
\end{equation*}
so that forms have tri-degree $(p,q,r)$. While $\Omega$ is purely
of type $(0,k,0)$, the complex 2-form $A=B+i\omega$ decomposes
into six components:
\begin{equation*}\begin{array}{ccc}
 A^{200}  &  &  \\
  A^{110} & A^{101} &  \\
  A^{020} & A^{011} & A^{002}
\end{array}.\end{equation*}
Only the components $A^{200},A^{101},A^{002}$ act nontrivially on
$\Omega$ in the expression $e^A\Omega$.  Hence we are free to
modify the other three components at will.  Note that
$\omega_0=-\tfrac{i}{2}(A^{200}-\overline{A^{200}})$, i.e. the
imaginary part of $A^{200}$ is precisely the symplectic structure
on $\Delta$. Therefore, define the real 2-form
\begin{equation*}
\tilde
B=\tfrac{1}{2}(A^{200}+\overline{A^{200}})+A^{101}+\overline{A^{101}}+A^{002}+\overline{A^{002}},
\end{equation*}
and observe that $e^{\tilde B +
i\omega_0}\Omega=e^{B+i\omega}\Omega$, demonstrating that
$\varphi_L=e^{\tilde B+i\omega_0}\Omega$, i.e. $\varphi_L$ is a
B-field transform of $e^{i\omega_0}\Omega$, which is a direct sum
of a symplectic and complex structure.
\end{proof}

\section{Almost structures and topological
obstructions}\label{sec2}

We now wish to transport generalized complex structures onto
manifolds.  In the case of complex or symplectic manifolds, this
involves two steps: the specification of an algebraic or `almost'
structure on the tangent bundle, as well as an integrability
condition imposed on this structure.  In the case of generalized
complex structures, our algebraic structure exists on the sum
$T\oplus T^*$ of the tangent and cotangent bundles, and our
integrability condition involves the Courant bracket.  In this
section we will describe the algebraic consequences of having a
generalized almost complex structure, as well as the topological
obstruction to its existence.

\begin{defn}
A generalized almost complex structure on a real $2n$-dimensional
manifold $M$ is given by the following equivalent data:
\begin{itemize}
\item
an almost complex structure $\JJ$ on $T\oplus T^*$ which is
orthogonal with respect to the natural inner product $\IP{,}$,
i.e. a reduction of structure for the $O(2n,2n)$-bundle $T\oplus
T^*$ to the group $U(n,n)$,
\item
a maximal isotropic sub-bundle $L<(T\oplus T^*)\otimes\CC$ of real
index zero, i.e. $L\cap\overline{L}=0$,
\item
A pure spinor line sub-bundle $U<\wedge^\bullet T^*\otimes\CC$,
called the canonical line bundle, satisfying
$\IPS{(\varphi,\overline\varphi)}\neq 0$ at each point $x\in M$
for any generator $\varphi\in U_x$
\end{itemize}
\end{defn}
The fact that $L$ is of real index zero leads to the important
fact that we have the decomposition
\begin{equation*}
(T\oplus T^*)\otimes\CC=L\oplus\overline{L},
\end{equation*}
hence we may use the inner product $\IP{,}$ to identify
$\overline{L}=L^*$.  In this way, we obtain an alternative
splitting into the sum of dual spaces
\begin{equation*}
(T\oplus T^*)\otimes\CC=L\oplus L^*.
\end{equation*}
This is particularly significant as it converts the filtration of
$\wedge^\bullet T^*\otimes\CC$ discussed in section~\ref{filtr}
into an actual grading.  That is, if $U=U_0$ is the canonical
bundle, then let $U_k=\wedge^k \overline{L}\cdot U_0$ for
$k=1,\ldots,2n$. Then we obtain an alternative grading for the
differential forms:
\begin{equation*}
\wedge^\bullet T^*\otimes\CC=U_0\oplus U_1\oplus\cdots\oplus
U_{2n},
\end{equation*}
where, depending on the parity of $\JJ$,
\begin{equation}\label{splitting}
U_0\oplus U_2\oplus \cdots \oplus
U_{2n}=\wedge^{ev/od}T^*\otimes\CC.
\end{equation}
Note also that there is a conjugation symmetry
$\overline{U_{k}}=U_{2n-k}$. Clifford multiplication by elements
in $L,\overline{L}$ is of degree $-1,+1$ respectively, in this
grading.  The bundle $U_k$ could alternatively be defined as the
eigenbundle of $\JJ$ (acting via the Spin representation) with
eigenvalue $i(n-k)$.

From a topological point of view, a generalized complex structure
is a reduction to $U(n,n)$, but this group is homotopic to its
maximal compact subgroup $U(n)\times U(n)$, and so the $U(n,n)$
structure may be further reduced to $U(n)\times U(n)$.  This
corresponds geometrically to the choice of a positive definite
sub-bundle $C_+<T\oplus T^*$ which is complex with respect to
$\JJ$.  The orthogonal complement $C_-=C_+^\bot$ is
negative-definite and also complex, and so we obtain the
orthogonal decomposition
\begin{equation*}
T\oplus T^*=C_+\oplus C_-.
\end{equation*}
Note that since $C_\pm$ are definite and $T$ is null, the
projection $\pi_T:C_\pm\rightarrow T$ is an isomorphism.  Hence we
can transport the complex structure on $C_\pm$ to $T$, obtaining
two almost complex structures $J_+, J_-$ on $T$.  Thus we see that
a generalized almost complex structure exists on a manifold if and
only if an almost complex structure does.
\begin{prop}
The obstruction to the existence of a generalized almost complex
structure is the same as that for an almost complex structure,
which itself is the same as that for a nondegenerate 2-form
(almost symplectic structure).
\end{prop}
The choice of an almost complex structure on a manifold is a
reduction of the structure group from $Gl(2n,\RR)$ to $GL(n,\CC)$.
Noone has yet described the sufficient cohomological conditions
for such a reduction to exist, (this has only been done in
dimensions $\leq 10$, see~\cite{MODULI}), however the known
necessary conditions in general dimension can be easily obtained.
\begin{prop}
The following are necessary conditions for the existence of a
(generalized) almost complex structure on a $2n$-manifold M:
\begin{itemize}
\item The odd Stiefel-Whitney classes of $T$ must be zero.
\item There must exist classes $c_i\in H^{2i}(M,\ZZ)$,
$i=0,\ldots,n$ whose mod 2 reductions are the even Stiefel-Whitney
classes of $T$.  Also, $c_n$ must be the Euler class of $T$, and
\begin{equation*}
\sum_{i=0}^{\lfloor n/2\rfloor}(-1)^ip_i=\sum_{j=0}^{n}c_j \cup
\sum_{k=0}^{n}(-1)^kc_k,
\end{equation*}
where $p_i$ are the Pontrjagin classes of $T$.
\end{itemize}
\end{prop}
\begin{proof}
The Stiefel-Whitney classes of a complex bundle are the mod $2$
reductions of the Chern classes, and so vanish in odd degree.
Therefore if $T$ is to admit a complex structure, it must have
vanishing odd Stiefel-Whitney classes.  Furthermore, the
Pontrjagin class $p_i(T)$ is equal to $(-1)^ic_{2i}(T\otimes\CC)$,
but if $T$ admits a complex structure $J$, then
$T\otimes\CC=T_{1,0}\oplus T_{0,1}$, where $T_{1,0}, T_{0,1}$ are
respectively the $+i, -i$-eigenbundles of $J$.  Therefore,
$c(T\otimes\CC)=c(T_{1,0})\cup c(T_{0,1})$, and since
$c_i(T_{0,1})=(-1)^ic_i(T_{1,0})$ we obtain finally that
\begin{equation*}
\sum_{i=0}^{\lfloor
n/2\rfloor}(-1)^ip_i=c(T\otimes\CC)=c(T_{1,0})\cup c(T_{0,1}).
\end{equation*}
Also, by definition of the Chern classes, $c_n$ is the Euler
class.  Hence we have the required result.
\end{proof}

\begin{remark}[Characteristic classes] Because of the decomposition $T\oplus T^*=C_+\oplus C_-$, a
generalized almost complex structure has two sets of Chern
classes; $\{c_k^+=c_k(C_+)\}$ and $\{c_k^-=c_k(C_-)\}$. While
$T\oplus T^*$ itself is a complex bundle and so has Chern classes,
they can be expressed in terms of $c_k^\pm$, as follows:
\begin{equation*}
c(T\oplus T^*,\JJ)=c(C_+)\cup c(C_-).
\end{equation*}
The canonical bundle $U$ also has a characteristic class $c_1(U)$,
and since by squaring the spinor (see Proposition
(\ref{squaredet}) we know that
\begin{equation*}
U\otimes U\cong\det L,
\end{equation*}
and since $L=(T\oplus T^*,\JJ)$ as complex bundles, we obtain the
fact that
\begin{equation*}
c_1(U)=\tfrac{1}{2}c_1(L)=\tfrac{1}{2}(c_1^++c_1^-).
\end{equation*}
Note here that $c_1^{\pm}$ is congruent to $w_2(C_\pm)=w_2(T)$
modulo $2$, so that $c_1^++c_1^-$ is even.
\end{remark}

We will explore the $U(n)\times U(n)$ reduction in greater detail
in chapter~\ref{GeneralizedKahler}, where we define an
integrability condition generalizing the K\"ahler condition.

\section{The Courant integrability condition}\label{sec3}

We now introduce the integrability condition on generalized almost
complex structures which interpolates between the symplectic
condition $d\omega=0$ and the complex condition that
$[T_{1,0},T_{1,0}]\subset T_{1,0}$.

\begin{defn}
The generalized almost complex structure $\JJ$ is said to be
integrable to a generalized complex structure when its
$+i$-eigenbundle $L<(T\oplus T^*)\otimes\CC$ is Courant
involutive.  In other words, a generalized complex structure is a
complex Dirac structure of real index zero.
\end{defn}
The requirement that $L$ be a complex Dirac structure and that
$L\cap\overline{L}=\{0\}$ leads us directly into the situation
described in Theorem~\ref{bial}. That is, $L,\overline{L}$ are
transverse Dirac structures of a Courant algebroid and therefore
form a \emph{Lie bialgebroid}.  Again, this means that the
exterior derivative $d_L$ on $\wedge^\bullet
L^*=\wedge^\bullet\overline{L}$ is a derivation of the Courant
bracket $[,]$ on $\wedge^\bullet\overline{L}$.

Furthermore, because of the decomposition $(T\oplus
T^*)\otimes\CC=L\oplus\overline{L}$, the projection $E=\pi_T(L)$
satisfies $E+\overline{E}=T\otimes\CC$, and we are in the
situation of Proposition~\ref{cxliedist}, which concludes that $E$
gives rise to a smooth integrable distribution $\Delta$, defined
by $\Delta\otimes\CC=E\cap\overline{E}$.  Recall that a point at
which $\dim\Delta$ is locally constant is called a regular point,
and from Proposition~\ref{transversecx} we conclude that near a
regular point of a generalized complex structure we obtain a
transverse complex structure to the foliation defined by $\Delta$.
The \emph{type}, $k\in\{0,\ldots,n\}$, of the generalized complex
structure at $x\in M$ is defined as the codimension of
$E_x<T_x\otimes\CC$, and therefore the leaves of the induced
foliation have dimension $\dim_{\RR}\Delta = 2n-2k$.  These leaves
inherit symplectic structure, as follows:  in the regular
neighbourhood, the complex Dirac structure $L$ may be expressed,
as in Proposition~\ref{indexzeros}, as $L(E,\eps)$, where
$E<T\otimes\CC$ is a sub-bundle and $\eps\in C^\infty(\wedge^2
E^*)$, such that $E+\overline E=T\otimes\CC$ and
$\omega_\Delta=\mathrm{Im}(\eps\big|_{E\cap\overline{E}})$ is a
nondegenerate real 2-form on $\Delta$.  The integrability of
$\omega_\Delta$ follows from the Courant involutivity of
$L(E,\eps)$:
\begin{prop}\label{regular}
Let $E<T\otimes\CC$ be a sub-bundle and $\eps\in C^\infty(\wedge^2
E^*)$.  Then the maximal isotropic $L(E,\eps)$ defines an
integrable generalized complex structure if and only if $E$ is
involutive and $d_E\eps=0$.
\end{prop}
\begin{proof}
Let $i:E\hookrightarrow T\otimes\CC$ be the inclusion. Then
$d_E:C^\infty(\wedge^k E^*)\rightarrow C^\infty(\wedge^{k+1}E^*)$
is defined by $i^*\circ d = d_E\circ i^*$.  Now let $\sigma\in
C^\infty(\wedge^2 T^*\otimes\CC)$ be a smooth extension of $\eps$,
i.e. $i^*\sigma=\eps$. Suppose that $X+\xi,Y+\eta\in C^\infty(L)$,
which means that $\xi|_{E}=i_X\eps$ and $\eta|_{E}=i_Y\eps$.
Consider the bracket $Z+\zeta=[X+\xi,Y+\eta]$; if $L$ is Courant
involutive, then $Z\in C^\infty(E)$, showing $E$ is involutive,
and the difference
\begin{align*}
\zeta|_{E}-i_Z\eps &=
i^*(\LL_X\eta-\LL_Y\xi-\tfrac{1}{2}d(i_X\eta-i_Y\xi)) - i_{[X,Y]}i^*\sigma\\
&=i_Xd_Ei^*\eta - i_Yd_Ei^*\xi
+\tfrac{1}{2}d_E(i_Xi_Y\eps-i_Yi_X\eps) - i^*[\LL_X,i_Y]\sigma\\
&=i_Xd_Ei^*\eta - i_Yd_Ei^*\xi +d_Ei_Xi_Y\eps - i^*(i_Xdi_Y\sigma + di_Xi_Y\sigma - i_Ydi_X\sigma-i_Yi_Xd\sigma)\\
&=i_Yi_Xd_E\eps
\end{align*}
must vanish for all $X+\xi,Y+\eta\in C^\infty(L)$, showing that
$d_E\eps=0$. Reversing the argument we see that the converse holds
as well.
\end{proof}
A corollary to this result is that the nondegenerate 2-form
$\omega_\Delta\in C^\infty(\wedge^2\Delta^*)$ is closed along the
leaves, showing that in a regular neighbourhood, a generalized
complex structure gives rise to a foliation with symplectic leaves
and a transverse complex structure.

We may now verify that the integrability condition on generalized
almost complex structures yields the classical conditions on
symplectic and complex structures. Recall that by \emph{type} we
mean the codimension of $E=\pi_T(L)$ in $T\otimes\CC$. As we are
projecting a bundle, this codimension may not be constant
throughout the manifold, a fact we investigate in
section~\ref{sec8}. In this section we consider only generalized
complex structures of constant type.
\begin{example}[Symplectic type $(k=0)$] The generalized almost complex
structure determined by a symplectic structure
\begin{equation*}
\JJ_\omega=\left(\begin{matrix}0&-\omega^{-1}\\\omega&0\end{matrix}\right)
\end{equation*}
has $+i$-eigenbundle
\begin{equation*}
L=\{X-i\omega(X)\ :\ X\in T\otimes\CC\},
\end{equation*}
which is Courant involutive if and only if $d\omega=0$.  Of course
we may apply a B-field transform ($B$ a real closed 2-form) to
$\JJ_\omega$, obtaining what we will call a B-symplectic
generalized complex structure. In fact, any generalized almost
complex structure which is everywhere of type $k=0$ must be of
this form: its $+i$-eigenbundle can be expressed as $L(E,\eps)$,
where $E=T\otimes\CC$ and $\eps=-B-i\omega$ is a complex 2-form
with $\omega$ non-degenerate. The maximal isotropic is Courant
involutive if and only if $d(B+i\omega)=0$.  Hence every
generalized complex structure of type zero is B-symplectic.
\end{example}

\begin{example}[Complex type $(k=n)$]The generalized complex
structure determined by a complex structure
\begin{equation*}
\JJ_J=\left(\begin{matrix}-J&0\\0&J^*\end{matrix}\right)
\end{equation*}
has maximal isotropic eigenbundle
\begin{equation*}
L=T_{0,1}\oplus T^*_{1,0},
\end{equation*}
which, as we saw in Example~\ref{cxdirac}, is Courant involutive
if and only if $J$ is integrable as a complex structure.

The general form for the $+i$-eigenbundle of a generalized almost
complex structure $\JJ$ of constant type $k=n$ is $L(E,\eps)$,
where $E\cap\overline{E}=\{0\}$ and $\eps\in \wedge^2 E^*$. In
this case $E$ determines an almost complex structure $J$ on $T$
(i.e. $E=T_{0,1}$) and $\eps\in C^\infty(\wedge^{2,0} T^*)$.  The
involutivity of $L(E,\eps)$ is equivalent to the condition that
$E=T_{0,1}$ is involutive and $\del\eps=0$. Hence we see that a
generalized complex structure of constant type $k=n$ must be the
conjugation of a bare complex structure by a $\del$-closed
$(2,0)$-form, i.e.
\begin{equation*}
\JJ=e^{-\eps}\JJ_Je^\eps=\left(\begin{matrix}-J&0\\\eps
J+J^*\eps&J^*\end{matrix}\right).
\end{equation*}
Note that $\eps$ is not necessarily closed and is not real, and so
$\JJ$ is not, in general, a B-field transform of a complex
structure. However, if $[\eps]=0$ in $H^{2,0}_{\del}(M)$, then we
can find a $(1,0)$-form $\varphi$ such that $\eps=\del\varphi$ and
then $\JJ$ is the B-field transform of $\JJ_J$, where
\begin{equation*}
B=\eps+\overline\eps+\delbar\varphi+\del\overline\varphi.
\end{equation*}
Alternatively, if the complex structure $J$ satisfies the
$\del\delbar$-lemma, we may again express $\eps$ as the
$(2,0)$-part of a real, closed 2-form $B$.
\end{example}
Summarizing, we obtain the following result:
\begin{prop}
On a $2n$-dimensional manifold, a generalized complex structure of
type zero is a B-symplectic structure, while a generalized complex
structure of type $n$ is the transform of a complex structure by a
$\del$-closed $(2,0)$-form.
\end{prop}

\section{Integrability and differential forms}\label{sec4}

The presence of a generalized complex structure on a manifold has
implications for complex differential forms, analogous to the fact
that a complex structure induces a $(p,q)$-decomposition of forms
and a splitting $d=\del+\delbar$.

Let $\JJ$ be a generalized complex structure, and let
$U<\wedge^\bullet T^*$ be the canonical line bundle of the
generalized complex structure. We have already observed that $\JJ$
determines an alternative grading for the differential forms
\begin{equation*}
\wedge^\bullet T^*\otimes\CC=U_0\oplus U_1\oplus\cdots\oplus
U_{2n},
\end{equation*}
where $U_0=U$ and $U_k=\wedge^k \overline{L}\cdot U_0$ for
$k=1,\ldots,2n$.  We now demonstrate that the integrability of
$\JJ$ is equivalent to the fact that the exterior derivative $d$
splits into the sum $d=\del+\delbar$ where for each
$k=0,\ldots,2n-1$,
\begin{equation*}
\xymatrix{C^\infty(U_k)\ar@<0.5ex>[r]^\delbar&C^\infty(U_{k+1})\ar@<0.5ex>[l]^\del}.
\end{equation*}

\begin{theorem}\label{deldelbar}
Let $\JJ$ be a generalized almost complex structure, and define
\begin{align*}
\delbar&=\pi_{k+1}\circ d:C^\infty(U_k)\lra C^\infty
(U_{k+1})\\
\del&=\pi_{k-1}\circ d:C^\infty(U_k)\lra C^\infty (U_{k-1}),
\end{align*}
where $\pi_k$ is the projection onto $U_k$, and $U_k=\{0\}$ for
$k<0$ and $k>2n$. Then $\JJ$ is integrable if and only if
$d=\del+\delbar$.
\end{theorem}
To prove this theorem we need a generalization of
formula~(\ref{dspinor}):
\begin{lemma}
For any differential form $\rho$ and any sections $A,B\in
C^\infty(T\oplus T^*)$, we have the following identity
\begin{equation}\label{stepwise}
A\cdot B\cdot d\rho = d(B\cdot A\cdot\rho) + B\cdot
d(A\cdot\rho)-A\cdot d(B\cdot\rho) + [A,B]\cdot\rho -
d\IP{A,B}\wedge\rho.
\end{equation}
\end{lemma}
\begin{proof}
Let $A=X+\xi, B=Y+\eta$.  then
\begin{align*}
A\cdot B\cdot d\rho&= (i_X+\xi\wedge)(i_Y+\eta\wedge)d\rho\\
&=i_Xi_Yd\rho+i_X\eta\wedge d\rho-\eta\wedge i_Xd\rho + \xi\wedge
i_Yd\rho + \xi\wedge\eta\wedge d\rho\\
&=di_Yi_X\rho + i_Ydi_X\rho-i_Xdi_Y\rho+i_{[X,Y]}\rho +
i_X\eta\wedge d\rho-\eta\wedge i_Xd\rho+\xi\wedge
i_Yd\rho+\xi\wedge\eta\wedge d\rho\\
&=d((i_Y+\eta\wedge)(i_X+\xi\wedge)\rho) +
(i_Y+\eta\wedge)d((i_X+\xi\wedge)\rho) -
(i_X+\xi\wedge)d((i_Y+\eta\wedge)\rho) \\
&\ \ \ \ \ +[X+\xi,Y+\eta]\rho-\tfrac{1}{2}d(i_Y\xi+i_X\eta)\\
&=d(B\cdot A\cdot\rho) + B\cdot d(A\cdot\rho)-A\cdot d(B\cdot\rho)
+ [A,B]\cdot\rho - d\IP{A,B}\wedge\rho.
\end{align*}
\end{proof}
\begin{proof}[Proof of Theorem \ref{deldelbar}]
By induction; recall that Clifford multiplication by
$L,\overline{L}$ is of degree $-1, +1$ respectively in the
alternative grading.  First let $\rho\in C^\infty(U_0)$. Then for
any $A,B\in C^\infty(L)$, equation~(\ref{stepwise}) implies that
$A\cdot B\cdot d\rho=[A,B]\cdot\rho$.  This shows that
$d(C^\infty(U_0))\subset C^\infty(U_1)$ if and only if $L$ is
Courant involutive, since $d$ is degree 1 in the usual grading of
forms and hence $d\rho$ can have no $U_0$ component (see
equation~(\ref{splitting})). Now assume that $d=\del+\delbar$ for
all $U_i$ such that $0\leq i<k$. Then take $\rho\in C^\infty(U_k)$
and for any $A,B\in C^\infty(L)$ we apply (\ref{stepwise}) to
obtain
\begin{equation*}
A\cdot B\cdot d\rho = d(B\cdot A\cdot\rho) + B\cdot
d(A\cdot\rho)-A\cdot d(B\cdot\rho) + [A,B]\cdot\rho.
\end{equation*}
On the right hand side, by induction each term is in
$C^\infty(U_{k-3}\oplus U_{k-1})$; therefore $d\rho$ is in
$C^\infty(U_{k-1}\oplus U_{k+1})$, again using the fact that $d$
is of degree 1 in the usual grading of forms, ensuring that
$d\rho$ has no $U_k$ component.  This completes the proof.
\end{proof}
\begin{example}
In the complex case, $U_0=\wedge^{n,0} T^*$ and
\begin{equation*}
U_k=\bigoplus_{p}\wedge^{n-p,k-p}T^*,
\end{equation*}
so that $\del$ and $\delbar$ are the usual operators on
differential forms for a complex manifold.
\end{example}

Just as in the complex case, once we have the decomposition
$d=\del+\delbar$ according to the grading, the fact that $d^2=0$
implies that $\del^2={\delbar}^2=0$, and
$\del\delbar=-\delbar\del$. Furthermore we can form another real
operator $d^{\JJ}=i(\delbar-\del)$, which can also be written
$d^{\JJ}=[d,\JJ]$, and which satisfies $(d^{\JJ})^2=0$.

\begin{remark}It is interesting to note that while in the complex case
$d^{\JJ}$ is just the usual $d^c$-operator $d^c=i(\delbar-\del)$,
in the symplectic case $d^\JJ$ is equal to the symplectic adjoint
of $d$ defined by Koszul~\cite{Koszul} and studied by
Brylinski~\cite{Brylinski} in the context of symplectic harmonic
forms. \end{remark}

The natural differential operator $\delbar$ can be viewed as an
operator
\begin{equation*}
\delbar: C^\infty(\wedge^k L^*\otimes U)\rightarrow
C^\infty(\wedge^{k+1}L^*\otimes U),
\end{equation*}
extended from $d:C^\infty(U)\rightarrow C^\infty(L^*\otimes U)$
via the rule
\begin{equation*}
\delbar(\mu\otimes s)=d_L\mu\otimes s + (-1)^{|\mu|}\mu\wedge ds,
\end{equation*}
for $\mu\in C^\infty(\wedge^k L^*)$ and $s\in C^\infty(U)$. As
such, $\delbar$ is an example of a Lie algebroid connection on
$U$, in the sense of Definition~\ref{algebroidconnection}.  Since
$\delbar^2=0$, this Lie algebroid connection has vanishing
curvature, and because of this, we say that $(U,\delbar)$ is a
\emph{generalized holomorphic bundle}.
\begin{defn}
Let $E$ be a complex vector bundle on a generalized complex
manifold with $+i$-eigenbundle $L$.  Then the data
$(E,\delbar_E)$, where $\delbar_E$ is a Lie algebroid connection
on $E$ with respect to $L$, is said to be a generalized
holomorphic bundle if and only if $\delbar_E^2=0$.
\end{defn}

Of course, the trivial complex line bundle is always generalized
holomorphic, using the operator $d_L:C^\infty(\wedge^k
L^*)\rightarrow C^\infty(\wedge^{k+1} L^*)$ as the Lie algebroid
connection.

In special cases, the canonical line bundle may be
\emph{holomorphically trivial}, in the sense that $(U,\delbar)$ is
isomorphic to the trivial bundle together with its Lie algebroid
connection $d_L$. This is equivalent to the existence of a
nowhere-vanishing section $\rho\in C^\infty(U)$ satisfying
$d\rho=0$.  In~\cite{Hitchin}, Hitchin calls these generalized
Calabi-Yau structures:
\begin{defn}
A generalized Calabi-Yau structure is a generalized complex
structure with holomorphically trivial canonical bundle, i.e.
there exists a nowhere-vanishing closed section $\rho\in
C^\infty(U)$.
\end{defn}

\section{Exotic examples of generalized complex
structures}\label{sec5}

In this section we describe examples of six-dimensional manifolds
which admit no known complex (type 0) or symplectic (type 3)
structures, and yet do admit generalized complex structures of
types 1 or 2.  These manifolds are \emph{nilmanifolds}; a more
extensive survey of generalized complex structures on
6-nilmanifolds has been carried out in collaboration with
Cavalcanti~\cite{MarcoGil}.

A nilmanifold is a homogeneous space $M=G/\Gamma$, where $G$ is a
simply-connected nilpotent real Lie group and $\Gamma$ is a
lattice of maximal rank in $G$. The simplest nilmanifold is the
torus $\RR^k/\ZZ^k$.  The differential graded algebra of
left-invariant forms on $G$ is quasi-isomorphic to the de Rham
complex of $M$, and serves as a rational minimal model for the
nilmanifold. Hence, up to equivalence by common finite covers, the
nilmanifolds can be fully distinguished simply by giving the
differentials of a set $\{e_1,\ldots, e_6\}$ of linearly
independent left-invariant 1-forms.  In the nilmanifold literature
this information is usually presented as in the following example:
the array $(0,0,0,12,13,14+35)$ describes a nilmanifold with de
Rham complex generated by 1-forms $e_1,\ldots e_6$ and such that
$de_1=de_2=de_3=0$, while $de_4=e_1\wedge e_2, de_5=e_1\wedge
e_3$, and $de_6=e_1\wedge e_4+e_3\wedge e_5$.

In six dimensions, there are 34 isomorphism classes of real
nilpotent Lie algebras (see~\cite{Salamon} for a detailed list and
literature review). For each of the 34 algebras there may be many
nilmanifolds, and these are distinguished by their fundamental
group.

Salamon studied the problem in \cite{Salamon} of which
6-nilmanifolds carry left-invariant complex or symplectic
structures.  Of course this reduces to determining which of the 34
nilpotent lie algebras admit invariant structures.  His results
are as follows: exactly 18 of the 34 admit complex structures,
exactly 26 admit symplectic structures, and 15 admit both, of
which the 6-torus is the only K\"ahler example. This leaves 5
classes of nilmanifold which admit no known complex or symplectic
structure. They are (as listed in~\cite{Salamon}):
\begin{itemize}
\item (0,0,12,13,14+23, 34+52)
\item (0,0,12,13,14,34+52)
\item (0,0,0,12,13,14+35)
\item (0,0,0,12,23,14+35)
\item (0,0,0,0,12,15+34)
\end{itemize}
In fact these 5 families admit generalized complex structures,
which we will exhibit by listing the pure spinor lines defining
them.  This is greatly simplified since we have an explicit handle
on the de Rham complex.  The situation is even more special since
in each case, the canonical bundle is holomorphically trivial: the
following examples are all generalized Calabi-Yau.
\begin{example}[0,0,12,13,14+23, 34+52]
Let $\rho = e^{B+i\omega}\Omega$, where
\begin{align*}
\Omega&=e_1+ie_2\\
B&=e_2\wedge e_6 -e_3\wedge e_5 + e_3\wedge e_6 - e_4\wedge e_5\\
\omega &= e_3\wedge e_6 + e_4\wedge e_5.
\end{align*}
As defined, $\rho$ is a closed pure spinor of real index zero, and
defines a generalized complex structure of type $k=1$.
\end{example}
\begin{example}[0,0,12,13,14,34+52]
Let $\rho = e^{B+i\omega}\Omega$, where
\begin{align*}
\Omega&=e_1+ie_2\\
B&=e_3\wedge e_6 - e_4\wedge e_5\\
\omega &= e_3\wedge e_6 + e_4\wedge e_5.
\end{align*}
As defined, $\rho$ is a closed pure spinor of real index zero, and
defines a generalized complex structure of type $k=1$.
\end{example}
\begin{example}[0,0,0,12,13,14+35]
Let $\rho = e^{B+i\omega}\Omega$, where
\begin{align*}
\Omega&=e_1+ie_2\\
B&=0\\
\omega &= e_3\wedge e_6 + e_4\wedge e_5.
\end{align*}
As defined, $\rho$ is a closed pure spinor of real index zero, and
defines a generalized complex structure of type $k=1$.
\end{example}
\begin{example}[0,0,0,12,23,14+35]
Let $\rho = e^{B+i\omega}\Omega$, where
\begin{align*}
\Omega&=e_1+ie_2\\
B&=-e_3\wedge e_6+e_4\wedge e_5\\
\omega &= e_3\wedge e_6 + e_4\wedge e_5.
\end{align*}
As defined, $\rho$ is a closed pure spinor of real index zero, and
defines a generalized complex structure of type $k=1$.
\end{example}
\begin{example}[0,0,0,0,12,15+34]
Let $\rho = e^{B+i\omega}\Omega$, where
\begin{align*}
\Omega&=(e_1+ie_2)\wedge(e_3+ie_4)\\
B&=0\\
\omega &= e_5\wedge e_6.
\end{align*}
As defined, $\rho$ is a closed pure spinor of real index zero, and
defines a generalized complex structure of type $k=2$.
\end{example}

Our purpose here is simply to demonstrate that there are examples
of generalized complex structures which are not simply $B$-field
transforms of products of complex and symplectic manifolds. A
fuller exploration of these structures and their moduli on
nilmanifolds will appear in our work with
Cavalcanti~\cite{MarcoGil}.

\section{Interpolation between complex and symplectic
structures}\label{sec6}

We learned from the case of linear generalized complex structures
that complex and symplectic structures have opposite parity in
dimensions $4k+2$ and the same parity in dimension $4k$.  We now
show that it is possible to interpolate smoothly between a complex
structure and a symplectic structure \emph{through} integrable
generalized complex structures when $M$ is a hyperk\"ahler
manifold, e.g. a K3 surface.

If $M$ is a K\"ahler manifold then it is equipped with an
integrable complex structure $J$ and a symplectic form $\omega$
which is of type $(1,1)$, which means that $\omega J=-J^*\omega$,
which implies that
\begin{equation*}
\left(\begin{matrix}-J&0\\0&J^*\end{matrix}\right)\left(\begin{matrix}0&-\omega^{-1}\\\omega&0\end{matrix}\right)=\left(\begin{matrix}0&-\omega^{-1}\\\omega&0\end{matrix}\right)\left(\begin{matrix}-J&0\\0&J^*\end{matrix}\right),
\end{equation*}
or in other words, the generalized complex structures commute.

On the other hand, suppose that $M$ is hyperk\"ahler, which means
it has a triple of K\"ahler complex structures $I,J,K$ with
coincident K\"ahler metric satisfying the quaternionic relations
$IJ=K=-JI$. Then we have the K\"ahler forms
$\omega_I,\omega_J,\omega_K$ and as before $\omega_J J
=-J^*\omega_J$.  However, since $IJ=-JI$, we see that $\omega_J I
=  I^*\omega_J$, which implies that the generalized complex
structures $\JJ_{\omega_J}$ and $\JJ_{I}$ \emph{anticommute}.
Hence form the one-parameter family of generalized almost complex
structures
\begin{equation*}
\JJ_t=\sin t \JJ_I + \cos t \JJ_{\omega_J},\ \
t\in[0,\tfrac{\pi}{2}].
\end{equation*}
Clearly $\JJ_t$ is a generalized almost complex structure; we now
check that it is integrable.
\begin{prop}
Let $M$ be a hyperk\"ahler manifold as above.  Then the
generalized complex structure $\JJ_t=\sin t \JJ_I + \cos t
\JJ_{\omega_J}$ is integrable $\forall t\in[0,\tfrac{\pi}{2}]$.
Therefore it is a family of generalized complex structures
interpolating between a symplectic structure and a complex
structure.
\end{prop}
\begin{proof}
Let $B=\tan t\ \omega_K$, a closed 2-form which is well defined
$\forall\ t\in[0,\tfrac{\pi}{2})$.  Noting that $\omega_K I =
I^*\omega_K=\omega_J$, we obtain the following expression:
\begin{equation*}
e^B\JJ_te^{-B}=\left(\begin{matrix}0&-(\sec t\  \omega_J)^{-1}\\
\sec t\ \omega_J&0\end{matrix}\right).
\end{equation*}
We conclude from this that for all $t\in[0,\tfrac{\pi}{2})$,
$\JJ_t$ is in fact the B-field transform of the symplectic
structure determined by $\sec t\ \omega_J$, with $B=\tan t\
\omega_K$, and is therefore integrable as a generalized complex
structure; at $t=\tfrac{\pi}{2}$, $\JJ_t$ is purely complex, and
is integrable as well, completing the proof.  Note that this
interpolation argument applies to holomorphic symplectic manifolds
as well.
\end{proof}

\section{Local structure: The generalized Darboux theorem}\label{sec7}

The Newlander-Nirenberg theorem tells us that an integrable
complex structure on a $2n$-manifold is locally equivalent, via a
diffeomorphism, to $\CC^n$.  Similarly, the Darboux theorem states
that a symplectic structure on a $2n$-manifold is locally
equivalent, via a diffeomorphism, to the standard symplectic
structure $(\RR^{2n},\omega)$, where
\begin{equation*}
\omega=dx_1\wedge dx_2+\cdots+dx_{2n-1}\wedge dx_n.
\end{equation*}
In this section we prove the analogous theorem in the generalized
context, which states that at a regular point of type $k$, a
generalized complex structure on a $2n$-manifold is locally
equivalent, via a diffeomorphism and a B-field, to the standard
product generalized complex structure
$\CC^k\times(\RR^{2n-2k},\omega)$.

We saw in Proposition~\ref{regular} that in a regular
neighbourhood, a generalized complex structure may be expressed as
$L(E,\eps)$ where $E<T\otimes\CC$ is an involutive sub-bundle and
$\eps\in C^\infty(\wedge^2 E^*)$ satisfies $d_E\eps=0$.  By
Proposition~\ref{transversecx}, the distribution $E$ determines a
foliation of the neighbourhood with transverse complex structure
isomorphic to an open set in $\RR^{2n-2k}\times\CC^k$, where $E$
is spanned by $\{\partial/{\partial x_1},\ldots,\partial/{\partial
x_{2n-2k}},\partial/{\partial z_1},\ldots,\partial/{\partial
z_{k}} \}$, where $\{x_i\}$ are coordinates for the leaves
$\RR^{2n-2k}$ and $\{z_i\}$ are transverse complex coordinates.
Therefore, by choosing $B+i\omega\in C^\infty(\wedge^2
T^*\otimes\CC)$ such that $i^*(B+i\omega)=\eps$, we may write a
generator for the canonical bundle defining $L(E,\eps)$ as
follows:
\begin{equation*}
\rho=e^{B+i\omega}\Omega,
\end{equation*}
where $\Omega=dz_1\wedge\cdots\wedge dz_{k}$; note $\rho$ is
independent of the choice of extension for $\eps$.  Furthermore,
we see that
\begin{equation*}
i^*d(B+i\omega)=d_Ei^*(B+i\omega)=d_E\eps = 0,
\end{equation*}
which means that $d(B+i\omega)\in\Ann^\bullet E$, implying finally
that
\begin{equation*}
d\rho = e^{B+i\omega}d(B+i\omega)\wedge\Omega = 0.
\end{equation*}

We have shown that in a regular neighbourhood, any generalized
complex structure on a $2n$-dimensional manifold may be expressed
as a closed complex differential form $\rho=e^{B+i\omega}\Omega$,
where $\Omega$ is decomposable of degree $0\leq k\leq n$ and such
that
\begin{equation*}
\omega^{n-k}\wedge\Omega\wedge\overline\Omega \neq 0.
\end{equation*}

Weinstein's proof of the Darboux normal coordinate theorem for a
family of symplectic structures (see~\cite{McDuff}) can be used to
find a leaf-preserving local diffeomorphism $\varphi$ taking
$\omega$ to a 2-form whose pullback to each leaf is the standard
Darboux symplectic form on $\RR^{2n-2k}$, i.e.
\begin{equation*}
\varphi^*\omega\big|_{\RR^{2n-2k}\times\{pt\}}=\omega_0=dx_1\wedge
dx_2+\cdots+dx_{2n-2k-1}\wedge dx_{2n-2k}.
\end{equation*}
Let us apply this diffeomorphism, obtaining new 2-forms
$\varphi^*B+i\varphi^*\omega$.  Note that $\Omega$ is unaffected
by this diffeomorphism, since $\{z_i\}$ are constant along the
leaves.

For convenience, let $K=\RR^{2n-2k}$ and $N=\CC^k$, so that
differential forms now have tri-degree $(p,q,r)$ for components in
$\wedge^pK^*\otimes\wedge^qN^*_{1,0}\otimes\wedge^rN^*_{0,1}$.
Furthermore, the exterior derivative decomposes into a sum of
three operators
\begin{equation*}
d=d_f+\del+\delbar,
\end{equation*}
each of degree 1 in the respective component of the tri-grading.
Note that $d_f$ is the leafwise exterior derivative. While
$\Omega$ is purely of type $(0,k,0)$, the complex 2-form
$A=\varphi^*B+i\varphi^*\omega$ decomposes into six components:
\begin{equation*}\begin{array}{ccc}
 A^{200}  &  &  \\
  A^{110} & A^{101} &  \\
  A^{020} & A^{011} & A^{002}
\end{array}\end{equation*}
Note that only the components $A^{200},A^{101},A^{002}$ act
nontrivially on $\Omega$ in the expression $e^A\Omega$.  Hence we
are free to modify the other three components at will.  Also, note
that the imaginary part of $A^{200}$ is simply $\omega_0$, so that
$d(A^{200}-\overline{A^{200}})=0$, since $\omega_0$ is in constant
Darboux form.

From the condition $d(B+i\omega)\wedge\Omega=0$, we obtain the
following four equations:
\begin{align}
\delbar A^{002}&=0\label{first}\\
\delbar A^{101}+d_fA^{002}&=0\label{second}\\
\delbar A^{200}+d_fA^{101}&=0\label{third}\\
d_f A^{200}&=0.\label{fourth}
\end{align}
The last equation simply states that the pull-back of $B+i\omega$
to any leaf is closed in the leaf, as we know already.

We will now endeavour to modify $A$ so that
$\varphi^*\rho=e^A\Omega$ is unchanged but $A$ is replaced with
$\tilde{A}=\tilde{B}+\tfrac{1}{2}(A^{200}-\overline{A^{200}})$,
where $\tilde{B}$ is a real closed 2-form.  This would demonstrate
that
\begin{equation*}
\varphi^*\rho=e^{\tilde{B}+i\omega_0}\Omega,
\end{equation*}
i.e. $\rho$ is equivalent, via the composition of a B-field
transformation and a diffeomorphism, to the product of a
symplectic with a complex structure.

In order to preserve $\varphi^*\rho$, the most general form for
$\tilde B$ is
\begin{equation*}
\tilde B =
\tfrac{1}{2}(A^{200}+\overline{A^{200}})+A^{101}+\overline{A^{101}}+A^{002}+\overline{A^{002}}
+ C,
\end{equation*}
where $C$ is a real 2-form of type $(011)$.  Then clearly
$\varphi^*\rho = e^{\tilde{B}+i\omega_0}\Omega$.  Requiring that
$d\tilde B=0$ imposes two constraint equations:
\begin{align}
(d\tilde B)^{012}&=\del A^{002}+\delbar C = 0.\label{cond1}\\
(d\tilde B)^{111}&=\del A^{101}+\overline{\del A^{101}} +
d_fC=0\label{cond2}
\end{align}
The question then becomes whether we can find a real $(011)$-form
$C$ such that these equations are satisfied.  The following are
all local arguments, making repeated use of the Dolbeault lemma.
\begin{itemize}
\item From equation (\ref{first}) we obtain that
$A^{002}=\delbar\alpha$ for some $(001)$-form $\alpha$.  Then
condition (\ref{cond1}) is equivalent to
$\delbar(C-\del\alpha)=0$, whose general solution is
\[C=\del\alpha+\overline{\del\alpha} + i\del\delbar\chi\]for any real
function $\chi$. We must now check that it is possible to choose
$\chi$ so that the second condition (\ref{cond2}) is satisfied by
this $C$.

\item From equation (\ref{second}) we obtain that
$\delbar(A^{101}-d_f\alpha)=0$, implying that
$A^{101}=d_f\alpha+\delbar\beta$ for some $(100)$-form $\beta$.
Condition (\ref{cond2}) then is equivalent to the fact that
\begin{equation*}
-id_f\del\delbar\chi=\del\delbar(\beta-\overline\beta),
\end{equation*}
which can be solved (for the unknown $\chi$) if and only if the
right hand side is $d_f$-closed.  From equation (\ref{third}) we
see that $\delbar(A^{200}-d_f\beta)=0$, showing that
$A^{200}=d_f\beta + \delta$, where $\delta$ is a $\delbar$-closed
$(200)$-form. Hence
\begin{equation*}
d_f\del\delbar(\beta-\overline\beta)=\del\delbar(A^{200}-\overline{A^{200}}),
\end{equation*}
and the right hand side vanishes precisely because
$A^{200}-\overline{A^{200}}=2\omega_0$, which is closed.  Hence
$\chi$ may be chosen to satisfy condition (\ref{cond2}), and so we
obtain a closed 2-form $\tilde B$.
\end{itemize}
Finally, we have proven the normal coordinate theorem for regular
neighbourhoods of generalized complex manifolds:
\begin{theorem}[Generalized Darboux theorem]
Any regular point in a generalized complex manifold has a
neighbourhood which is equivalent, via a diffeomorphism and a
B-field transformation, to the product of an open set in $\CC^k$
with an open set in the standard symplectic space
$(\RR^{2n-2k},\omega_0)$.
\end{theorem}

\section{The jumping phenomenon}\label{sec8}

While we have fully characterized generalized complex structures
in regular neighbourhoods, it remains an essential feature of
generalized complex geometry that the \emph{type} of the structure
may vary throughout the manifold. In view of the normal form
theorem, the type can be thought of as the number of transverse
complex directions, and this is an upper semi-continuous function
on the manifold, i.e. each point has a neighbourhood in which it
does not increase. The most generic type is zero, when there are
only symplectic directions, and the most special type is $n$, when
all directions are complex. Note that the type may jump up, but
always by an even number, since types of the same parity are in
the same connected component of the space of linear generalized
complex structures (see section~\ref{maxiso} for details).

In this section we present a simple example of a generalized
complex structure on $\RR^4$ which is of symplectic type $(k=0)$
outside a codimension 2 hypersurface and jumps up to complex type
$(k=2)$ along the hypersurface.

Consider the differential form
\begin{equation*}
\rho=z_1 + dz_1\wedge dz_2,
\end{equation*}
where $z_1,z_2$ are the standard coordinates on $\CC^2\cong\RR^4$.
When $z_1=0$, $\rho=dz_1\wedge dz_2$ and so defines the standard
complex structure, whereas when $z_1\neq 0$, $\rho$ defines a
$B$-symplectic structure since
\begin{equation*}
\rho=z_1 e^{\frac{dz_1\wedge dz_2}{z_1}}.
\end{equation*}
Hence, algebraically the form $\rho$ defines a generalized almost
complex structure whose type jumps along $z_1=0$.

To check the integrability of this structure, we take the exterior
derivative:
\begin{equation*}
d\rho=dz_1=i_{-\partial_{z_2}}(z_1+dz_1\wedge dz_2)=
(-\partial_{z_2})\cdot\rho,
\end{equation*}
showing that $\rho$ indeed satisfies the integrability condition
of Theorem~\ref{deldelbar}, and defines a generalized complex
structure on all of $\RR^4$.  In this case it is easy to see that
although the canonical line bundle is topologically trivial, it is
not holomorphically trivial, i.e. there is no nowhere-vanishing
closed section of $U$.

In the next chapter we will produce more general examples of the
jumping phenomenon, and on compact manifolds as well.

\section{Twisted generalized complex
structures}\label{twistcx}

Following on from section~\ref{twistedcour}, where we described
the twisted Courant bracket, we now define the notion of
\emph{twisted} generalized complex structure.  The underlying
algebraic structure is the same as the generalized complex case:
\begin{defn}
A generalized almost complex structure $\JJ$ is said to be twisted
generalized complex with respect to the closed 3-form $H$ when its
$+i$-eigenbundle $L$ is involutive with respect to the $H$-twisted
Courant bracket.
\end{defn}
\begin{remark}
It is important to note that given any $H$-twisted generalized
complex structure $\JJ$, the conjugate $e^b\JJ e^{-b}$, for $b$
any smooth 2-form, is integrable with respect to the
$H+db$-twisted Courant bracket.  This means that the space of
twisted generalized complex structures depends only on the
cohomology class $[H]\in H^3(M,\RR)$.
\end{remark}
As we described in Proposition~\ref{twistprop}, the integrability
condition on the differential forms defining the generalized
almost complex structure is simply that $d_H=d+H\wedge\cdot$ maps
$C^\infty(U_0)$ to $C^\infty(U_1)$. A Darboux theorem for
$H$-twisted structures follows from the methods developed in
section~\ref{sec7}: in a regular neighbourhood an $H$-twisted
generalized complex structure can be expressed as the B-transform
of a product of a symplectic by a complex structure, except that
in the twisted case $B$ is not a closed form. Instead, $dB=H$ in
the neighbourhood.

We will provide interesting examples of twisted generalized
complex structures in section~\ref{twistkahler}.

\clearpage
\thispagestyle{empty}
\cleardoublepage
\chapter{Deformations of generalized complex
structures}\label{defm}

In the deformation theory of complex manifolds developed by
Kodaira, Spencer, and Kuranishi, one begins with a compact complex
manifold $(M,J)$ with holomorphic tangent bundle $\TT$, and
constructs an analytic subvariety $\zed\subset H^1(M,\TT)$
(containing $0$) which is the base space of a family of
deformations $\MM=\{\eps_z\ :\ z\in \zed,\ \eps_0=0\}$ of the
original complex structure $J$.  This family is \emph{locally
complete} (also called \emph{miniversal}), in the sense that any
family of deformations of $J$ can be obtained, up to equivalence,
by pulling $\MM$ back by a map $f$ to $\zed$, as long as the
family is restricted to a sufficiently small open set in its base.
The subvariety $\zed\subset H^1(M,\TT)$ is defined as the zero set
of a holomorphic map $\Phi:H^1(M,\TT)\rightarrow H^2(M,\TT)$, and
so the base of the miniversal family is certainly smooth when this
obstruction map vanishes.

In this section we extend these results to the generalized complex
setting, following the method of Kuranishi~\cite{Kuranishi}.  In
particular, we construct, for any generalized complex manifold, a
locally complete family of deformations.  We then proceed to
produce new examples of generalized complex structures by
deforming known ones.

\section{The deformation complex}\label{deform}

The generalized complex structure $\JJ$ is determined by its
$+i$-eigenbundle $L<(T\oplus T^*)\otimes\CC$ which is isotropic,
satisfies $L\cap\overline{L}=\{0\}$, and is closed under the
Courant bracket. Recall that since $(T\oplus
T^*)\otimes\CC=L\oplus\overline L$, we use the natural metric
$\IP{,}$ to identify $\overline L$ with $L^*$.

To deform $\JJ$ we will vary $L$ in the Grassmannian of maximal
isotropics. Any maximal isotropic having zero intersection with
$\overline{L}$ (this is an open set containing $L$) can be
uniquely described as the graph of a homomorphism $\epsilon:L\lra
\overline L$ satisfying $\IP{\epsilon X,Y }+\IP{X,\epsilon Y}=0\ \
\forall X,Y\in C^\infty(L)$, or equivalently $\epsilon\in
C^\infty(\wedge^2 L^*)$. Therefore the new isotropic is given by
$L_\epsilon=(1+\epsilon)L$. As the deformed $\JJ$ is to remain
real, we must have $\overline
L_\epsilon=(1+\overline\epsilon)\overline L$.  Now we observe that
$L_\epsilon$ has zero intersection with its conjugate if and only
if the endomorphism we have described on $L\oplus L^*$, namely
\begin{equation}\label{matrixa}
A_\epsilon=\left(\begin{matrix}1&\overline\epsilon\\\epsilon&1\end{matrix}\right),
\end{equation}
is invertible; this is the case for $\epsilon$ in an open set
around zero.

So, providing $\epsilon$ is small enough,
$\JJ_\epsilon=A_\epsilon\JJ A_\epsilon^{-1}$ is a new generalized
almost complex structure, and all nearby almost structures are
obtained in this way.  Note that while $A_\epsilon$ itself is not
an orthogonal transformation, of course $\JJ_\epsilon$ is.

To describe the integrability condition for $\epsilon\in
C^\infty(\wedge^2 L^*)$ which guarantees that $\JJ_\epsilon$ is
integrable, we notice that we are precisely in the situation of
Theorem~\ref{masterequation}, which states that $\JJ_\epsilon$ is
integrable if and only if $\epsilon\in C^\infty(\wedge^2 L^*)$
satisfies the equation
\begin{equation}\label{master}
d_L\epsilon +\frac{1}{2} [\epsilon,\epsilon]=0.
\end{equation}
In this way, we can interpret sufficiently small solutions of
equation~(\ref{master}) as deformations of a genuine geometrical
structure, thereby solving the open problem stated in~\cite{LWX}
of interpreting the equation as deriving from a deformation
theory.

Infinitesimally, this means that nearby generalized complex
structures are in the kernel of the $d_L:C^\infty(\wedge^2
L^*)\rightarrow C^\infty(\wedge^3 L^*)$, and so we need to know
the nature of this linear operator.
\begin{prop}
If $L$ is the Lie algebroid deriving from a generalized complex
structure, then the Lie algebroid differential complex
\begin{equation*}
d_L:C^\infty(\wedge^p L^*)\rightarrow C^\infty(\wedge^{p+1}L^*)
\end{equation*}
is elliptic.  Hence its cohomology groups, which we denote by
$H^p_L(M)$, are finite dimensional complex vector spaces.
\end{prop}
\begin{proof}
The principal symbol $s(d_L):(T^*\otimes\CC)\otimes\wedge^p
L^*\rightarrow \wedge^{p+1}L^*$ is given by
$\pi^*:T^*\otimes\CC\rightarrow L^*$ (where $\pi:L\rightarrow
T\otimes\CC$ is the projection) composed with wedge product, i.e.
\begin{equation*}
s_\xi(d_L)=\pi^*(\xi)\wedge\cdot\ \ ,
\end{equation*}
where $\xi\in T^*\otimes\CC$. If $\xi$ is a nonzero real 1-form,
then since $(T\oplus T^*)\otimes\CC=L\oplus\overline{L}$, we have
the decomposition $\xi=x+\overline x$ for a nonzero element $x\in
L$, and then $s_\xi(d_L)=\overline x\wedge\cdot$ defines an
elliptic Koszul complex. Hence $(C^\infty(\wedge^\bullet
L^*),d_L)$ is an elliptic differential complex.  See~\cite{Wells}
for a proof that the cohomology groups of an elliptic differential
complex on a compact manifold are finite dimensional.
\end{proof}

Now that we have described which tensors $\eps\in
C^\infty(\wedge^2 L^*)$ are integrable deformations, we must
explain when two such deformations are considered equivalent.  In
the case of complex or symplectic geometry, two deformations are
considered equivalent if they are related by a small
diffeomorphism.   In the case of generalized complex geometry,
however, the Courant bracket on $T\oplus T^*$ has a larger group
of symmetries, and so we will consider two deformations to be
equivalent if they are related by a diffeomorphism (connected to
the identity) and an \emph{exact} B-field transformation.  A
special case of such a transformation is one generated by a vector
field $X$ and a 1-form $\xi$:
\begin{equation*}
F_{X+\xi}=e^{d\xi}\circ {e^X}
\end{equation*}
for $X+\xi\in C^\infty(T\oplus T^*)$, where by $e^{tX}$ we mean
the one-parameter group of diffeomorphisms generated by the vector
field $X$.  Note that there is redundancy in expressing a symmetry
as a section of $T\oplus T^*$, since the 1-form $\xi$ could be
exact, in which case it has no action whatsoever.
\begin{remark}
In the presence of a generalized complex structure $\JJ$, certain
infinitesimal symmetries $X+\xi\in C^\infty(T\oplus T^*)$ preserve
the tensor $\JJ$; these \emph{generalized holomorphic symmetries}
are given by sections of $T\oplus T^*$ whose $L^*$ component lies
in the kernel of $d_L$.  In the complex case these are holomorphic
vector fields together with $\delbar$-closed $(0,1)$-forms, while
in the symplectic case they arise from symplectic vector fields,
i.e. sections $X\in C^\infty(T)$ such that $\LL_X\omega=0$.  Note
that any complex-valued function $f\in C^\infty(M,\CC)$ generates
a symmetry $X+\xi=d_Lf+\overline{d_Lf}$; such holomorphic
symmetries could be called \emph{Hamiltonian} symmetries of the
generalized complex structure, and do coincide with the notion of
Hamiltonian vector field in the symplectic case.  In the complex
case they are given by 1-forms $\xi=\delbar f + \del
\overline{f}$, and generate $B$-field transformations with
$B=\del\delbar(f-\overline{f})$.
\end{remark}
By differentiating, we now show the relationship between sections
in the image of $d_L:C^\infty(L^*)\rightarrow C^\infty(\wedge^2
L^*)$ and equivalent deformations.
\begin{prop}
Let $\JJ$ be a generalized complex structure with $+i$-eigenbundle
$L$, and let $\eps_0\in C^\infty(\wedge^2 L^*)$ be a deformation
of this structure.  If $X+\xi\in C^\infty(T\oplus T^*)$, and if
$t\in\RR$ is in a sufficiently small neighbourhood of $0$,
then we have the following expression for the equivalent
deformation induced by $F_{t(X+\xi)}=e^{td\xi}\circ e^{tX}$ by its
action on the graph of $\eps_0$:
\begin{equation}\label{transdefm}
F_{t(X+\xi)}(\eps_0)=\eps_0 +
td_L((X+\xi)_{L^*})+R(\eps_0,t(X+\xi)),
\end{equation}
where $(X+\xi)_{L^*}$ is the component in $L^*$ according to the
splitting $L\oplus L^*$, and $R$ satisfies
\begin{equation*}
R(t\eps_0,t(X+\xi))=t^2\widetilde R(\eps_0,X+\xi,t),
\end{equation*}
where $\widetilde R(\eps_0,X+\xi,t)$ is smooth. In this sense,
$R(\eps_0,X+\xi)$ is $O(t^2)$ in the deformations.
\end{prop}
\begin{proof}
Let $\{s\eps_0\ :\ s\in\RR\}$ be a
straight line in the space of smooth sections of $\wedge^2 L^*$
passing through $0$ and the deformation $\eps_0$.  Let
$F_{t(X+\xi)}$, for $t$ in some neighbourhood of $0\in\RR$, be the
family of automorphisms of $T\oplus T^*$ defined by
$F_{t(X+\xi)}=e^{td\xi}\circ e^{tX}$. The combined action of the
section $s\eps_0$ and the automorphism $F_{t(X+\xi)}$ on the
$+i$-eigenbundle $L$ is given by the composition
\begin{equation*}
F_{t(X+\xi)}A_{s\eps_0}=\begin{pmatrix}\sigma&\overline{\tau}\\\tau&\overline{\sigma}\end{pmatrix},
\end{equation*}
where $A_{s\eps_0}$ is as in (\ref{matrixa}) and the right hand
side is written in the splitting $L\oplus\overline{L}$.  Assuming
$t$ is small enough, $\sigma$ is invertible, and we may factorise
\begin{equation*}
F_{t(X+\xi)}A_{s\eps_0}=\begin{pmatrix}\sigma&\\&\overline{\sigma}\end{pmatrix}\begin{pmatrix}1&\overline{\eps(s,t)}\\\eps(s,t)&1\end{pmatrix}
=C_\sigma A_{\eps(s,t)},
\end{equation*}
where $\eps(s,t)\in C^\infty(\wedge^2 L^*)$ is the new section of
$\wedge^2 L^*$ given by the action of $F_{t(X+\xi)}$ on $s\eps_0$,
i.e.
\begin{equation*}
F_{t(X+\xi)}(s\eps_0)=\eps(s,t).
\end{equation*}

Differentiating $A_{\eps(s,t)}=C^{-1}_\sigma
F_{t(X+\xi)}A_{s\eps_0}$ and evaluating at $(s,t)=(0,0)$, we
obtain
\begin{equation*}
\dot{A}_{\eps(s,t)}\big|_{(0,0)}=-\dot{C}_\sigma\big|_{(0,0)} +
\dot{F}_{t(X+\xi)}\big|_{(0,0)} + \dot{A}_{s\eps_0}\big|_{(0,0)}.
\end{equation*}
Now let $A+a,B+b\in C^\infty(L)$, so that
\begin{align*}
\eps(s,t)(A+a,B+b)=\IP{A_{\eps(s,t)}(A+a),B+b}.
\end{align*}
Differentiating, we see that since $C_\sigma$ is simply an
automorphism of $L$, the $\dot{C}_\sigma$ term has no effect,
leaving two terms:
\begin{align*}
\dot{\eps}(s,t)(A+a,B+b)=\IP{\dot{F}_{t(X+\xi)}\big|_{(0,0)}A+a,B+b}+\IP{\dot{A}_{s\eps_0}\big|_{(0,0)}A+a,B+b}.
\end{align*}
We compute each term separately:
\begin{align*}
\tfrac{\partial}{\partial
s}\eps(s,t)\big|_{t=s=0}(A+a,B+b)&=\tfrac{d}{ds}\IP{A_{s\eps_0}(A+a),B+b}\big|_{s=0}\\
&=\eps_0(A+a,B+b),
\end{align*}
while
\begin{align*}
\tfrac{\partial}{\partial
t}\eps(s,t)\big|_{t=s=0}(A+a,B+b)&=\tfrac{d}{dt}\IP{F(t(X+\xi))(A+a),B+b}\big|_{t=0}\\
&=\tfrac{d}{dt}\IP{e^{tX}_*A+{(e^{tX})^*}^{-1}a+i_{e^{tX}_*A}(td\xi),B+b}\big|_{t=0}\\
&=-\IP{[X,A]+\LL_Xa+i_Ad\xi,B+b}.
\end{align*}
Now note that
\begin{align*}
d_L(X+\xi)_{L^*}(A+a,B+b)&=i_Ad\IP{X+\xi,B+b}-i_Bd\IP{X+\xi,A+a}-\IP{X+\xi,[A+a,B+b]}\\
&=-\tfrac{1}{2}(i_B\LL_Xa+i_Bi_Ad\xi+[\LL_X,i_A]b)\\
&=-\IP{[X,A]+\LL_Xa+i_Ad\xi,B+b}\\
&=\tfrac{\partial}{\partial t}\eps(s,t)\big|_{t=s=0}(A+a,B+b).
\end{align*}
Therefore, by Taylor's theorem,
\begin{equation}\label{epseps}
\eps(s,t)=s\eps_0 + td_L(X+\xi)_{L^*} + R(s\eps_0,t(X+\xi)),
\end{equation}
where $R(s\eps_0,t(X+\xi))$ is of order $s^2, st$, and $t^2$, so
that, taking $s=t$,
\begin{equation*}
R(t\eps_0,t(X+\xi))=t^2\widetilde R(\eps_0,X+\xi,t)
\end{equation*}
with $\widetilde R(\eps_0,X+\xi,t)$ smooth.  Setting $s=1$ in
(\ref{epseps}), we obtain the required result.
\end{proof}
This proposition shows us that, infinitesimally, deformations
which differ by sections which lie in the image of $d_L$ in
$C^\infty(\wedge^2 L^*)$ are equivalent by transformations
generated by vector fields and 1-forms. Hence we expect the
tangent space to the moduli space to lie in $H^2_L(M)$. In the
remainder of this section we develop the Hodge theory for the
elliptic complex $(C^\infty(\wedge^\bullet L^*),d_L)$ so that we
may prove this assertion rigorously.

We follow the usual treatment of Hodge theory as described
in~\cite{Wells}. Choose a Hermitian metric on the complex Lie
algebroid $L$ and let $(\varphi,\psi)_k$ be the $L^ 2_k$ Sobolev
inner product on sections $\varphi,\psi\in C^\infty(\wedge^p L^*)$
induced by the metric.  Recall that $L^2_0$ is the usual $L^2$
inner product. We will use $|\ |_k$ to denote the Sobolev norm,
defined by
\begin{equation*}
|u|_k=((u,u)_k)^{1/2}.
\end{equation*}
We define the $L^2$ adjoint $d_L^*$ of $d_L$ via the formula
\begin{equation*}
(d_L\varphi,\phi)=(\varphi,d_L^*\phi),
\end{equation*}
and obtain the elliptic, self-adjoint Laplacian
\begin{equation*}
\Delta_L=d_Ld^*_L+d^*_Ld_L.
\end{equation*}
Let $\Har^p$ be the space of $\Delta_L$-harmonic forms, which is
isomorphic to $H^p_L(M)$ by the standard argument, and let $H$ be
the orthogonal projection of $C^\infty(\wedge^p L)$ onto the
closed subspace $\Har^p$.  Also, let $G$ be the Green smoothing
operator quasi-inverse to $\Delta_L$, i.e.
\begin{equation*}
G:L^2_k\rightarrow L^2_{k+2},
\end{equation*}
so that we have
\begin{equation*}
\mathrm{Id}=H+\Delta G=H+G\Delta.
\end{equation*}
Also, $G$ satisfies $[G,d_L]=[G,d^*_L]=0$.  We will find it
useful, as Kuranishi did, to define the once-smoothing operator
\begin{equation*}
Q=d_L^*G:L^2_k\rightarrow L^2_{k+1},
\end{equation*}
which then satisfies
\begin{align}\label{Qrel}
\mathrm{Id}&=H+d_LQ+Q d_L,\\
Q^2=d^*_LQ&=Qd^*_L=HQ=QH=0.
\end{align}
We now have the algebraic and analytical tools we need to describe
the deformation theory of generalized complex manifolds.

\section{The deformation theorem}

Let $\JJ$ be a generalized complex structure on the compact
manifold $M$ with $+i$-eigenbundle $L$, and let $\eps\in
C^\infty(\wedge^2 L^*)$. Then recall that since $(L,\overline L)$
is a Lie bialgebroid, we have compatibility between the Courant
bracket on $L^*$ and the differential $d_L$.

\begin{theorem}[Deformation theorem for generalized complex
structures] There exists an open neighbourhood $U\subset H^2_L(M)$
containing zero, a smooth family $\widetilde\MM=\{\eps_u\ :\ u\in
U,\ \eps_0=0\}$ of generalized almost complex deformations of
$\JJ$, and an analytic obstruction map $\Phi:U\rightarrow
H^3_L(M)$ with $\Phi(0)=0$ and $d\Phi(0)=0$, such that the
deformations in the sub-family $\MM=\{\eps_z\ :\
z\in\zed=\Phi^{-1}(0)\}$  are precisely the integrable ones.
Furthermore, any sufficiently small deformation $\eps$ of $J$ is
equivalent to at least one member of the family $\MM$. Finally, in
the case that the obstruction map vanishes, we show that $\MM$ is
a smooth locally complete family.
\end{theorem}
\begin{proof} The proof is divided into two parts: first, we
construct a smooth family $\widetilde\MM$, and show it contains
the family of integrable deformations $\MM$ defined by the map
$\Phi$; second, we describe its miniversality property.  We follow
the paper of Kuranishi~\cite{Kuranishi} closely.

{\bf Part I:} For sufficiently large $k$, $L^2_k(M,\RR)$ is a
Banach algebra (see~\cite{Palais}), and the map $F:\eps\mapsto
\eps+\tfrac{1}{2}Q[\eps,\eps]$ extends to a map
\begin{equation*}
F:L^2_k(\wedge^2 L^*)\lra L^2_k(\wedge^2 L^*),
\end{equation*}
which is a smooth map of the Hilbert space into itself whose
derivative at the origin is clearly the identity mapping. Hence by
the Banach space inverse function theorem, $F^{-1}$ maps a
neighbourhood of the origin in $L^2_k(\wedge^2 L^*)$ smoothly and
bijectively to another neighbourhood of the origin. Hence, if we
choose a sufficiently small $\delta>0$ then the following
finite-dimensional subset of harmonic sections
\begin{equation*}
U=\{u\in\Har^2<L^2_k(\wedge^2L^*)\ : \ |u|_k<\delta\}
\end{equation*}
is taken by $F^{-1}$ into another smooth finite-dimensional set,
i.e. defines a family of sections
\begin{equation*}
\widetilde M=\{\eps(u)=F^{-1}(u)\ :\ u\in U\}
\end{equation*}
where $\eps(u)$ depends smoothly (in fact, holomorphically) on
$u$, and satisfies $F(\eps(u))=u$. Applying the Laplacian to this
equation, we obtain
\begin{equation*}
\Delta_L\eps(u) + \tfrac{1}{2}\Delta_LQ[\eps(u),\eps(u)]=0,
\end{equation*}
and using $\Delta_L Q =\Delta_L d^*_LG=d^*_L\Delta_LG=d^*_L$, we
see that
\begin{equation*}
\Delta_L\eps(u) + \tfrac{1}{2}d^*_L[\eps(u),\eps(u)]=0.
\end{equation*}
This is a quasi-linear elliptic PDE, and by the standard result of
Morrey~\cite{Morrey}, we conclude that the solutions $\eps(u)$ of
this equation are actually smooth sections, i.e.
\begin{equation*}
\eps(u)\in C^\infty(\wedge^2 L^*).
\end{equation*}
Hence we have constructed a smooth family of generalized almost
complex deformations of $\JJ$, over an open set $U\subset
H^2_L(M)$.  We emphasize that we have exhibited this family as a
genuine finite-dimensional submanifold of a Hilbert space whose
tangent space at the origin is $\Har^2=H^2_L(M)$.

We now ask which of these deformations are integrable, i.e.
satisfy the equation $d_L\eps+\tfrac{1}{2}[\eps,\eps]=0$. Since
$\eps(u)+\tfrac{1}{2}Q[\eps(u),\eps(u)]=u$ and $d_Lu=d^*_Lu=0$, we
see that $d_L\eps(u)=-\tfrac{1}{2}d_LQ[\eps(u),\eps(u)]$, and
using~(\ref{Qrel}) we obtain
\begin{align*}
d_L\eps(u)+\tfrac{1}{2}[\eps(u),\eps(u)]&=-\tfrac{1}{2}d_LQ[\eps(u),\eps(u)]+\tfrac{1}{2}[\eps(u),\eps(u)]\\
&=\tfrac{1}{2}(Qd_L+H)[\eps(u),\eps(u)].
\end{align*}
Since the images of $Q$ and $H$ are orthogonal in the $L_2$ inner
product, we see that $\eps(u)$ is integrable if and only if
$H[\eps(u),\eps(u)]=Qd_L[\eps(u),\eps(u)]=0$.

Now we claim that $H[\eps(u),\eps(u)]=0$ implies that
$Qd_L[\eps(u),\eps(u)]=0$: using the compatibility of $[,]$ and
$d_L$ we obtain
\begin{align*}
Qd_L[\eps(u),\eps(u)]&=2Q[d_L\eps(u),\eps(u)]\\
&=-Q[d_LQ[\eps(u),\eps(u)],\eps(u)]\\
&=-Q[(\mathrm{Id}-Qd_L-H)[\eps(u),\eps(u)],\eps(u)].
\end{align*}
So, assuming that $H[\eps(u),\eps(u)]=0$, we obtain
\begin{align*}
Qd_L[\eps(u),\eps(u)]&=-Q[(\mathrm{Id}-Qd_L)[\eps(u),\eps(u)],\eps(u)]\\
&=Q[Qd_L[\eps(u),\eps(u)],\eps(u)].
\end{align*}
Letting $\zeta(u)=Qd_L[\eps(u),\eps(u)]$, we have that
\begin{equation*}
\zeta(u)=Q[\zeta(u),\eps(u)],
\end{equation*}
and since for sufficiently large $k$ the map
$(\alpha,\beta)\mapsto Q[\alpha,\beta]$ satisfies
$|Q[\alpha,\beta]|_k\leq c|\alpha|_k|\beta|_k$ for some $c>0$, we
have
\begin{equation*}
|\zeta(u)|_k\leq c|\zeta(u)|_k|\eps(u)|_k
\end{equation*}
for some $c>0$.  Therefore, if we take $\delta$ to be so small
that $|\eps(u)|_k<\tfrac{1}{c}$ for all $|u|_k<\delta$, we obtain
that $\zeta(u)=0$.

Hence, we have shown that $\eps(u)$ is integrable precisely when
$u$ lies in the vanishing set of the analytic mapping
$\Phi:U\rightarrow H^3_L(M)$ defined by
$\Phi(u)=H[\eps(u),\eps(u)]$.  Note that $\Phi(0)=d\Phi(0)=0$.
Before we proceed to the second part of the proof, we wish to give
an alternative characterisation of the family $\MM=\{\eps(z)\ :\
z\in \zed=\Phi^{-1}(0) \}$.  We claim that $\MM$ is actually a
neighbourhood around zero in the set
\begin{equation*}
\MM'=\left\{ \eps\in C^\infty(\wedge^2 L^*)\ :\
d_L\eps+\tfrac{1}{2}[\eps,\eps]=d^*_L\eps=0\right\}.
\end{equation*}
To show this, let $\eps(u)\in \MM$. Then since
$\eps(u)=u-\tfrac{1}{2}Q[\eps(u),\eps(u)]$ and $d^*_LQ=0$, we see
that $d^*_L\eps(u)=0$, showing that $\MM\subset\MM'$.  Conversely,
let $\eps\in\MM'$.  Then since $d^*_L\eps=0$, applying $d^*_L$ to
the equation $d_L\eps+\tfrac{1}{2}[\eps,\eps]=0$ we obtain
$\Delta_L\eps+\tfrac{1}{2}d^*_L[\eps,\eps]=0$, and applying
Green's operator we see that
$\eps+\tfrac{1}{2}Q[\eps,\eps]=H\eps$, i.e.
$F(\eps)=H\eps\in\Har^2$, proving that a small open set in $\MM'$
is contained in $\MM$, completing the argument.

{\bf Part II:} Let $P<C^\infty(L^*)$ be the $L^2$ orthogonal
complement of the $d_L$-closed sections $\ker d_L<C^\infty(L^*)$,
or in other words, sections in the image of $d^*_L$. We show that
there exist neighbourhoods of the origin $V\subset
C^\infty(\wedge^2 L^*)$ and $W\subset P$ such that for any
$\eps\in V$ there is a unique $X+\xi\in C^\infty(T\oplus T^*)$
such that $(X+\xi)_{L^*}\in W$ and
\begin{equation}
d^*_L(e^{d\xi}e^X(\eps))=0.
\end{equation}
This would imply that any sufficiently small solution to
$d_L\eps+\tfrac{1}{2}[\eps,\eps]=0$ is equivalent to another
solution $\eps'$ such that $d^*_L\eps'=0$, i.e. a solution in
$\MM$. Extended to smooth families, this result would prove local
completeness.

Restricting to a sufficiently small neighbourhood in
$C^\infty(T\oplus T^*)$ so that we may take $t=1$ in
Equation~(\ref{transdefm}), we see that
$d^*_L(e^{d\xi}e^X(\eps))=0$ if and only if
\begin{equation*}
d^*_L\eps +d^*_Ld_L(X+\xi)_{L^*}+d^*_LR(\eps,X+\xi)=0.
\end{equation*}
Assuming $(X+\xi)_{L^*}\in P$, we see that
$d^*_L(X+\xi)_{L^*}=H(X+\xi)_{L^*}=0$, so that
\begin{equation*}
d^*_L\eps +\Delta_L(X+\xi)_{L^*}+d^*_LR(\eps,X+\xi)=0,
\end{equation*}
and applying $G$,
\begin{equation*}
(X+\xi)_{L^*}+Qd^*_L\eps +QR(\eps,X+\xi)=0.
\end{equation*}
Since $R(\eps,X+\xi)$ involves one derivative of $X+\xi$, the map
\begin{equation*}
F:(\eps,X+\xi)\mapsto (X+\xi)_{L^*}+Qd^*_L\eps +QR(\eps,X+\xi)
\end{equation*}
is continuous from a neighbourhood of the origin $V_0\times W_0$
in $C^\infty(\wedge^2 L^*)\times P$ (where $R(\eps,X+\xi)$ is
defined) to $P$, where all spaces are endowed with the $L^2_k$
norm, $k$ sufficiently large. $F$ can therefore be extended to a
continuous map from the completion of the domain,
$\widehat{V_0}\times \widehat{W_0}$, to the completion of $P$. The
derivative of $F$ with respect to $X+\xi$ is the identity map, and
so by the implicit function theorem there are neighbourhoods
$V\subset V_0, W_1\subset \widehat{W_0}$ such that given $\eps\in
V$, $F(\eps,X+\xi)=0$ is satisfied for a unique $(X+\xi)_{L^*}\in
W_1$, and which depends smoothly on $\eps\in V$. Furthermore,
since $\eps\in V$ is itself smooth, the unique solution $X+\xi$
satisfies the quasi-linear elliptic PDE
$\Delta_L(X+\xi)_{L^*}+d^*_L\eps+d^*_LR(\eps,X+\xi)=0$, implying
that $X+\xi$ is smooth as well, hence $(X+\xi)_{L^*}$ lies in the
neighbourhood $W=W_1\cap P$. Therefore we have shown that every
sufficiently small deformation of the generalized complex
structure is equivalent to one in our finite-dimensional family
$\MM$.

If the obstruction map $\Phi$ vanishes, so that $\MM$ is a smooth
family, then given any other smooth family $\MM_S=\{\eps_s\ :\
s\in S,\ \eps_{s_0}=0\}$ with basepoint $s_0\in S$, the above
argument provides a smooth family of equivalences $(X+\xi)_s$
taking each $\eps_s$ for $s$ in some neighbourhood $T$ of $s_0$ to
$\eps_{f(s)}$, $f(s)\in U\subset H^2_L(M)$, defining a smooth map
$f:T\rightarrow U$, $f(s_0)=0$, so that $f^*\MM=\MM_S$.  Thus we
establish that $\MM$ is a locally complete family of deformations.
\end{proof}
\begin{remark}
The natural complex structure on $H^2_L(M)$ and on the vanishing
set of the holomorphic obstruction map $\Phi$ raises the question
of whether there is a notion of holomorphic family of generalized
complex structures.  There is: if $S$ is a complex manifold then a
holomorphic family of generalized complex structures on $M$ is a
generalized complex structure on $M\times S$ which can be pushed
down via the projection to yield the complex structure on $S$. The
family $\MM$ is actually such a holomorphic family.
\end{remark}

\section{Examples of deformed structures}

Consider deforming a complex manifold $(M,J)$ as a generalized
complex manifold.  Then since
\begin{equation*}
L=T_{0,1}\oplus T^*_{1,0},
\end{equation*}
the deformation complex is actually
\begin{equation*}
(\oplus_{p+q=k}\Omega^{0,q}(\wedge^p T_{1,0}),\ \ \delbar),
\end{equation*}
so that the base of the Kuranishi family lies in
\begin{equation*}
H^2_L(M)=\oplus_{p+q=2}H^q(M,\wedge^p T_{1,0}).
\end{equation*}
The image of the obstruction map lies in
\begin{equation*}
H^3_L(M)=\oplus_{p+q=3}H^q(M,\wedge^p T_{1,0})
\end{equation*}
Therefore we see immediately that generalized complex manifolds
provide a solution to the problem of finding a geometrical
interpretation of the ``extended complex deformation space''
defined by Kontsevich and Barannikov~\cite{KontsevichBarannikov}.
Any deformation $\eps$ has three components
\begin{equation*}
\beta\in H^0(M,\wedge^2T_{1,0}),\ \ \ \varphi\in H^1(M,T_{1,0}),\
\ \ B\in H^2(M,\mathcal{O}).
\end{equation*}
The component $\varphi$ is a usual deformation of the complex
structure, as discovered by Kodaira and Spencer. The component $B$
is a complex $B$-field action as we have discussed.  The component
$\beta$, however, is a completely new type of deformation for
complex manifolds.  The integrability condition on such a
deformation $\beta\in C^\infty(\wedge^2 T_{1,0})$ is simply that
\begin{equation*}
\delbar\beta + \tfrac{1}{2}[\beta,\beta]=0,
\end{equation*}
which is satisfied if and only if the bivector $\beta$ is
holomorphic and Poisson.  As we saw in section~\ref{symm}, a
$\beta$-transform acts by shearing $T\oplus T^*$ in the $T$
direction, and hence may change the type of the generalized
complex structure.  As an example, let us explore such
deformations on $\CC P^2$.

\begin{example}[Deformed generalized complex structure on $\CC
P^2$]\label{defmcp2} On $\CC P^2$, $\wedge^2 T_{1,0}$ is simply
the anticanonical bundle $\mathcal{O}(3)$, whose nonzero
holomorphic sections vanish along cubics. Any holomorphic bivector
$\beta\in H^0(M,\mathcal{O}(3))$ is automatically Poisson since we
are in complex dimension 2, and hence any sufficiently small
holomorphic section of $\mathcal{O}(3)$ defines an integrable
deformation of the complex structure into a generalized complex
structure.

$\beta$ takes $X+\xi\in T_{0,1}\oplus T^*_{1,0}$ to
$X+\xi+i_{\xi}\beta$ and so whenever $\beta$ is nonzero, the
deformed Lie algebroid projects surjectively onto $T\otimes\CC$.
Hence the deformed structure is of B-symplectic type (type 0)
outside the cubic vanishing set and of complex type (type 2) along
the cubic.  The B-field and symplectic form go to infinity as one
approaches the cubic curve. Therefore we obtain a generalized
complex structure on $\CC P^2$ which is clearly inequivalent to
either a complex or a symplectic structure.  This compact
generalized complex manifold exhibits a jumping phenomenon along a
codimension 2 subvariety.

To check the smoothness of the locally complete family, note that
since we are in complex dimension 2, $H^0(\wedge^3
T_{1,0})=H^3(\mathcal{O})=0$. Also, by Serre duality,
$H^2(T_{1,0})\cong H^0(T^*_{1,0}\otimes \mathcal{O}(-3))=0$, since
$T^*_{1,0}$ has no holomorphic sections.  Furthermore,
$H^1(\wedge^2 T_{1,0})=H^1(\mathcal{O}(3))=0$ by the Bott
formulae.  Hence we see that the obstruction space vanishes for
$\CC P^2$, and so we conclude that there is a smooth locally
complete family of deformations for $\CC P^2$ as a generalized
complex manifold.
\end{example}

\begin{example}
One can of course deform $\CC^2$ in the same way that we have
deformed $\CC P^2$; we choose the holomorphic bivector
\begin{equation*}
\beta=z_1\del_{z_1}\wedge\del_{z_2},
\end{equation*}
where $z_1,z_2$ are the usual complex coordinates. Then applying a
$\beta$-transform to the usual complex structure defined by the
spinor $\Omega = dz_1\wedge dz_2$, we obtain
\begin{equation*}
e^{\beta}\Omega= dz_1\wedge dz_2 + z_1,
\end{equation*}
which is precisely the example we provided in section~\ref{sec8}
of a jumping generalized complex structure on $\CC^2$.  We see now
that it is actually a deformation of the usual complex structure
by a holomorphic Poisson structure.  To see the behaviour of the
$B$-symplectic form as we approach the vanishing set $z_1=0$ of
the bivector, express the differential form for $z_1\neq 0$ as
\begin{equation*}
z_1 e^{\frac{dz_1\wedge dz_2}{z_1}},
\end{equation*}
showing that as $z_1$ approaches zero,
$B+i\omega=\tfrac{1}{z_1}dz_1\wedge dz_2$ approaches infinity.
\end{example}

These deformations by holomorphic Poisson bivectors can be thought
of as non-commutative deformations of the complex manifold, in the
sense of quantization of Poisson structures.  The connection
between non-commutative geometry and these deformations of
generalized complex structure is explored by Kapustin
in~\cite{Kap2}, where the relations to topological string theory
are described as well.

\clearpage
\thispagestyle{empty}
\cleardoublepage
\chapter{Generalized K\"ahler
geometry}\label{GeneralizedKahler}

As we have seen, a generalized complex structure is an integrable
reduction of the structure group of $T\oplus T^*$ to $U(n,n)$.
This structure group may always be further reduced to its maximal
compact subgroup $U(n)\times U(n)$ by the choice of an appropriate
metric $G$ on $T\oplus T^*$. In this section we show that there is
an integrability condition which applies to such $U(n)\times U(n)$
structures, which generalizes the usual K\"ahler condition.  We
then show that this generalized K\"ahler geometry is equivalent to
a geometry first discovered by physicists (see~\cite{Rocek})
investigating supersymmetric nonlinear sigma-models.  Aspects of
this geometry, in particular the fact that it involves a
bi-Hermitian structure, were later studied (in the
four-dimensional case) by mathematicians
(see~\cite{Gauduchon},\cite{Pontecorvo},\cite{Kobak}).  The main
open problem in this field, as stated in~\cite{Gauduchon}, is to
determine whether or not there exist bi-Hermitian structures on
complex surfaces not admitting any anti-self-dual metric, for
example, $\CC P^2$.  Using generalized K\"ahler structures, we are
able to provide an affirmative solution to this problem. Finally,
we define twisted generalized K\"ahler structures, and describe an
interesting class of examples: the even-dimensional semi-simple
Lie groups.

\section{Definition}

Since the bundle $T\oplus T^*$ has a natural inner product
$\IP{,}$, it has structure group $O(2n,2n)$ in a natural way.  A
reduction from $O(2n,2n)$ to its maximal compact subgroup
$O(2n)\times O(2n)$ is equivalent to the choice of a
$2n$-dimensional subbundle $C_+$ which is positive definite with
respect to the inner product. Let $C_-$ be the (negative definite)
orthogonal complement to $C_+$. Note that the splitting $T\oplus
T^* = C_+\oplus C_-$ defines a positive definite metric on
$T\oplus T^*$ via
\[
G=\IP{,}\big|_{C_+}-\IP{,}\big|_{C_-}.
\]
Using the inner product to identify $T\oplus T^*$ with its dual,
the metric $G$ may be viewed as an automorphism of $T\oplus T^*$
which is symmetric, i.e. $G^*=G$, and which squares to the
identity, i.e. $G^2=1$. Note that $C_\pm$ are the $\pm 1$
eigenspaces of $G$. Hence we have the following:
\begin{prop}
A reduction to $O(2n)\times O(2n)$ is equivalent to specifying a
positive definite metric on $T\oplus T^*$ which is compatible with
the pre-existing inner product, i.e. $G^2=1$.
\end{prop}

Now suppose that we have a generalized almost complex structure
$\JJ$ defining a further reduction to $U(n,n)\subset O(n,n)$.  To
now reduce to $U(n)\times U(n)$ is equivalent to choosing a metric
$G$ as above, which commutes with the generalized complex
structure $\JJ$.  This is the same as choosing the space $C_+$ to
be stable under $\JJ$. Since $C_+$ is stable under $\JJ$ we see
that $\JJ$ is orthogonal with respect to $G$ and so we obtain a
Hermitian structure on $T\oplus T^*$ compatible with the
pre-existing inner product.

Note that since $G^2=1$ and $G\JJ =\JJ G$, the map $G\JJ$ squares
to $-1$, and since $G$ is symmetric while $\JJ$ is skew, $G\JJ$ is
also skew, and therefore defines a generalized almost complex
structure.  Hence we have the following:

\begin{prop}
A reduction to $U(n)\times U(n)$ is equivalent to the existence of
two generalized almost complex structures $\JJ_1,\JJ_2$ as well as
a positive definite metric $G$ satisfying $G^2=1$, which are
related by the following commuting diagram:
\begin{equation*}
\xymatrix{T\oplus T^*\ar[rr]^{G}&&T\oplus T^*\\
&T\oplus T^*\ar[ul]^{\JJ_1} \ar[ur]_{\JJ_2}&}
\end{equation*}
Note that the conditions on $G$ are equivalent to requiring that
$\JJ_1$ and $\JJ_2$ commute and that $-\JJ_1\JJ_2$ is positive
definite.
\end{prop}

We are now in a position to impose the integrability condition on
the $U(n)\times U(n)$ structure which defines generalized K\"ahler
structure. We simply require that both $\JJ_1$ and $\JJ_2$ are
integrable.
\begin{defn}
A generalized K\"ahler structure is a pair $(\JJ_1,\JJ_2)$ of
commuting generalized complex structures such that $G=-\JJ_1\JJ_2$
is a positive definite metric on $T\oplus T^*$.
\end{defn}
Our first example of a generalized K\"ahler structure justifies
the nomenclature.
\begin{example} Let $(g,J,\omega)$ be a usual K\"ahler structure
on a manifold, i.e. a Riemannian metric $g$, a complex structure
$J$, and a symplectic structure $\omega$ such that the following
diagram commutes.
\begin{equation*}
\xymatrix{T\ar[rr]^{g}&&T^*\\
&T\ar[ul]^{J} \ar[ur]_{\omega}&}
\end{equation*}
Then forming the generalized complex structures
\begin{equation*}
\JJ_J=\left(\begin{matrix}J&\\&-J^*\end{matrix}\right),\ \ \
\JJ_\omega=\left(\begin{matrix}&-\omega^{-1}\\\omega&\end{matrix}\right),
\end{equation*}
we see immediately that $\JJ_J,\JJ_\omega$ commute and
\begin{equation*}
G=-\JJ_J\JJ_\omega=\left(\begin{matrix}&g^{-1}\\g&\end{matrix}\right)
\end{equation*}
is a positive definite metric on $T\oplus T^*$.  Hence
$(\JJ_J,\JJ_\omega)$ defines a generalized K\"ahler structure.
\end{example}

\begin{example}
Given any generalized K\"ahler structure $(\JJ_1,\JJ_2)$, we may
transform it by a B-field, for $B$ any closed 2-form:
$(\JJ_1^B,\JJ_2^B)=(\BB\JJ_1\BB^{-1},\BB\JJ_2\BB^{-1})$ is also
generalized K\"ahler.  Applying such a transformation to the first
example $(\JJ_J,\JJ_\omega)$, we obtain the following generalized
complex structures
\begin{equation}\label{btrans}
\JJ_J^B=\left(\begin{matrix}J&\\BJ+J^*B&-J^*\end{matrix}\right),\
\ \
\JJ_\omega^B=\left(\begin{matrix}\omega^{-1}B&-\omega^{-1}\\\omega+B\omega^{-1}B&-B\omega^{-1}\end{matrix}\right).
\end{equation}
Similarly, the metric $G$ becomes
\begin{equation*}
G^B=\left(\begin{matrix}-g^{-1} B&g^{-1}\\g - Bg^{-1} B
&Bg^{-1}\end{matrix}\right),
\end{equation*}
showing that the metric $G$ in a generalized K\"ahler structure
need not be diagonal.  Note also that the restriction of $G^B$ to
the tangent bundle is the component $g-Bg^{-1}B$, which is indeed
a Riemannian metric for any two-form $B$.
\end{example}

\section{Torsion and the generalized K\"ahler metric}

The last example gives an indication of the general form of a
generalized K\"ahler metric.  Let $(\JJ_1,\JJ_2)$ be any
generalized K\"ahler structure and write the metric
$G=-\JJ_1\JJ_2$ as follows.
\begin{equation*}
G=\left(\begin{matrix} A & g^{-1} \\
\sigma & A^* \end{matrix}\right),
\end{equation*}
where $g,\sigma$ are genuine Riemannian metrics on the manifold
and $A$ is an endomorphism of $T$.  The condition $G^2=1$ implies
that $A$ is skew-symmetric with respect to both metrics $g$ and
$\sigma$, and that if we define the 2-form $b=-gA$, we can write
\begin{equation*}
G=\left(\begin{matrix}-g^{-1}b&g^{-1}\\g-bg^{-1}b&bg^{-1}\end{matrix}\right)
=\left(\begin{matrix}1&\\b&1\end{matrix}\right)
\left(\begin{matrix}&g^{-1}\\g&\end{matrix}\right)
\left(\begin{matrix}1&\\-b&1\end{matrix}\right).
\end{equation*}
We see from this argument that any generalized K\"ahler metric is
uniquely determined by a Riemannian metric $g$ together with a
2-form $b$.  We may also see $g$ and $b$ as follows: The
$+1$-eigenbundle $C_+$ of $G$ is positive definite in the natural
inner product, and since the tangent bundle is isotropic, $C_+$
can be expressed as the graph of a positive definite linear map
from $T$ to $T^*$.
\begin{prop}\label{graph}
$C_\pm$ is the graph of $b\pm g :T\lra T^*$
\end{prop}
\begin{proof}
Let $X+\xi\in C_+$, so that
\[
\left(\begin{matrix}-g^{-1}b&g^{-1}\\g-bg^{-1}b&bg^{-1}\end{matrix}\right)
\left(\begin{matrix}X\\
\xi\end{matrix}\right)=\left(\begin{matrix}X\\
\xi\end{matrix}\right).
\]
The equation in vector fields states that $-g^{-1}bX+g^{-1}\xi=X$,
i.e. $\xi=(b+g)X$, as required. The equation in 1-forms is
automatically satisfied. Similarly for $C_-$.
\end{proof}

It may seem from the discussion above that any generalized
K\"ahler metric is the B-field transform of a bare Riemannian
metric, but this is not the case, as the 2-form $b$ need not be
closed. We will present an example of this in a later section. The
derivative of $b$ actually plays an important role in generalized
K\"ahler geometry.
\begin{defn}
The \emph{torsion} of a generalized K\"ahler structure is the
3-form $h=db$.
\end{defn}

\section{Courant integrability}

We wish to describe the meaning of the generalized K\"ahler
condition in terms of subbundles of $(T\oplus T^*)\otimes\CC$.  As
usual, $\JJ_1$ and $\JJ_2$ engender a decomposition into $\pm i$
eigenbundles:
\begin{align*} (T\oplus
T^*)\otimes\CC&=L_1\oplus \overline{L_1}\\
&=L_2\oplus \overline{L_2}.
\end{align*}
Since $\JJ_1$ and $\JJ_2$ commute, $L_1$ must decompose into $\pm
i$ eigenbundles of $\JJ_2$, which we denote $L_1^{\pm}$. Then we
have the following decomposition into four isotropic subbundles:
\[
(T\oplus T^*)\otimes\CC = L_1^+\oplus L_1^-\oplus
\overline{L_1^+}\oplus\overline{L_1^-}.
\]
Note that $L_2=L_1^+\oplus \overline{L_1^-}$. Since $C_\pm$ is the
$\pm 1$ eigenbundle of $G=-\JJ_1\JJ_2$, we see that
\[
C_\pm\otimes\CC=L_1^\pm\oplus \overline{L_1^\pm},
\]
which, incidentally, proves that $\rk L_1^+=\rk L_1^-=n$, yielding
the following result.
\begin{prop}\label{evenod}
The generalized complex structures $\JJ_1,\JJ_2$ in a generalized
K\"ahler pair must have the same parity if $n$ is even and must
have opposite parity if $n$ is odd. For example, on a
4-dimensional manifold, the two generalized complex structures
comprising the $U(2)\times U(2)$ structure must have the same
parity.
\end{prop}
\begin{proof}
The rank $n$ bundle $L^+_1$ is the intersection of $L_1$ and
$L_2$, hence if $n$ is even, the maximal isotropics $L_1$ and
$L_2$ must have the same parity, and similarly for $n$ odd they
must have opposite parity.
\end{proof}
The positive-definiteness condition, in terms of $L_1^{\pm}$,
becomes
\[
\IP{p,\bar{p}}-\IP{q,\bar{q}}\geq 0,
\]
for all $p\in L_1^+$, $q\in L_1^-$, with equality if and only if
$p=q=0$. In other words, $\pm \IP{x,\bar{x}}>0$ for all nonzero
$x\in L_1^\pm$.

Furthermore, since $L_1^+=L_1\cap L_2$ and $L_1^-=L_1\cap
\overline{L_2}$, each of $L^\pm_1$ is closed under the Courant
bracket.  Using this information, we obtain a useful
characterization of a generalized K\"ahler structure.  First we
need a useful lemma which generalizes the decomposition
$d=\partial+\overline{\partial}$ for complex manifolds:

\begin{defn}
Let $L$ be a complex Lie algebroid with bracket $[,]$ and
differential $d_L:C^\infty(\wedge^kL^*)\lra
C^\infty(\wedge^{k+1}L^*)$. If $L=L^+\oplus L^-$ then we can
define $\wedge^{p,q}(L^*)=\wedge^p(L^+)^*\otimes\wedge^q(L^-)^*$,
as well as the operators
\begin{align*}
\partial^+_L&=\pi_{p+1,q}\circ d_L:C^\infty(\wedge^{p,q}L^*)\lra C^\infty
(\wedge^{p+1,q}L^*)\\
\partial^-_L&=\pi_{p,q+1}\circ d_L:C^\infty(\wedge^{p,q}L^*)\lra C^\infty
(\wedge^{p,q+1}L^*),
\end{align*}
where $\pi_{p,q}$ is the projection $\wedge^{p+q}L^*\lra
\wedge^{p,q}L^*$. If $L^\pm$ are closed under the Lie bracket,
then we have the equality
\begin{equation*}
d_L=\partial^+_L+\partial^-_L.
\end{equation*}
\end{defn}

\begin{prop}\label{integ}
A generalized K\"ahler structure on a real $2n$-dimensional
manifold is equivalent to the specification of two complex rank
$n$ subbundles $L_1^+, L_1^-$ of $(T\oplus T^*)\otimes\CC$
satisfying
\begin{itemize}
\item $L^\pm_1$ are isotropic, and
\item $L^+_1\bot L^-_1$ and $L^+_1\bot\overline{L^-_1}$, where $\bot$ indicates orthogonality with respect to the inner product, and
\item $L^\pm_1$ are $\pm$-definite, in the sense that
\begin{equation}\label{nondeg} \pm\IP{x,\bar{x}}>0\ \ \forall
x\in L_1^\pm,\ x\neq 0,
\end{equation}
\end{itemize}
with the integrability condition that both $L_1^\pm$ are closed
under the Courant bracket, and also that $L_1^+\oplus L_1^-$ is
closed under the Courant bracket.
\end{prop}
\begin{proof}
As we saw above, a generalized K\"ahler structure certainly
provides two subbundles with the required properties. We show the
converse. Given subbundles as above, condition~(\ref{nondeg})
implies that $L_1^\pm,\overline{L_1^\pm}$ are all mutually
transverse, and that $L_1 = L_1^+\oplus L_1^-$ and
$L_2=L_1^+\oplus \overline{L_1^{-}}$ are maximally isotropic,
defining two commuting almost generalized complex structures
$\JJ_1$ and $\JJ_2$ respectively. The condition that $L_1^+\oplus
L_1^-$ is integrable implies that $\JJ_1$ is integrable, and we
must show that $\JJ_2$ is integrable as well.  For this we use the
expression for the Courant bracket in the presence of the dual Lie
algebroid splitting $L_1\oplus\overline{L_1}=L_1\oplus L_1^*$ (we
identify $\overline{L_1}=L_1^*$ using the inner product):
\begin{align*}
[A+\alpha,B+\beta]&= [A,B] + \LL_\alpha B - \LL_\beta A
-d_{\overline{L}_1}(\IP{A,\beta}-\IP{B,\alpha})\\
&+[\alpha,\beta] + \LL_A\beta- \LL_B\alpha +
d_{L_1}(\IP{A,\beta}-\IP{B,\alpha}),
\end{align*}
where $A,B\in C^\infty(L_1)$ and $\alpha,\beta\in
C^\infty(\overline{L_1})$. Now, as a special case, take $A\in
C^\infty(L_1^+)$ and $\beta\in C^\infty(\overline{L_1^-})$,
setting $B=\alpha=0$. Then, using the fact that $L_1^+$ and
$\overline{L_1^-}$ are orthogonal,
\begin{align*}
[A,\beta]&=\LL_A\beta - \LL_\beta A\\
         &=i_A d_{L_1}\beta - i_\beta d_{L_1^*}A.
\end{align*}
where we recall that $\LL_A$ is defined by the Cartan formula
$\LL_A=i_A d_{L_1} + d_{L_1}i_A$.  Since $L_1^\pm$ are closed
under the bracket, we can use the decomposition of $d_{L_1}$
described in the previous lemma, obtaining
\begin{align*}
[A,\beta ] & = i_A (\partial^+_{L_1}+\partial^-_{L_1})\beta +
i_\beta (\partial^+_{L^*_1}+\partial^-_{L^*_1})A\\
&= i_A\partial^+_{L_1}\beta + i_\beta\partial^-_{L_1^*}A
\end{align*} showing that $[A,\beta]\in
C^\infty(L_1^+\oplus \overline{L_1^-})$, and therefore that
$L_2=L_1^+\oplus \overline{L_1^-}$ is also closed under the
Courant bracket, i.e. $\JJ_2$ is integrable.
\end{proof}

\section{Relation to bi-Hermitian geometry}\label{relat}

In this section we begin with a generalized K\"ahler structure
$(\JJ_1,\JJ_2)$ and `project' it to the tangent bundle, obtaining
more familiar types of structures, but with interesting
integrability conditions.

\subsection{The algebraic structure}
Since the bundle $C_+$ is positive definite while $T$ is null, the
projection $\pi:T\oplus T^*\lra T$ induces isomorphisms
\begin{equation*}
\xymatrix{\pi:C_\pm\ar[r]^{\ \cong}&T}.
\end{equation*}
This means that any algebraic structure existing on $C_\pm$ may be
transported to the tangent bundle.
\begin{itemize}
\item(Metric and 2-form)  As $C_\pm$ are definite subspaces of $T\oplus T^*$, and
as there are natural symmetric and skew-symmetric inner products
on $T\oplus T^*$, namely
\begin{align*}
\IP{X+\xi,Y+\eta}_+&=\frac{1}{2}(\xi(Y)+\eta(X))\\
\IP{X+\xi,Y+\eta}_-&=\frac{1}{2}(\xi(Y)-\eta(X)),
\end{align*}
we obtain natural Riemannian metrics and 2-forms on both of
$C_\pm$ by restriction.  If we transport these structures via
$\pi$ to the tangent bundle, we obtain precisely the metric and
2-form discussed in proposition~\ref{graph}, that is, projecting
from $C_\pm$ yields $b\pm g$, with $g$ a Riemannian metric on $T$
and $b$ a 2-form.

\item(Compatible almost complex structures) Since $C_\pm$ are stable
under both $\JJ_1$ and $\JJ_2$, they carry complex structures
which are each compatible with the metric $\IP{,}_+$, in the sense
that each is an orthogonal transformation.  Note that
$\JJ_1=\JJ_2$ on $C_+$ and $\JJ_1=-\JJ_2$ on $C_-$, so we only
need to project one of them, say $\JJ_1$.  By projection from
$C_\pm$, $\JJ_1$ induces two almost complex structures on $T$,
which we denote $J_\pm$, and these are compatible with the induced
Riemannian metric $g$, so we will call them Hermitian almost
complex structures.
\end{itemize}
In fact, as we now show, algebraically a $U(n)\times U(n)$
structure is equivalent to the specification of the quadruple
$(g,b,J_+,J_-)$, that is, a Riemannian metric $g$, a 2-form $b$,
and two Hermitian almost complex structures $J_\pm$. One could
call this structure an `almost bi-Hermitian structure with
b-field'.
\begin{defn}
Let $\omega_\pm$ be the 2-forms associated to the Hermitian almost
complex structures $J_\pm$, i.e.
\begin{equation*}
\omega_\pm=gJ_\pm.
\end{equation*}
\end{defn}

\begin{prop}
The generalized K\"ahler structure $(\JJ_1,\JJ_2)$ can be
reconstructed from the data $(g,b,J_+,J_-)$.
\end{prop}
\begin{proof}
The maps $b\pm g$ determine the metric $G$ by the formula
\begin{equation*}
G=\left(\begin{matrix}-g^{-1}b&g^{-1}\\g-bg^{-1}b&bg^{-1}\end{matrix}\right)
=\left(\begin{matrix}1&\\b&1\end{matrix}\right)
\left(\begin{matrix}&g^{-1}\\g&\end{matrix}\right)
\left(\begin{matrix}1&\\-b&1\end{matrix}\right),
\end{equation*}
and therefore they determine $C_\pm$. Using $\pi$ we can
reconstruct $\JJ_1$ by defining it to be the transport of $J_+$ on
$C_+$ and the transport of $J_-$ on $C_-$. To reconstruct $\JJ_2$
we use $J_+$ on $C_+$ and $-J_-$ on $C_-$. In formulae:
\begin{align*}
\JJ_1&=\pi\big|_{C_+}^{-1}J_+\pi P_++\pi\big|_{C_-}^{-1}J_-\pi P_-\\
\JJ_2&=\pi\big|_{C_+}^{-1}J_+\pi P_+-\pi\big|_{C_-}^{-1}J_-\pi
P_-,
\end{align*}
where $P_\pm$ are the projections from $T\oplus T^*$ to $C_\pm$,
namely $P_\pm=\frac{1}{2}(1\pm G)$.  Using these formulae we are
able to write $\JJ_1$ and $\JJ_2$ explicitly:
\begin{equation}\label{reconstruct}
\JJ_{1/2}=\frac{1}{2}\left(\begin{matrix}1&\\b&1\end{matrix}\right)
\left(\begin{matrix}J_+\pm J_- & -(\omega_+^{-1}\mp\omega_-^{-1}) \\
\omega_+\mp\omega_-&-(J^*_+\pm J^*_-)\end{matrix}\right)
\left(\begin{matrix}1&\\-b&1\end{matrix}\right).
\end{equation}
\end{proof}
\begin{remark}
From equation~(\ref{reconstruct}) we see that the degenerate case
where the almost complex structures $J_+$,$J_-$ are equal or
conjugate ($J_+ = \pm J_-$) corresponds to $(\JJ_1,\JJ_2)$ being
the B-field transform of a genuine K\"ahler structure.
\end{remark}

Before proceeding to the integrability conditions, we mention here
a fact about the relationship between the parity of the
generalized K\"ahler pair and the orientations of the almost
complex structures $J_\pm$.
\begin{remark}
As we have seen, in dimension $4k$ the generalized K\"ahler pair
must have the same parity, i.e. must have parity (even,even) or
(odd,odd). In the former case, the complex structures $J_\pm$
induce the same orientation on the manifold. In the latter case,
$J_\pm$ induce opposite orientations.

In dimension $4k+2$, the parities of the generalized complex
structures must be different, and this places no restriction on
the orientations of $J_\pm$; indeed changing the sign of either
$J_\pm$ reverses its orientation in this dimension.
\end{remark}

\subsection{The integrability condition: part I}

Now that we can express the generalized K\"ahler structure in
terms of an almost bi-Hermitian structure with b-field, we must
describe what the Courant integrability condition implies for the
quadruple $(g,b,J_\pm)$.

The first observation is that since the map $\pi:T\oplus T^*\lra
T$ is the Courant algebroid anchor, it satisfies the condition
\begin{equation*}
\pi[X+\xi,Y+\eta]=[\pi(X+\xi),\pi(Y+\eta)]=[X,Y].
\end{equation*}
Therefore we see that if a subbundle of $T\oplus T^*$ is closed
under the Courant bracket, then its image under $\pi$ will be
closed under the Lie bracket.  Immediately we have the following
result.

\begin{prop}
The complex structures $J_+,J_-$ coming from a generalized
K\"ahler structure are integrable, and therefore $(g,J_+,J_-)$ is
a bi-Hermitian structure.
\end{prop}
\begin{proof}
The integrability condition for the complex structure $J_+$ is
that the $+i$ eigenbundle in the complexified tangent bundle is
closed under the Lie bracket.  We denote this subbundle by
$T^{1,0}_+ < T\otimes\CC$. Since $J_+$ is obtained from the
complex structure $\JJ_1\big|_{C_+}$ via the isomorphism
$\pi|_{C_+}$, we see that the preimage of $T^{1,0}_+$ is the
bundle $L_1^+ < (T\oplus T^*)\otimes\CC$:
\begin{equation*}
\xymatrix{\pi:L_1^+\ar[r]^{\ \cong}&T^{1,0}_+}.
\end{equation*}
But we know that the bundle $L_1^+$ is closed under the Courant
bracket.  Hence $T^{1,0}_+$ is closed under the Lie bracket, as
required. Similarly for $J_-$.
\end{proof}

To obtain the full integrability conditions it will be useful to
explicitly describe the data entering into proposition~\ref{integ}
in terms of the bi-Hermitian structure with b-field.  To this end,
we observe from proposition~\ref{graph} that the bundle $L^+_1$
may be described as the graph of $g+b$ thought of as a map from
$T^{1,0}_+$ to $T^*\otimes\CC$:
\begin{equation}\label{Lplus}
\begin{split}
L^+_1&=\{X+(b+g)X\ |\ X\in C^\infty(T^{1,0}_+) \}\\
&=\{X+(b-i\omega_+)X\ |\ X\in C^\infty(T^{1,0}_+) \},
\end{split}
\end{equation}
where in the second line we use the fact that $g=-\omega_+J_+$. We
can use the same argument for $L^-_1$, obtaining the expression
\begin{equation}\label{Lminus}
\begin{split}
L^-_1&=\{X+(b-g)X\ |\ X\in C^\infty(T^{1,0}_-) \}\\
&=\{X+(b+i\omega_-)X\ |\ X\in C^\infty(T^{1,0}_-) \}.
\end{split}
\end{equation}
Now that we have expressed the bundles $L^\pm_1$ in terms of the
bi-Hermitian data, we must discover the meaning of the three
integrability conditions required by proposition~\ref{integ},
namely that $L_1^+$, $L_1^-$, and $L_1^+\oplus L_1^-$ are all
Courant integrable. The first two conditions are easily
understood, since each bundle $L_1^\pm$ has been expressed as the
graph of a complex two-form, and we may use the following result
about Courant integrability:
\begin{prop}\label{cx2form}
The subbundle of $(T\oplus T^*)\otimes\CC$ defined by
\begin{equation*}
F=\{X+cX\ |\ X\in C^\infty(E)\},
\end{equation*}
for some complex 2-form $c$ and subbundle $E$ of $T\otimes\CC$, is
Courant integrable if and only if $E$ is Lie integrable and the
form $c$ satisfies
\begin{equation*}
i_Yi_Xdc=0\ \ \forall\ X,Y\in C^\infty(E).
\end{equation*}
\end{prop}
Based on this, we obtain equivalent conditions to the Courant
integrability of $L_1^\pm$:
\begin{prop}\label{project}
Given an almost generalized K\"ahler structure $(\JJ_1,\JJ_2)$ and
its related almost bi-Hermitian structure $(g,b,J_\pm)$, the
bundles $L_1^\pm$ are Courant integrable if and only if
\begin{itemize}
\item
$T^{1,0}_\pm$ are Lie integrable, i.e. $J_\pm$ are integrable
complex structures, and
\item
The three 2-forms $\omega_\pm,b$ satisfy the conditions
\begin{equation}\label{condition}
d^c_-\omega_-=-d^c_+\omega_+=db,
\end{equation}
where $d^c_\pm = i(\bar\partial_\pm-\partial_\pm)$, and
$\partial_\pm$ is the $\partial$-operator for the complex
structure $J_\pm$. Condition~(\ref{condition}) may also be written
as follows:
\begin{equation}\label{noc}
db(X,Y,Z)=d\omega_+(J_+ X, J_+ Y, J_+ Z)=-d\omega_-(J_- X, J_- Y,
J_- Z)
\end{equation}
for all vector fields $X,Y,Z$.
\end{itemize}
\end{prop}
\begin{proof}
Proposition~\ref{cx2form} implies that $L_1^\pm$ is Courant
integrable if and only if $J_\pm$ is integrable and the following
condition holds:
\begin{equation*}
i_Yi_Xd(b\mp i\omega_\pm)=0,\ \ \forall X,Y\in
C^\infty(T^{1,0}_\pm).
\end{equation*}
Using the $(p,q)$ decomposition of forms determined by $J_\pm$, we
obtain the equivalent condition
\begin{equation*}
db\pm d^c_\pm\omega_\pm=0,
\end{equation*}
as required.  To obtain the final equation~(\ref{noc}), note that
since $\omega_\pm$ is of type $(1,1)$, $d\omega_\pm$ is of type
$(2,1)+(1,2)$. The $(2,1)$ forms have eigenvalue $i$ under
$\wedge^3 J_\pm^*$ while the $(1,2)$ forms have eigenvalue $-i$.
Hence when acting on a $(1,1)$ form,
$d^c_\pm=i(\bar\partial_\pm-\partial_\pm)$ is the same as
$-\wedge^3 J_\pm^*\circ d$, that is,
\begin{equation*}
d^c_\pm\omega_\pm(X,Y,Z)=-d\omega_\pm(J_\pm X, J_\pm Y,J_\pm Z),
\end{equation*}
as required.
\end{proof}
An immediate result of these relations is a constraint on the
\emph{torsion} $h=db$ of a generalized K\"ahler structure, namely:
\begin{corollary}
The torsion $h$ of a generalized K\"ahler structure is of type
$(2,1)+(1,2)$ with respect to both complex structures $J_\pm$;
equivalently, it satisfies the condition
\begin{equation}\label{twoone}
h(X,Y,Z)=h(X,JY,JZ)+h(JX,Y,JZ) + h(JX, JY, Z)
\end{equation}
for all vector fields $X,Y,Z$, and for both complex structures
$J=J_\pm$.  Note that this is equivalent to
\begin{equation}
h(J X, J Y, J Z)=h(JX, Y, Z)+h(X,JY, Z)+h(X, Y, JZ).
\end{equation}
\end{corollary}
Also, we obtain a description of the degenerate case when $h=0$:
\begin{corollary}
For a generalized K\"ahler structure with data $(g,b,J_\pm)$, the
following are equivalent:
\begin{itemize}
\item $h=db=0$
\item $(J_+,g)$ is K\"ahler
\item $(J_-,g)$ is K\"ahler.
\end{itemize}
\end{corollary}

Now that we have proven that the integrability of $L_1^\pm$
implies the existence of a bi-Hermitian structure with the
additional constraint~(\ref{condition}), we ask what additional
constraint arises from the final integrability condition in
Proposition~\ref{integ}, namely, that $L_1^+\oplus L_1^-$ must be
Courant integrable.  In fact we intend to show that no additional
constraint arises; indeed, this last condition is redundant.
However, to do this, we must develop some of the Riemannian
geometry on the generalized K\"ahler manifold.

\subsection{The Bismut connection versus the Levi-Civita
connection}\label{bis}

Unless a Hermitian manifold is K\"ahler, the Levi-Civita
connection is not a $U(n)$ connection.  However there are several
canonically defined $U(n)$ connections on a Hermitian manifold,
including the Chern connection (the unique $U(n)$ connection with
torsion $\tau$ such that $\tau(JX,JY)=-\tau(X,Y)$), and more
importantly for us, the Bismut connection, defined as follows.
\begin{defn}
The Bismut connection associated to a Hermitian structure is the
unique $U(n)$ connection with totally skew-symmetric torsion.
See~\cite{Bismut} for details.
\end{defn}
In the case of a K\"ahler structure, all these connections
coincide, as we now recall.
\begin{prop}\label{kahler}
Let $(g,J,\omega)$ be an almost Hermitian structure.  Then the
following are equivalent:
\begin{itemize}
\item $\nabla J=0$, where $\nabla$ is the Levi-Civita connection
\item $N_J=0$ and $d\omega=0$, i.e. $(g,J,\omega)$ is K\"ahler.
\end{itemize}
\end{prop}
It will be useful to precisely describe the difference between the
Levi-Civita and the Bismut connection for a general Hermitian
structure. To do this, we generalize Proposition~\ref{kahler}. We
need two lemmas:

\begin{lemma}
Let $(g,J,\omega)$ be an almost Hermitian structure.  Then if
$N_J$ is the Nijenhuis tensor of $J$ and $\nabla$ is the
Levi-Civita connection,
\begin{align}\label{nijenh}
N_J(X,Y)=(\nabla_{JX}J)Y-(\nabla_{JY}J)X + (\nabla_XJ)JY
-(\nabla_YJ)JX,
\end{align}
for any vector fields $X,Y$.
\end{lemma}
\begin{proof}
By definition, the Nijenhuis tensor is
\begin{align*}
N_J(X,Y)=[JX,JY]-J[JX,Y]-J[X,JY]-[X,Y].
\end{align*}
Using the fact that $\nabla$ is torsion-free, we obtain
\begin{align*}
N_J(X,Y)&=\nabla_{JX}JY-\nabla_{JY}JX-J(\nabla_{JX}Y-\nabla_YJX)-J(\nabla_{X}JY-\nabla_{JY}X)
-(\nabla_XY-\nabla_YX)\\
&=(\nabla_{JX}J)Y-(\nabla_{JY}J)X - J(\nabla_XJ)Y + J(\nabla_YJ)X.
\end{align*}
Differentiating $J^2=-1$ we obtain $(\nabla_XJ)J+J(\nabla_XJ)=0$,
which provides the result.
\end{proof}

\begin{lemma}
Let $(g,J,\omega)$ be an almost Hermitian structure. Then if
$\nabla$ is the Levi-Civita connection,
\begin{equation}\label{omeganabla}
d\omega(X,Y,Z)=g((\nabla_XJ)Y,Z)+g((\nabla_YJ)Z,X)+g((\nabla_ZJ)X,Y),
\end{equation}
for any vector fields $X,Y,Z$.
\end{lemma}
\begin{proof}
The exterior derivative is defined such that
\begin{align*}
d\omega(X,Y,Z)=\nabla_X\omega(Y,Z)-\omega([X,Y],Z)\ +\
\text{c.p.},
\end{align*}
where `c.p.' stands for cyclic permutations.  Using the facts that
$gJ=\omega$, that $\nabla$ is metric, and finally that $\nabla$ is
torsion free, we obtain
\begin{align*}
d\omega(X,Y,Z)&=g(\nabla_X(JY),Z) + g(JY,\nabla_XZ) -
\omega([X,Y],Z)\ + \ \text{c.p.}\\
&=g((\nabla_XJ)Y,Z) + g(J\nabla_XY,Z) + g(JY,\nabla_XZ) -
\omega([X,Y],Z)\ + \ \text{c.p.}\\
&=g((\nabla_XJ)Y,Z) + g(J\nabla_XY,Z) + g(JZ,\nabla_YX) -
\omega([X,Y],Z)\ + \ \text{c.p.}\\
&=g((\nabla_XJ)Y,Z) + g(J[X,Y],Z) -\omega([X,Y],Z)\ + \ \text{c.p.}\\
&=g((\nabla_XJ)Y,Z)\ +\ \text{c.p.},
\end{align*}
as required.
\end{proof}

\begin{prop}
Let $h$ be any 3-form and let $(g,J,\omega)$ be an almost
Hermitian structure. Consider the connection
\begin{equation}
\nabla^h=\nabla + \tfrac{1}{2}g^{-1}h.
\end{equation}
This is a metric connection with torsion $g^{-1}h$ and the
following are equivalent:
\begin{itemize}
\item $\nabla^h J = 0$ (or $\nabla^h\omega=0$).
\item $N_J=4g^{-1}h^{(3,0)+(0,3)}$ and $d\omega^{(2,1)+(1,2)}=
-ih^{(2,1)}+ih^{(1,2)}$.
\end{itemize}
\end{prop}
\begin{proof}
Substituting the definition of $\nabla^h$ into
equation~(\ref{nijenh}), we obtain the expression
\begin{equation*}
\begin{split}
g(N_J(X,Y),Z)&=g((\nabla^h_{JX}J)Y-(\nabla^h_{JY}J)X +
(\nabla^h_XJ)JY
-(\nabla^h_YJ)JX, Z) \\
&\quad - h(JX,JY,Z)- h(JX,Y,JZ)- h(X,JY,JZ)+ h(X,Y,Z)\\
&=g((\nabla^h_{JX}J)Y-(\nabla^h_{JY}J)X + (\nabla^h_XJ)JY
-(\nabla^h_YJ)JX, Z)+4h^{(3,0)+(0,3)}(X,Y,Z),
\end{split}
\end{equation*}
showing that if $\nabla^h J=0$, then $N_J=4g^{-1}h^{(3,0)+(0,3)}$.
Substituting the definition of $\nabla^h$ into
equation~(\ref{omeganabla}), we obtain the expression
\begin{equation*}
\begin{split}
d\omega(X,Y,Z)&=g((\nabla^h_XJ)Y,Z)+g((\nabla^h_YJ)Z,X)+g((\nabla^h_ZJ)X,Y)\\
&\quad -h(JX,Y,Z)-h(X,JY,Z)-h(X,Y,JZ),
\end{split}
\end{equation*}
showing that if $\nabla^hJ=0$, then
$d\omega^{(2,1)+(1,2)}=-ih^{(2,1)}+ih^{(1,2)}$, as required.

To show the converse, we combine equations~(\ref{nijenh}) and
(\ref{omeganabla}) to express $\nabla J$ in terms of $N_J$ and
$d\omega$:
\begin{equation*}
\begin{split}
g(N_J(X,Y),Z)&=g((\nabla_{JX}J)Y-(\nabla_{JY}J)X + (\nabla_XJ)JY
-(\nabla_YJ)JX,Z)\\
&=d\omega(JX,Y,Z)+d\omega(X,JY,Z) - 2g((\nabla_ZJ)X,JY).
\end{split}
\end{equation*}
Now, using the definition of $\nabla^h$:
\begin{equation*}
\begin{split}
2g((\nabla^h_ZJ)X,JY)&=2g(((\nabla_Z-\tfrac{1}{2}g^{-1}h_Z)J)X,JY)\\
&=d\omega(JX,Y,Z)+d\omega(X,JY,Z)-g(N_J(X,Y),Z) -h(JX,JY,Z)+
h(X,Y,Z),
\end{split}
\end{equation*}
which vanishes if we use the expressions for $N_J$ in terms of
$h^{(3,0)+(0,3)}$ and $d\omega^{(2,1)+(1,2)}$ in terms of
$h^{(2,1)+(1,2)}$, as well as the standard fact that for any
almost Hermitian structure, the skew part of the Nijenhuis tensor
is related to $d\omega^{(3,0)}$, i.e.
\begin{equation}
(i(d\omega)^{(3,0)}-i(d\omega)^{(0,3)})(X,Y,Z)=\tfrac{1}{4}(g(N_J(X,Y),Z)+\
\text{c.p.})
\end{equation}
This completes the proof.
\end{proof}

This proposition finally allows us to identify the Bismut
connection of any Hermitian structure:
\begin{corollary}
Let $(g,J,\omega)$ be a Hermitian structure.  Then the Bismut
connection is given by $\nabla^h=\nabla +\tfrac{1}{2}g^{-1}h$,
where $h=-d^c\omega=-i(\bar\partial-\partial)\omega$, and $\nabla$
is the Levi-Civita connection.
\end{corollary}

On a generalized K\"ahler manifold there are two Hermitian
structures $(g,J_+,J_-)$ and therefore two Bismut connections,
which we denote by $\nabla^\pm$. By our work on Hermitian
structures, we see that the torsion 3-form of the generalized
K\"ahler structure serves also as the torsion of the Bismut
connections:
\begin{theorem}
Let $(g,b,J_\pm)$ derive from an almost generalized K\"ahler
structure. Define $h=db$ and connections as follows
\begin{equation*}
\nabla^\pm=\nabla \pm \tfrac{1}{2}g^{-1}h.
\end{equation*}
These are metric connections with torsion
\begin{equation*}
\tau(\nabla^\pm) =\pm g^{-1}h.
\end{equation*}
The following are equivalent:
\begin{itemize}
\item $J_\pm$ are integrable and
$d^c_-\omega_-=-d^c_+\omega_+=h$,
\item $\nabla^\pm J_\pm=0$ and $h$ is of type $(2,1)+(1,2)$ with
respect to both $J_\pm$.
\end{itemize}
Finally, if these conditions hold, $\nabla^\pm$ become the Bismut
connections.
\end{theorem}
\begin{proof}
Let $J_\pm$ be integrable.  Then we have two Hermitian structures,
with Bismut connections $\nabla^\pm=\nabla -
\tfrac{1}{2}g^{-1}d^c\omega_\pm$ such that $\nabla^\pm J_\pm=0$.
But the condition $d^c_-\omega_-=-d^c_+\omega_+=h$ implies that
$\nabla^\pm=\nabla\pm g^{-1}h$, and also that $h$ is of type
$(2,1)+(1,2)$, as required. For the converse, $\nabla^\pm J_\pm=0$
implies that $N_{J_\pm}=\pm4g^{-1}h^{(3,0)+(0,3)}$ and
$d\omega_{\pm}^{(2,1)+(1,2)}= \mp ih^{(2,1)}\pm ih^{(1,2)}$.  If
$h$ is of type $(2,1)+(1,2)$, then $N_J=0$ and
$d^c_\mp\omega_{\mp}=\pm h$, as required.
\end{proof}
\begin{remark}
We now see that generalized K\"ahler geometry naturally induces a
bi-Hermitian structure whose integrability can be characterized by
two connections with totally skew-symmetric torsion which average
to the Levi-Civita connection:
\begin{equation*}
\frac{\nabla^++\nabla^-}{2}=\nabla.
\end{equation*}
This kind of geometry was first introduced in~\cite{Rocek} and
arises in physics as the natural geometry present on the target of
a N=(2,2) non-linear sigma model with torsion.
\end{remark}

\subsection{The integrability condition: part II}\label{integ2}

Having described the relevance of torsion connections to the
generalized K\"ahler geometry, we may now dispense with the
remaining integrability condition, that $L_1^+\oplus L_1^-$ be
Courant integrable, obtaining the full set of conditions for
generalized K\"ahler structures.
\begin{theorem}\label{equivbih}
Let $(\JJ_1,\JJ_2)$ be an almost generalized K\"ahler structure,
with associated data $(g,b,J_\pm)$, and let $h=db$. Any one of the
following conditions is equivalent to the integrability of the
generalized K\"ahler structure:
\begin{itemize}
\item[i)] The bundles $L^\pm_1$ are Courant integrable.
\item[ii)] $J_\pm$ are integrable and
$d^c_-\omega_-=-d^c_+\omega_+=h$.
\item[iii)] $\nabla^\pm J_\pm=0$ and $h$ is of type $(2,1)+(1,2)$ with respect
to both $J_\pm$ (the latter condition is implied by $J_\pm$
integrable).
\end{itemize}
\end{theorem}
\begin{proof}
We have shown that integrability implies condition $i)$, and we
have shown that $i),ii),iii)$ are equivalent.  All that remains is
to prove that one of $i),ii),iii)$ implies that $L^+_1\oplus
L^-_1$ is Courant integrable.

Let $X+\xi\in C^\infty(L^+_1)$ and $Y+\eta\in C^\infty(L^-_1)$. We
now show that $iii)$ implies that $[X+\xi,Y+\eta]\in
C^\infty(L^+_1\oplus L^-_1)$, which would complete the proof.
Recall that by expressions~(\ref{Lplus}) and (\ref{Lminus}),
$X+\xi=X+(b+g)X$ and $Y+\eta=Y+(b-g)Y$, and so
\begin{equation*}
\begin{split}
[X+\xi,Y+\eta]&=[X+gX+bX,Y-gY+bX]\\
&=[X+gX,Y-gY] + b([X,Y]) + i_Yi_Xdb\\
&=[X,Y]-\LL_X(gY)-\LL_Y(gX)+d(g(X,Y))+ b([X,Y]) + i_Yi_Xh.
\end{split}
\end{equation*} First we consider the tangent vector component.
Since $X+\xi\in C^\infty(L^+_1)$, by projection $X\in
C^\infty(T^{1,0}_+)$. Similarly, $Y\in C^\infty(T^{1,0}_-)$. Since
the Levi-Civita connection is torsion free, we have
\begin{align*}
[X,Y]&=\nabla_XY-\nabla_YX\\
&=(\nabla^-_X + \tfrac{1}{2}g^{-1}h_X)Y-(\nabla^+_Y -
\tfrac{1}{2}g^{-1}h_Y)X\\
&=\nabla^-_XY-\nabla^+_YX +\tfrac{1}{2}g^{-1}(h(X,Y)+h(Y,X))\\
&=\nabla^-_XY-\nabla^+_YX.
\end{align*}
If we assume $\nabla^\pm J_\pm=0$, then $\nabla^-_XY\in
C^\infty(T^{1,0}_-)$ and $\nabla^+_YX\in C^\infty(T^{1,0}_+)$.

Now, consider the 1-form component
\begin{equation*}
\begin{split}
i_Z(-\LL_X(g&Y)-\LL_Y(gX)+d(g(X,Y)))=-i_Z\LL_X(gY)-i_Z\LL_Y(gX)+
g(\nabla_ZX,Y)+g(X,\nabla_ZY)\\
&=-i_Z\LL_X(gY)-i_Z\LL_Y(gX)+
g([Z,X],Y)+g(X,[Z,Y])+g(\nabla_XZ,Y)+g(X,\nabla_YZ)\\
&=-\LL_Xi_Z(gY)-\LL_Yi_Z(gX)+i_Xd(i_Z(gY))+i_Yd(i_Z(gX))-g(Z,\nabla_XY+\nabla_YX)\\
&=-i_Zg(\nabla_XY+\nabla_YX).
\end{split}\end{equation*}
Therefore, we obtain the expression
\begin{equation*}
\begin{split}
[X+\xi,Y+\eta]&= [X,Y]-g(\nabla_XY+\nabla_YX)+ b([X,Y]) +i_Yi_Xh\\
&=\nabla^-_XY-\nabla^+_YX+ b(\nabla^-_XY-\nabla^+_YX)- g(\nabla^-_XY+\nabla^+_YX)\\
&=(-\nabla^+_YX + (b+g)(-\nabla^+_YX))\ +\ ((\nabla^-_XY +
(b-g)(\nabla^-_XY)),
\end{split}\end{equation*}
which is clearly in $C^\infty(L^+_1\oplus L^-_1)$, as required.
This completes the proof.
\end{proof}

\section{Examples of generalized K\"ahler structures}

We will now provide some examples of generalized K\"ahler
four-manifolds. Because of the equivalence between generalized
K\"ahler structures and bi-Hermitian structures, this places us
within the study of complex structures on Riemannian 4-manifolds.
In four dimensions there are two kinds of generalized K\"ahler
structure, as the pair of generalized complex structures may be of
parity (even,even) or (odd,odd), as we showed in
Proposition~\ref{evenod}. We will concentrate on the case
(even,even), which corresponds to the situation that both induced
complex structures $J_\pm$ in the bi-Hermitian structure have the
same orientation. A remarkable result of Salamon and
Pontecorvo~\cite{Pontecorvo} is that once a Riemannian 4-manifold
admits three distinct Hermitian complex structures with the same
orientation, it must admit infinitely many. Note that by
``distinct'' we mean in the sense that $\{J_i\}_{i=1}^n$ are
distinct when for each $i\neq  j$ there exists at least one point
$p\in M$ such that $J_i(p)\neq \pm J_j(p)$.
\begin{theorem}[\cite{Pontecorvo}]
Let $(M,g)$ be a Riemannian 4-manifold. Then it may admit either
$0$, $1$, $2$, or infinitely many distinct integrable Hermitian
complex structures with the same orientation.  Furthermore, if it
admits infinitely many, then it must be hyperhermitian and
therefore admit an $S^2$ family of orthogonal complex structures:
it must be either a flat torus, a $K3$ surface with its Ricci-flat
K\"ahler metric, or a hyperhermitian Hopf surface.
\end{theorem}

Of course any K\"ahler surface provides us with an example of a
generalized K\"ahler 4-manifold; in this case the generalized
K\"ahler pair is simply $(\JJ_J,\JJ_\omega)$, where $J,\omega$ are
the complex and symplectic structures. It is not difficult to see
that any Hyperk\"ahler structure provides us with another example:
\begin{example}[Hyperk\"ahler]\label{hkbih}
Let $(M,g,I,J,K)$ be a hyperk\"ahler structure. Then clearly
$(g,I,J)$ is a bi-Hermitian structure, and since
$d\omega_I=d\omega_J=0$, we see that $(g,I,J)$ defines a
generalized K\"ahler structure with $b=0$.  From
formula~(\ref{reconstruct}), we reconstruct the generalized
complex structures:
\begin{equation}\label{HKbih}
\JJ_{1/2}=\frac{1}{2}
\left(\begin{matrix}I\pm J & -(\omega_I^{-1}\mp\omega_J^{-1}) \\
\omega_I\mp\omega_J&-(I^*\pm J^*)\end{matrix}\right).
\end{equation}
From the expression~(\ref{btrans}) for B-field transforms, we
notice that~(\ref{HKbih}) describes two generalized complex
structures, each a $B$-field transform of a symplectic structure:
\begin{align*}
\JJ_{1}&=\left(\begin{matrix}1&\\\omega_K&1\end{matrix}\right)
\left(\begin{matrix} 0& -\tfrac{1}{2}(\omega_I^{-1}-\omega_J^{-1}) \\
\omega_I-\omega_J&0\end{matrix}\right)
\left(\begin{matrix}1&\\-\omega_K&1\end{matrix}\right)\\
\JJ_{2}&=\left(\begin{matrix}1&\\-\omega_K&1\end{matrix}\right)
\left(\begin{matrix} 0& -\tfrac{1}{2}(\omega_I^{-1}+\omega_J^{-1}) \\
\omega_I+\omega_J&0\end{matrix}\right)
\left(\begin{matrix}1&\\\omega_K&1\end{matrix}\right).
\end{align*}
In other words, the generalized K\"ahler structure is made up of
two $B$-symplectic structures defined by the differential forms
\begin{equation*}
\varphi_1=e^{B+i\omega_1}\ \ \ \varphi_2=e^{-B+i\omega_2},
\end{equation*}
where $B=\omega_K$, $\omega_1=\omega_I-\omega_J$, and
$\omega_2=\omega_I+\omega_J$.
\end{example}

The bi-Hermitian structure $(g,I,J)$ obtained from a hyperk\"ahler
structure is an example of a \emph{strongly} bi-Hermitian
structure, i.e. a bi-Hermitian structure $(g,J_+, J_-)$ such that
$J_+$ is nowhere equal to $\pm J_-$.  From expression
(\ref{reconstruct}), it is clear that in 4 dimensions, strongly
bi-Hermitian structures with equal orientation correspond exactly
to generalized K\"ahler structures where both generalized complex
structures are $B$-symplectic. Without loss of generality, the
$B$-field can be chosen to be $\pm B$ for the pair $\JJ_{1/2}$:

\begin{example}[strongly bi-Hermitian] A strongly bi-Hermitian
structure in 4 dimensions with equal orientation is equivalent to
a generalized K\"ahler structure $(\JJ_1,\JJ_2)$ such that both
$\JJ_1$ and $\JJ_2$ are $B$-symplectic, i.e. given by differential
forms
\begin{equation*}
\varphi_1=e^{B+i\omega_1},\ \ \ \varphi_2=e^{-B+i\omega_2}.
\end{equation*}
The algebraic conditions on these forms are as follows: first the
requirements $L_1\cap L_2\neq 0$ and $L_1\cap\overline{L_2}\neq 0$
imply that
\begin{equation}\label{speq}
\IPS{(\varphi_1,\varphi_2)}=\IPS{(\varphi_1,\overline{\varphi_2})}=0.
\end{equation}
In 4 dimensions there is only one additional constraint:
$-\JJ_1\JJ_2$ must be positive definite.  This is equivalent to
the requirement that $\omega_1^2, \omega_2^2$ have the same sign.
Equation~(\ref{speq}) is satisfied if and only if
$\IPS{(e^{B+i\omega_1},e^{-B\pm i\omega_2})}=0$, i.e. if and only
if $(2B+i(\omega_1\mp\omega_2))^2=0$.  Therefore we conclude that
a strongly bi-Hermitian structure in 4 dimensions with equal
orientation is equivalent to the specification of three closed
2-forms $B,\omega_1,\omega_2$ such that
\begin{equation}\label{strongbi}
\begin{split}
B\omega_1=B\omega_2&=\omega_1\omega_2=\omega_1^2+\omega_2^2-4B^2=0,\\
\omega_1^2&=\lambda\omega_2^2,\ \ \lambda>0.
\end{split}\end{equation}
\end{example}
This result appeared, of course without the generalized complex
interpretation, in~\cite{Gauduchon}, and it was used to show, as
we now do, that one can deform the hyperk\"ahler example mentioned
above in such a way to produce a genuine bi-Hermitian structure,
i.e. a metric which admits exactly two distinct Hermitian complex
structures:
\begin{example}[Non-hyperhermitian bi-Hermitian structure]
Let $B=\omega_K$, $\omega_1=\omega_I-\omega_J$, and
$\omega_2=\omega_I+\omega_J$ be the forms defining the
4-dimensional bi-Hermitian structure from Example~\ref{hkbih}.
Then let $F_t$ be a 1-parameter family of diffeomorphisms
generated by a $\omega_K$-Hamiltonian vector field.  For $t$
sufficiently small the forms
\begin{equation*}
B=\omega_K,\ \ \ \omega_1=\omega_I - F^*_t\omega_J,\ \ \
\omega_2=\omega_I+F^*_t\omega_J
\end{equation*}
satisfy the equations~(\ref{strongbi}), and hence define a
strongly bi-Hermitian structure.  In~\cite{Gauduchon} it is shown
that the Hamiltonian vector field can be chosen so that the
resulting deformed metric is not anti-self-dual, and hence does
not admit more than two orthogonal complex structures.  The idea
of deforming a Hyperk\"ahler structure in this way is due to D.
Joyce.
\end{example}

The question remains as to whether one can find bi-Hermitian
structures on manifolds which admit no hyperhermitian structure;
this problem has up to now remained unsolved.  The authors
of~\cite{Gauduchon} were unable to rule this out, but they were
able to obtain much information about examples, should they exist.
They reasoned that such a structure would not be strongly
bi-Hermitian; indeed $J_+$ and $\pm J_-$ would coincide along a
subvariety which would be an anticanonical divisor for both
$J_\pm$.  They went on to show that if the first Betti number is
even, then the 4-manifold must be either $\CC P^2$ or a minimal
ruled surface of genus $\leq 1$, or obtained from these by blowing
up points along an anti-canonical divisor.  In what follows we
show that there is in fact a bi-Hermitian structure on $\CC P^2$,
answering the question posed in~\cite{Gauduchon}.

\subsection{Deformation of generalized K\"ahler structures}

Let $(\JJ_1,\JJ_2)$ be a generalized K\"ahler structure.  In this
section we study deformations of this structure where $\JJ_1$
varies and $\JJ_2$ is kept fixed.  As we saw in
section~\ref{defm}, deformations of $\JJ_1$ as a generalized
complex structure are given by small sections $\eps$ of $\wedge^2
L_1$ which satisfy
\begin{equation*}
d_L\epsilon + \tfrac{1}{2}[\epsilon,\epsilon]=0.
\end{equation*}
In the presence of $\JJ_2$ in a generalized K\"ahler structure,
$L_1$ decomposes as $L_1=L_1^+\oplus L_1^-$.  $L_2$ can then be
written as $L_2=L_1^+\oplus \overline{L_1^-}$.  It is easy to see,
therefore, that small deformations of $\JJ_1$ keeping $\JJ_2$
fixed correspond to graphs of maps taking $L_1^+$ to
$\overline{L_1^-}$, i.e. small sections of the bundle
\begin{equation*}
\overline{L_1^+}\otimes \overline{L_1^-}<\wedge^2(L_1^+\oplus
L_1^-)^*,
\end{equation*}
satisfying the same condition as above.  As we proved in
section~\ref{defm}, assuming the obstruction vanishes, the locally
complete family of deformations of $\JJ_1$ has tangent space given
by $H^2_L(M)$.  Therefore if the infinite-dimensional subspace
$C^\infty(\overline{L_1^+}\otimes \overline{L_1^-})$ intersects
the space of $d_L$-closed sections of $\wedge^2 L^*$ in a subspace
with nontrivial image in $H^2_L(M)$, we can immediately conclude
that there exist deformations of the generalized K\"ahler
structure keeping $\JJ_2$ fixed.

As an example of this, we begin with a usual K\"ahler structure
and attempt to deform it in the generalized sense.  Let
$\JJ_J,\JJ_\omega$ be a usual K\"ahler structure.  Then
\begin{equation*}
L_1^+=\{X-i\omega X\ :\ X\in T_{1,0}\},
\end{equation*}
whereas
\begin{equation*}
L_1^-=\{Y+i\omega Y\ :\ Y\in T_{1,0}\}.
\end{equation*}
As a consequence, the deformation $\epsilon$ which we seek is a
linear combination of elements of the form
\begin{equation*}
\epsilon = (X+i\omega X)\wedge(Y-i\omega Y),\ \ X,Y\in T_{1,0}.
\end{equation*}
In particular, $\overline{L_1^+}\otimes \overline{L_1^-}$ has
sections of the form
\begin{equation*}
\tfrac{1}{2}((X+i\omega X)\wedge(Y-i\omega Y)-(Y+i\omega
Y)\wedge(X-i\omega X))=X\wedge Y + \omega X\wedge\omega Y,
\end{equation*}
i.e. sections of the form $\eps=\beta+\omega^{-1}\beta\omega$ for
$\beta$ a $(2,0)$-bivector. Recall that in this situation the
differential operator $d_L$ is just $\overline{\partial}$, and so
the equation $d_L\eps=0$ is satisfied if and only if $\beta$ is
holomorphic.  So we see that if the K\"ahler manifold admits a
holomorphic bivector, there must be a solution to the deformation
equation within the subspace of sections
$\beta+\omega^{-1}\beta\omega$.  Such deformations are clearly not
trivial, since the bivector component $\beta$ will reduce the type
of the complex structure wherever it is nonvanishing.  Hence we
conclude:
\begin{prop}
Let $(M,J,\omega)$ be a K\"ahler manifold, and assume that the
obstruction to generalized deformations of $J$ vanishes.  If the
K\"ahler manifold admits a holomorphic bivector $\beta$, then it
has a deformation as a generalized K\"ahler structure such that
the complex structure is deformed while the symplectic structure
remains unchanged.
\end{prop}

\begin{example}[Existence of bi-Hermitian structure on $\CC P^2$]
As we saw in Example~\ref{defmcp2}, the complex structure of $\CC
P^2$ has vanishing obstruction map. Furthermore $\CC P^2$ admits
holomorphic bivectors as $\wedge^2 T_{1,0}=\mathcal{O}(3)$.  Hence
we conclude from the deformation theorem that there exists a
nonzero generalized K\"ahler deformation keeping the symplectic
structure constant, within the class of sections
$\beta+\omega^{-1}\beta\omega$, where $\beta\in C^\infty(\wedge^2
T_{1,0})$.  Hence we obtain a non-K\"ahler generalized K\"ahler
structure on the real 4-manifold underlying $\CC P^2$, and by
Theorem~\ref{equivbih}, we see that this produces a genuine
bi-Hermitian structure, i.e. a metric with exactly two distinct
Hermitian complex structures with equal orientation.
\end{example}

\section{Twisted generalized K\"ahler
structures}\label{twistkahler}

In the list of bi-Hermitian structures studied by Apostolov,
Gauduchon, and Grantcharov, there are some which do not satisfy
the additional constraints $\mp d^c_\pm\omega_\pm=db$ of
generalized K\"ahler geometry.  For example, the hyperhermitian
Hopf surface:

\begin{example}[The Hopf surface: not generalized
K\"ahler]\label{HopfnotGK} Express $M=S^3\times S^1$ as the
quotient $(\CC^2-\{0\})/\{x\mapsto 2x\}$.  As a complex surface
this is known as the Hopf surface.  The hyperk\"ahler structure
$(\tilde g,I,J,K)$ on $\CC^2$ descends to a hyperhermitian
structure on the Hopf surface once we rescale the metric $\tilde
g\mapsto g=r^{-2}\tilde g$, where $r(x)=|x|$, the Euclidean length
in $\CC^2$. Then the associated 2-form $\omega_I$ is not closed
(The Hopf surface has odd first Betti number and hence cannot be
K\"ahler). The 3-form $d^c_I\omega_I$ is closed, however, and is
proportional to $i_{\del_r}(\omega_I\wedge\omega_I)$, the
pull-back of the volume form on $S^3$. Hence we obtain that
\begin{equation}\label{integral}
\int_{S^3}d^c_I\omega_I\neq 0.
\end{equation}
Therefore we see that $d^c_I\omega_I$ is not exact. Hence the
standard hyperhermitian structure on the Hopf surface is not
generalized K\"ahler.

We wish to show that no Hermitian metric could make the Hopf
surface generalized K\"ahler.  Suppose that $\omega$ is any other
positive $(1,1)$-form with respect to the standard complex
structure $I$. If $I$ is to be one complex structure in a
generalized K\"ahler pair, then $d^c\omega=db$, implying that
$dd^c\omega=0$.  Hence $\del\omega$ defines a cohomology class in
the Dolbeault group $H^1(\Omega^2)$.  This class must be nonzero
for the following reason: if $\del\omega=\delbar\sigma$ for some
$(2,0)$-form $\sigma$, then let $T$ be an elliptic curve in the
Hopf surface (which must be homologically trivial, bounding the
3-manifold $D$).  We have
\begin{equation*}
\begin{split}
0\neq\int_T\omega &= \int_D
d\omega=\int_D\del\omega+\delbar\omega\\
&=\int_D\delbar\sigma+\del\overline{\sigma}=\int_Dd(\sigma+\overline{\sigma})\\
&=\int_T\sigma+\overline{\sigma}=0,
\end{split}
\end{equation*}
a contradiction.  Furthermore since $H^1(\Omega^2)\cong \CC$ for
the Hopf surface (see~\cite{Barth}), $[\del\omega]$ must be a
nonzero multiple of $[\del\omega_I]$, i.e.
\begin{equation*}
\del\omega=k\del\omega_I + \delbar\tau,
\end{equation*}
for some $(2,0)$-form $\tau$ and nonzero $k$.  Summing with the
conjugate, we obtain
\begin{equation*}
d\omega=k\del\omega_I +\overline{k\del\omega_I}+
d(\tau+\overline{\tau}).
\end{equation*}
We have shown (\ref{integral}) that $\int_{S^3} \del\omega_I=c\neq
0$, and so by integrating the last equation over $S^3$ we see that
$kc$ is imaginary.  But recall that $d^c\omega=db$, i.e.
\begin{equation*}
db=i(\overline{k\del\omega_I} - k\del\omega_I + d(\overline{\tau}
- \tau)).
\end{equation*}
Therefore, integrating the last equation over $S^3$, we obtain
that $kc$ is real, a contradiction. We conclude that the Hopf
surface may never enter into a generalized K\"ahler bi-Hermitian
pair.
\end{example}

We now define a natural extension of generalized K\"ahler
geometry, which allows non-exact torsion $h=db$, and which will
accomodate the Hopf surface:
\begin{defn}[Twisted generalized K\"ahler]
An almost generalized K\"ahler structure $(\JJ_1,\JJ_2)$ is said
to be twisted generalized K\"ahler with respect to the closed
3-form $H$ if both generalized complex structures $\JJ_1,\JJ_2$
are twisted generalized complex structures, as defined in
section~\ref{twistcx}.
\end{defn}
In sections~\ref{bis} and~\ref{integ2} we proved that the
generalized K\"ahler integrability condition is equivalent to
certain conditions on the bi-Hermitian data $(g,b,J_\pm)$ it
defines.  Using the same methods, together with our results on the
twisted Courant bracket, it is easy to prove the following
generalization of that result.

\begin{theorem}
Let $(\JJ_1,\JJ_2)$ be an almost generalized K\"ahler structure,
and let $(g,b,J_\pm)$ be the associated almost bi-Hermitian data.
Also, let $H$ be a closed 3-form. Then the following are
equivalent:
\begin{itemize}
\item $(\JJ_1,\JJ_2)$ is an $H$-twisted generalized K\"ahler
structure.
\item $J_\pm$ are integrable complex structures and
\begin{equation*}
d^c_-\omega_-=-d^c_+\omega_+=H+db.
\end{equation*}
\item $\nabla^\pm J_\pm=0$, where
\begin{equation*}
\nabla^\pm=\nabla \pm \tfrac{1}{2}g^{-1}(H+db)
\end{equation*}
and $H+db$ is of type $(2,1)+(1,2)$ with respect to both $J_\pm$
(this last condition is implied by $J_\pm$ integrable).
\end{itemize}
\end{theorem}

Given this result, we revisit the Hopf surface to discover our
first example of a twisted generalized K\"ahler structure:
\begin{example}[The Hopf surface: twisted generalized K\"ahler]
Let $I$ be the standard complex structure on the Hopf surface
defined by $(\CC^2-\{0\})/\{x\mapsto 2x\}$.  We use the notation
of example~\ref{HopfnotGK}.  We noted there that $H=d^c_I\omega_I$
was a closed, non-exact 3-form.  Now let $\check{I}$ be the
complex structure on the Hopf surface induced by the complex
structure on $\CC^2$ defined by holomorphic coordinates
$(z_1,\overline{z_2})$; note that $\check{I},I$ have opposite
orientations, but are both Hermitian with respect to $g$. Also it
is clear that
$-d^c_{\check{I}}\omega_{\check{I}}=d^c_I\omega_I=H$, showing that
\begin{equation*}
(g, I, \check{I}, H)
\end{equation*}
defines an $H$-twisted generalized K\"ahler structure of type
(odd,odd) and with $b=0$.  Note that the complex structures
$I,\check{I}$ commute.  Also, it is remarkable that this geometry
was discovered by Ro\v{c}ek, Schoutens, and
Sevrin~\cite{RocekSchoutensSevrin} in the context of a
supersymmetric $SU(2)\times U(1)$ Wess-Zumino-Witten model.
\end{example}

Finally, we provide a striking family of examples of twisted
generalized K\"ahler manifolds: the even-dimensional semi-simple
real Lie groups.
\begin{example}[The even-dimensional compact semi-simple Lie groups]
It has been known since the work of Samelson and
Wang~\cite{Samelson}, \cite{Wang} that any compact
even-dimensional Lie group admits left- and right-invariant
complex structures $J_L,J_R$, and that if the group is
semi-simple, these can be chosen to be Hermitian with respect to
the bi-invariant metric induced from the Killing form $\IP{,}$.
The idea, then, is to use $(\IP{,},J_L,J_R)$ as a bi-Hermitian
structure with $b=0$ and to show that it is integrable as an
$H$-twisted generalized K\"ahler structure with respect to
$H(X,Y,Z)=\IP{[X,Y],Z}$, the bi-invariant Cartan 3-form. To see
that this works, let us compute $d^c_{J_L}\omega_{J_L}$:
\begin{equation*}
\begin{split}
A=d^c_{J_L}\omega_{J_L}(X,Y,Z)&=d\omega_{J_L}(J_LX,J_LY,J_LZ)\\
&=-\omega_{J_L}([J_LX,J_LY],J_LZ)+c.p.\\
&=-\IP{[J_LX,J_LY],Z}+c.p.\\
&=-\IP{J_L[J_LX,Y]+J_L[X,J_LY]+[X,Y],Z}+ c.p.\\
&=(2\IP{[J_LX,J_LY],Z}+c.p.) - 3H(X,Y,Z)\\
&=-2A - 3H(X,Y,Z),
\end{split}
\end{equation*}
Proving that $d^c_{J_L}\omega_{J_L}=-H$.  Since the right Lie
algebra is anti-isomorphic to the left, the same calculation with
$J_R$ yields $d^c_{J_R}\omega_{J_R}=H$, and finally we have
\begin{equation*}
-d^c_{J_L}\omega_{J_L}=d^c_{J_R}\omega_{J_R}=H,
\end{equation*}
i.e. $(\IP{,},J_L,J_R)$ forms an $H$-twisted generalized K\"ahler
structure.
\end{example}

\section{Generalized Calabi-Yau metrics}

While we will not explore this geometry in detail in this thesis,
we wish to define \emph{generalized Calabi-Yau metric} geometry,
which is a further reduction beyond generalized K\"ahler, and
which is the analog of a Calabi-Yau manifold.  This geometry has
been studied in the 4-dimensional case on $K3$ by
Huybrechts~\cite{Huybrechts}, and turns out to be precisely the
geometrical structure parametrized by the moduli space of
$N=(2,2)$ superconformal field theories with $K3$ target, as
studied by Aspinwall and Morrison~\cite{Aspinwall}.

\begin{defn}[Generalized Calabi-Yau metric] A generalized
Calabi-Yau metric geometry is defined by a generalized K\"ahler
structure $(\JJ_1,\JJ_2)$ where each generalized complex structure
has holomorphically trivial canonical bundle, i.e. their canonical
bundles have nowhere-vanishing closed sections $\rho_1,\rho_2\in
C^\infty(\wedge^\bullet T^*\otimes\CC)$. Finally, we require that
the lengths of these sections are related by a constant, i.e.
\begin{equation}\label{CY}
\IPS{(\rho_1,\overline{\rho_1})}=c\IPS{(\rho_2,\overline{\rho_2})},
\end{equation}
where $c\in\RR$ is a constant.  By rescaling the differential
forms by a constant, $c$ may be chosen to be $\pm 1$.

Of course, generalized Calabi-Yau metrics may also be twisted in
the presence of a closed 3-form $H$, by requiring that
$(\JJ_1,\JJ_2)$ is an $H$-twisted generalized K\"ahler structure
defined by $d_H$-closed forms $\rho_1,\rho_2$
satisfying~(\ref{CY}).
\end{defn}

\begin{example}[Calabi-Yau manifolds]
A usual Calabi-Yau manifold is a K\"ahler manifold of complex
dimension $m$ with symplectic form $\omega$ and holomorphic volume
form $\Omega$, satisfying
\begin{equation*}
\omega^m=\frac{i^mm!}{2^m}\Omega\wedge\overline{\Omega}.
\end{equation*}
In terms of the differential forms $e^{i\omega},\Omega$ defining
the generalized complex structures, this condition is simply that
\begin{equation*}
\IPS{(e^{i\omega},e^{-i\omega})}=(-1)^{\frac{m(m-1)}{2}}
\IPS{(\Omega,\overline{\Omega})},
\end{equation*}
i.e. $c=(-1)^{\frac{m(m-1)}{2}}$. So we see that Calabi-Yau
manifolds provide the basic examples.
\end{example}

\clearpage
\thispagestyle{empty}
\cleardoublepage
\chapter{Generalized complex submanifolds}\label{subm}

In this section we introduce the notion of generalized complex
submanifold, which generalizes both the idea of complex
submanifold and that of Lagrangian submanifold.  In fact, even in
the case of a usual symplectic manifold, there are generalized
complex submanifolds besides the Lagrangian ones: these are the
so-called co-isotropic A-branes discovered recently by Kapustin
and Orlov~\cite{Kapustin} in the context of D-branes of
topological sigma-models.  We show that covariance with respect to
$B$-field transformations dictates that generalized complex
submanifolds are not simply submanifolds but carry extra data,
consisting of a complex line bundle with connection, or more
generally, a section of a gerbe with connection, which gives rise
to a 2-form $F$ on the submanifold.  Such objects have been
described by string theorists as D-branes; we show here that they
arise completely naturally in the context of the generalized
geometry of $T\oplus T^*$.  We are thankful to Anton Kapustin for
emphasizing the importance of gauge invariance in determining the
correct definition of structured submanifold.

We should mention at the outset that generalized complex
submanifolds do not necessarily inherit a generalized complex
structure of their own; indeed they may be odd-dimensional.  In
this respect they resemble Lagrangian submanifolds.  It is
possible to study the generalization of symplectic submanifolds,
as is done in~\cite{Ben-Bassat}, however we will not explore this
topic here.

\section{The generalized tangent bundle}

The tangent bundle of a submanifold $M\subset N$ is a natural
sub-bundle $T_M\leq T_N|_{M}$ of the restriction of the ambient
tangent bundle.  This sub-bundle has an associated annihilator
$\Ann T_M\leq T_N^*$, also known as the conormal bundle of the
submanifold.  Taking the sum of these natural sub-bundles, we
obtain a natural real maximal isotropic sub-bundle of $(T_N\oplus
T_N^*)|_{M}$, which we call the \emph{generalized tangent bundle}:
\begin{defn}
The sub-bundle
\begin{equation*}
\tau_M=T_M\oplus \Ann T_M\leq (T_N\oplus T_N^*)\big|_{M}
\end{equation*}
is called the generalized tangent bundle of the submanifold
$M\subset N$.
\end{defn}

If the ambient manifold has a (possibly $H$-twisted) generalized
complex structure $\JJ$, a natural condition on a submanifold $M$
would be that its generalized tangent bundle is stable under
$\JJ$.  Indeed this definition reduces to familiar conditions in
the complex and symplectic cases:
\begin{example}[Complex submanifold] Let $(N,\JJ_J)$ be a complex
manifold.  Then $\tau_M=T_M\oplus\Ann T_M$ is stable under
\begin{equation*}
\JJ_J=\left(\begin{matrix}J&\\&-J^*\end{matrix}\right)
\end{equation*}
if and only if $T_M$ is stable under $J$, which is equivalent to
the condition that $M$ be a complex submanifold.
\end{example}
\begin{example}[Lagrangian submanifold] Let $(N,\JJ_\omega)$ be a
symplectic manifold.  Then $\tau_M=T_M\oplus\Ann T_M$ is stable
under
\begin{equation*}
\JJ_\omega=\left(\begin{matrix}&-\omega^{-1}\\\omega
&\end{matrix}\right)
\end{equation*}
if and only if
\begin{itemize}
\item $\omega$ takes $T_M$ into $\Ann T_M$, i.e. $M$ is
an isotropic submanifold, and
\item $\omega^{-1}$ takes $\Ann T_M$ into $T_M$, i.e. $M$ is
coisotropic,
\end{itemize}
i.e. $M$ must be isotropic and co-isotropic, hence Lagrangian.
\end{example}

The problem with this possible definition of a generalized complex
submanifold is that it does not behave well with respect to
$B$-field transformations.  That is, if $\tau_M$ is stable under
$\JJ$, then applying a $B$-field will modify $\JJ$ but not
$\tau_M$, and so under this definition $M$ would cease to be a
generalized complex submanifold when applying a $B$-field, which
is supposed to be an underlying symmetry of the whole geometry.

The answer to this problem is to modify the definition of a
generalized tangent space so that it is acted upon naturally by
$B$-field transformations.

Let $(M,F)$ be a pair consisting of a submanifold $M\subset N$ and
a real 2-form $F\in\Omega^2(M)$ on $M$, and suppose that
$dF=H|_M$, where $H$ is the closed 3-form on $N$ defining the
twist.  Then we will call this a generalized submanifold of
$(N,H)$.  Using the language of gerbes (see section~\ref{gerbes}),
$(M,F)$ is really a submanifold on which the gerbe is
trivializable, together with a trivialization with connection
relative to the gerbe connection (whose curvature is $H$). While
this language is more precise, we simplify matters in the
following definition:
\begin{defn}[Generalized submanifold]
Let $(N,H)$ be a pair consisting of a manifold $N$ and a closed
3-form $H$.  Then the pair $(M,F)$ of a submanifold $M\subset N$
together with a 2-form $F\in\Omega^2(M)$ is said to be a
generalized submanifold of $(N,H)$ iff $dF=H|_{M}$.
\end{defn}
The advantage of the gerbe interpretation is that a special case
of a generalized submanifold (when $H|_M=0$ and $F$ is integral)
is a triple $(M,L,\nabla)$ consisting of a submanifold together
with a line bundle with unitary connection.  The simplified
definition above only sees the curvature $F^\nabla\in\Omega^2(M)$.

Now we define the generalized tangent space of the generalized
submanifold:
\begin{defn}[Generalized tangent bundle]
The generalized tangent bundle $\tau_M^F$ of the generalized
submanifold $(M,F)$ is
\begin{equation*}
\tau_M^F=\left\{X+\xi\in T_M\oplus T_N^*\big|_{M}\ :\
\xi\big|_{M}=i_XF \right\},
\end{equation*}
a real, maximal isotropic sub-bundle $\tau^F_M<(T_N\oplus
T_N^*)|_{M}$.
\end{defn}
It is clear that since the generalized tangent bundle sits in
$(T_N\oplus T_N^*)|_M$, it is acted upon naturally by $B$-field
transforms of the ambient space; as a result the action of
$B$-fields on generalized submanifolds is as follows:
\begin{equation*}
e^B(M,F)=(M,F+B).
\end{equation*}
Clearly this transformation does not interfere with the condition
$dF=H|_M$.

\section{Generalized complex submanifolds}

Finally we are able to define generalized complex submanifolds in
such a way that the property of being one is covariant under
$B$-field transformations.
\begin{defn}[Generalized complex submanifold]
Let $(N,\JJ,H)$ be a $H$-twisted generalized complex manifold,
where $H$ is a real closed $3$-form.  Then the generalized
submanifold $(M,F)\subset(N,\JJ,H)$ is said to be a generalized
complex submanifold when $\tau_M^F$ is stable under $\JJ$.
\end{defn}

We now determine what the generalized complex submanifolds are in
the standard cases of complex and symplectic geometry.

\begin{example}[Complex case]
Let $(M,F)\subset(N,\JJ_J)$ be a generalized submanifold of a
complex manifold. Then it is generalized complex if and only if
\begin{equation*}
\tau_M^F=\left\{X+\xi\in T_M\oplus T_N^*\big|_{M}\ :\
\xi\big|_{M}=i_XF \right\}
\end{equation*}
is stable under the action of
\begin{equation*}
\JJ_J=\left(\begin{matrix}J&\\&-J^*\end{matrix}\right),
\end{equation*}
which happens if and only if
\begin{itemize}
\item $T_M$ is closed under $J$, i.e. $M$ is a complex
submanifold, and
\item $J^*i_XF+i_{JX}F\in\Ann T_M\ \forall X\in T_M$, i.e. $F$ is
of type $(1,1)$ on $M$.
\end{itemize}
In the special case where $H|_{M}=0$ and $F$ is integral, then we
see that the generalized complex submanifold is a complex
submanifold together with a unitary holomorphic line bundle on it
(a line bundle with curvature of type $(1,1)$ has a unique
compatible holomorphic structure).
\end{example}

\begin{example}[Symplectic case]
Let $(M,F)\subset(N,\JJ_\omega)$ be a generalized submanifold of a
symplectic manifold. Then it is generalized complex if and only if
\begin{equation*}
\tau_M^F=\left\{X+\xi\in T_M\oplus T_N^*\big|_{M}\ :\
\xi\big|_{M}=i_XF \right\}
\end{equation*}
is stable under the action of
\begin{equation*}
\JJ_\omega=\left(\begin{matrix}&-\omega^{-1}\\\omega&\end{matrix}\right).
\end{equation*}
To aid with the calculation we choose a 2-form $B$ on $N$ such
that $B|_M=F$ (only locally).  Then $\tau_M^F$ is
$\JJ_\omega$-stable if and only if $\tau_M^0=T_M\oplus\Ann M$ is
stable under
\begin{equation*}
e^{-B}\JJ_\omega
e^{B}=\left(\begin{matrix}-\omega^{-1}B&-\omega^{-1}\\\omega+B\omega^{-1}B&B\omega^{-1}\end{matrix}\right),
\end{equation*}
which happens if and only if
\begin{itemize}
\item $\omega^{-1}(\xi)\in T_M\ \forall \xi\in\Ann T_M$, i.e.
$M$ is a coisotropic submanifold (note that $\Ann T_M$ always sits
inside $\tau_M^F$), and
\item $\omega^{-1}(i_XB)\in T_M\ \forall X\in T_M$, i.e.
$i_X F=0\ \forall X\in T_M^{\bot}$ ($T_M^\bot$ is the symplectic
orthogonal bundle), i.e. $F$ descends to $T_M/T_M^{\bot}$, and
\item $\omega+ B\omega^{-1}B$ sends $T_M$ into $\Ann T_M$, i.e.
$(\omega|_M)^{-1}F$ is an almost complex structure on
$T_M/T_M^\bot$.
\end{itemize}
Note that $F+i\omega|_M$ defines a nondegenerate form of type
$(2,0)$ on $T_M/T_M^\bot$, and so the complex dimension of
$T_M/T_M^\bot$ must be even.  Hence if $\dim N=2n$, then $\dim
M=n+2k$ for some positive integer $k$.  In the case where
$H|_M=0$, $F$ is closed and we obtain a holomorphic symplectic
structure transverse to the Lagrangian foliation of the
coisotropic submanifold.

In the case that $M$ is a Lagrangian submanifold, then the second
condition implies that $F=0$.  This means that Lagrangian
submanifolds can only be generalized complex submanifolds when
$H|_M=dF=0$, and in this case they always carry a flat line bundle
on them.

It is remarkable that in the symplectic case, the notion of
generalized complex submanifold coincides exactly with the
recently discovered ``coisotropic A-branes'' of Kapustin and Orlov
~\cite{Kapustin}.  This connection is explored in further detail
by Kapustin in~\cite{Kap2}.
\end{example}

Before we move on to discuss when generalized complex submanifolds
can be \emph{calibrated}, we provide one example of generalized
complex submanifold in a case which is neither complex nor
symplectic.
\begin{example}[Deformed $\CC P^2$]
In example~\ref{defmcp2}, we studied a generalized complex
structure on $\CC P^2$ obtained by deformation.  This structure is
$B$-symplectic outside a cubic, and along the cubic the
generalized complex structure is none other than the original
complex structure on $\CC P^2$.  Hence it is easy to see that
since the cubic began as a complex submanifold and is fixed by the
deformation, this means that it remains a generalized complex
submanifold of the deformed structure.
\end{example}

\section{Generalized calibrations}\label{calibrate}

In this section we wish to describe briefly an idea which can be
used to generalize the notion of special Lagrangian submanifold.
In the usual Calabi-Yau case, a Lagrangian submanifold is said to
be special Lagrangian when it is calibrated with respect to the
real part of the Calabi-Yau volume form.  That is,
$\text{Re}(\Omega)$ restricts to the Lagrangian to yield the
volume form induced by the Riemannian metric.

Having established the machinery of generalized complex
submanifolds and generalized Calabi-Yau metrics, the
generalization is clear:

\begin{defn}[Calibration in generalized Calabi-Yau manifolds]
Let $(\JJ_1,\JJ_2)$ define a generalized Calabi-Yau metric
structure, as in definition~\ref{CY}.  It is defined by two global
closed complex differential forms $\rho_1,\rho_2$.  By squaring
the spinors we obtain top-degree skew forms
\begin{equation*}
\Omega_1\in C^\infty(\det L_1),\ \ \Omega_2\in C^\infty(\det L_2).
\end{equation*}
Using the natural isomorphism $L_i=(T\oplus T^*,\JJ_i),\ i=1,2$,
we see that $\Omega_i$ are complex top-degree forms on $T\oplus
T^*$.

Let $(M,F)$ be a generalized complex submanifold with respect to
$\JJ_1$.  Then it is said to be calibrated with respect to $\JJ_2$
when the real part of $\Omega_2$, restricted to the generalized
tangent bundle of $M$, yields the natural volume induced by the
K\"ahler metric $G=-\JJ_1\JJ_2$ on $T\oplus T^*$, i.e.
\begin{equation*}
\text{Re}(\Omega_2)\big|_{\tau_M^F}=\det G\big|_{\tau_M^F}.
\end{equation*}

This condition can also be expressed directly in terms of
$\rho_2$, by requiring that
\begin{equation*}
\tau_M^F\cdot \text{Im}(\rho_2)=0,
\end{equation*}
i.e. the Clifford product of any generalized tangent vector with
$\text{Im}(\rho_2)$ vanishes.
\end{defn}
This definition carries over without change to the twisted case,
and specializes to the usual notion of SLAG in the presence of a
usual Calabi-Yau structure. Of course this definition engenders
many natural questions about moduli of such calibrated
submanifolds as well as about their intersection theory, among
others.  It will be the subject of future work.

\clearpage
\thispagestyle{empty}
\cleardoublepage

\chapter{Speculations on mirror symmetry}

In this highly speculative chapter, we outline some thoughts
concerning an approach to mirror symmetry through generalized
complex geometry.  It is clear that many aspects of mirror
symmetry manifest themselves naturally in the context of
generalized complex geometry, and so it is reasonable to expect
that the ``mirror relation'', if it indeed exists, should be
stated in the language of generalized complex manifolds.

Let $(M_A,G_A,\nabla_A)$ and $(M_B, G_B,\nabla_B)$ be manifolds
with gerbes and gerbe connections, for which the curvatures are
$H_A,H_B$: closed 3-forms, defining integral cohomology classes.
Then on $M_A\times M_B$ we can put a gerbe $\pi_A^*G_A^{-1}
\otimes \pi_B^*G_B$ (which we will denote $G_A^{-1}G_B$ for short)
with the tensor product connection $\nabla$. It has curvature
$-\pi_A^*H_A+\pi_B^*H_B$.  In general we needn't consider
$M_A\times M_B$ but perhaps another gerbed manifold fibering over
$M_A$ and $M_B$ in the appropriate way.

Suppose there is a generalized submanifold $(\MM,F)\subset
M_A\times M_B$ which submerges onto both $M_A$ and $M_B$ via
$\pi_A,\pi_B$:
\begin{equation*}
\xymatrix{ (\MM,F)\ar@{^{(}->}[r]&M_A\times M_B\ar[dl]_{\pi_A}\ar[dr]^{\pi_B}& \\
M_A & &M_B}
\end{equation*}
Recall that this means that the gerbe $G_A^{-1}G_B$ is
trivializable over $\MM$, i.e.
\begin{equation*}
[-\pi_A^*H_A+\pi_B^*H_B]\big|_{\MM}=0,
\end{equation*}
and that we have a choice of trivialization with connection over
$\MM$, giving rise to a 2-form $F\in\Omega^2(\MM)$ satisfying
\begin{equation*}
(-\pi_A^*H_A+\pi_B^*H_B)\big|_{\MM}=dF.
\end{equation*}
Now $(\MM,F)$ has generalized tangent bundle sitting in
$T_{M\times N}\oplus T^*_{M\times N}$:
\begin{equation*}
\tau_\MM=\{X+\xi\in T_\MM\oplus T^*_{M\times N}\ :\
\xi|_{T_{\MM}}=i_XF\},
\end{equation*}
Of course $M$ and $N$ have generalized tangent bundles in
$T_{M\times N}\oplus T^*_{M\times N}$ as well, namely $T_M\oplus
T^*_M$ and $T_N\oplus T^*_N$.  We require that $\tau_\MM$ be
transverse to both these.  This is an analog of a transversality
condition on the submanifold $\MM$.

Then we propose that $(\MM,F)$ defines a ``mirror relation''
between the manifolds $\MM_A$ and $\MM_B$, and there is a
transform-like mapping (of Dirac structures, differential forms,
and other objects, like generalized holomorphic bundles) defined
by
\begin{equation*}
\FF=(\pi_A)_*\circ e^F\circ (\pi_B\big|_{\MM})^*.
\end{equation*}

Now suppose we have generalized complex structures $L_A$, $L_B$,
and $\LL$ on the spaces $(M_A,G_A)$, $(M_B,G_B)$, and $(M_A\times
M_B, G_A^{-1}G_B)$, such that $(\pi_A)_*\LL=L_A$ and
$(\pi_B)_*\LL=L_B$.  And suppose further that $(\MM,F)$ is a
generalized complex submanifold of $(M_A\times
M_B,G_A^{-1}G_B,\nabla,\LL)$. Then we propose that $\MM$ defines a
mirror relation between the generalized complex structures on
$\MM_A$, $\MM_B$ by the same transform $\FF$. This transform can
be thought of as a generalization of the following ideas:
\begin{itemize}
\item T-duality, where $M_A$ and $M_B$ are $S^1$-bundles via $\varphi_A,\varphi_B$ over some
base $B$ and $\MM$ is the correspondence space
\begin{equation*}
\{(x,y)\in M_A\times M_B\ :\ \varphi_A(x)=\varphi_B(y)\}\subset
M_A\times M_B
\end{equation*}
on which is defined a natural $2$-form $F$ coming from the
T-duality construction.  This construction is described in the
physics literature and especially in the remarkable paper of
Bouwknegt, Evslin, and Mathai~\cite{Bouwknegt}.  Understanding
T-duality in the generalized complex framework is joint work
between the author and Gil Cavalcanti; this will appear
in~\cite{GilMarco2}.

\item The Fourier-Mukai transform, where both generalized complex
structures are complex.  Also bi-meromorphic mappings in algebraic
geometry.

\item B\"acklund transformations, showing that a PDE system on
$M_A$ is equivalent to a PDE system on $M_B$ through the use of a
correspondence space.

\item Poisson dual pairs and canonical transformations, where
$L_A$, $L_B$ are taken to be real Dirac structures.
\end{itemize}

In the case where both spaces $M_A$, $M_B$ are generalized
Calabi-Yau metric geometries, it may be possible to find a
submanifold $\MM\subset M_A\times M_B$ which satisfies the
properties of a mirror relation with respect to both pairs of
generalized complex structures, and furthermore is calibrated in
the sense of section~\ref{calibrate}. The hope then would be that
the known instances of mirror ensembles of Calabi-Yau 3-folds
could be phrased using this description.
It is with this wild speculation that we end the thesis.

\bibliographystyle{plain}
\bibliography{mgsingle}
\end{document}